\documentclass[12pt]{amsart}
\usepackage{amssymb,amscd}
\usepackage{pb-diagram}
\usepackage{verbatim}

\usepackage[dvipdfmx]{graphicx,color}  %これは、platex(ptex2pdf)用　これの[]を消してこのすぐしたの行作った　

\usepackage{graphicx}
\usepackage{mathrsfs}

\usepackage{here}   %zu wooku  \begin{figure}[H]

\usepackage[english]{babel}
\usepackage{amsfonts,amssymb,amsmath,amstext,amsbsy,amsopn,amscd,amsthm,graphicx,euscript,graphicx}
\usepackage{graphics}
\usepackage[all]{xy}

\newcommand{\tart}[2]{\begin{tabular}{c} #1 \\ #2\end{tabular}}
\newcommand{\incg}[2]{\includegraphics[height=#1em]{#2.eps}}%{pics/#2.eps}}
\newcommand{\targ}[3]{\tart{\incg{#1}{#2}}{#3}}

\newtheorem*{Whitney towers}{Theorem~\ref{Whitney towers}}
\newtheorem*{h-towers}{Theorems ~\ref{half} \& \ref{$(n)$-solvable}}

\newtheorem*{surgery curves}{Theorem~\ref{surgery curves}}
\newtheorem*{cg=0}{Theorem~\ref{vanish}}

\newtheorem{thm}{Theorem}[section]
\newtheorem{mth}[thm]{Main Theorem}
\newtheorem{pr}[thm]{Proposition}

\newtheorem{fact}[thm]{Fact}

\theoremstyle{definition}
\newtheorem{defn}[thm]{Definition}
\newtheorem{note}[thm]{Note}

\newtheorem{exa}[thm]{Example}

\numberwithin{equation}{section}
\numberwithin{figure}{section}
\numberwithin{table}{section}

\newcommand{\bb}{\bigbreak}
\newcommand{\bs}{\smallbreak}
\newcommand{\h}{\noindent}
\newcommand{\x}{\times}
\newcommand{\np}{\newpage}
\newcommand{\Z}{\mathbb{Z}}

\newcommand{\C}{\mathbb{C}}

\newcommand{\R}{\mathbb{R}}

\def\yen{{\setbox0=\hbox{Y}Y\kern-.97\wd0\vbox{hrule height.lex width.98%
\wd0\kern.33ex\hrule height.lex width.98\wd0\kern.45ex}}}%%

\def\np{\newpage}

\begin{document}
\pagestyle{plain}

\title{
Khovanov-Lipshitz-Sarkar homotopy type %\\
%and Steenrod square \\
 for links in thickened higher genus surfaces
}
\author{Louis H. Kauffman,
%Vassily Olegovich Manturov, \\
Igor Mikhailovich Nikonov,
and
Eiji Ogasa}

%\thanks{\hskip-4mm %E-mail: %kauffman@uic.edu\quad
%ogasa@mail1.meijigakuin.ac.jp \newline
%Keywords: Virtual knots, Tori in $S^4$ }
% {\bf PACS nos.} 11-25w, 11-25Uv.
%\newline MSC2000 57N10, 57N13, 57N15.
%\newline Differential Topology

\date{}

%significantの段落取ってなおす  書いた2021-7-21

%Virtual のKh　と　links in surfaceのvirtual はちゃう　と　強調  書いた2021-7-21

%ここではMantNikoのみやる　それの定義は復習する　書いた2021-7-21

%VirtualのKhはtrivialでもtrivial knotでないの有る  書いた2021-7-21

\begin{abstract}
We discuss links in thickened surfaces.
We define the Khovanov-Lipshitz-Sarkar stable homotopy type and
the Steenrod square for
the homotopical Khovanov homology of links in thickened surfaces with genus$>1$ .

A surface means a closed oriented surface unless otherwise stated.
Of course, a surface may or may not be the sphere.
A thickened surface means
a product manifold of a surface and the interval.
A link in a thickened surface (respectively, a 3-manifold) means
a submanifold of a thickened surface (respectively, a 3-manifold)
which is diffeomorphic to a disjoint collection of circles.
%We discuss links in thickened surfaces.

%We define the Khovanov-Lipshitz-Sarkar stable homotopy type  and the Steenrod square   for the homotopical Khovanov homology of links in thickened surfaces with genus$>1$.

Our Khovanov-Lipshitz-Sarkar stable homotopy type  and our Steenrod square
of links in thickened surfaces with genus$>1$
are stronger
than
the homotopical Khovanov homology
of links in thickened surfaces with genus$>1$.
%Each of the former two invariants, which we define in this paper,  is stronger than the latter one, which had known before this paper, as invariants of links in thickened surfaces.

It is the first meaningful Khovanov-Lipshitz-Sarkar stable homotopy type of links in
 3-manifolds other than the 3-sphere.

We point out that our theory has a different feature in the torus case.
\end{abstract}
\maketitle

\tableofcontents

\bb
\section{\bf Introduction}\label{intro}
%{\Large
\h
In this paper, we discuss links in thickened surfaces. We define the Khovanov-Lipshitz-Sarkar stable homotopy type and the Steenrod square for the homotopical Khovanov homology of links in thickened surfaces with genus$>1$.
\\

A surface means a closed oriented surface unless otherwise stated.
Of course, a surface may or may not be the sphere.
\\

A thickened surface means
a product manifold of a surface and the interval.
A link in a thickened surface means
a submanifold of a thickened surface (respectively, a 3-manifold)
which is diffeomorphic to a disjoint collection of circles.
%We discuss links in thickened surfaces.
If $\mathcal L$ is a link in a thickened surface,
then a link diagram $L$ which represents  $\mathcal L$  is in the surface.
\\

Our theory has a special behavior at genus one as explained in \S\ref{hiyoko}.
In this paper, a higher genus surface means a surface with genus greater than one unless otherwise stated.
 We will discuss the torus case in a sequel of this paper \cite{zoku}.
\\

In \cite{K},  Khovanov defined the Khovanov homology for  links in $S^3$,
 and proved that  its graded Euler characteristic is the Jones polynomial of the link.
In \cite{B}, Bar-Natan proved  that
the Khovanov homology  is stronger than  the Jones polynomial
as invaraints of  links in $S^3$.

In \cite{LSk}, Lipshitz and Sarkar defined
the Khovanov-Lipshitz-Sarkar stable homotopy type for  links in $S^3$,
and proved that
the cohomology group of
the Khovanov-Lipshitz-Sarkar stable homotopy type of any  link $L$ in $S^3$
is the Khovanov homology of $L$.

\bb
\h{\bf Note.}
Khovanov-Lipshitz-Sarkar stable homotopy type
is sometimes abbreviated to
Khovanov-Lipshitz-Sarkar homotopy type or
Khovanov homotopy type,
in this paper and in other papers,
when it is clear from the context.
\bb

In \cite{LSs}, Lipshitz and Sarkar found
a method to calculate
the second Steenrod square operator on the Khovanov homology for links in $S^3$.
In \cite{Seed}, Seed made a computer program of the above method, used it, and showed that
the second Steenrod square operator  and
the Khovanov-Lipshitz-Sarkar stable homotopy type
are stronger than the Khovanov homology
as invaraints of  links in $S^3$.
\\

In \cite{APS},
 Asaeda,  Przytycki, and Sikora
 extended Khovanov homology for links in $S^3$ to
 thickened surfaces.
In \cite{MN}, Manturov and Nikonov made
an alternative definition of that in \cite{APS},
and obtained a new result by using it.
 There, the homology is called the {\it homotopical Khovanov homology}.
 We review the definition  in \S\ref{KMNhomology}.

\bb
In this paper,
we construct Khovanov-Lipshitz-Sarkar stable homotopy type
for the homotopical Khovanov homology for links in thickened surfaces with genus$>1$.
It is our main result, Main Theorem \ref{main},  below.
It is the first meaningful example
of Khovanov-Lipshitz-Sarkar stable homotopy type
for links in other 3-manifolds than the 3-sphere.

\begin{mth}\label{main}
$(1)$ We define Khovanov-Lipshitz-Sarkar stable homotopy type for
the homotopical Khovanov homology for links in thickened surfaces with genus$>1$.

\bs
\h $(2)$ We define
the second Steenrod square acting on
the  homotopical Khovanov homology for links in thickened surfaces  with genus$>1$
by using the Khovanov-Lipshitz-Sarkar stable homotopy type in $(1)$.

\bs
\h $(3)$ Each of
the Khovanov-Lipshitz-Sarkar stable homotopy type in $(1)$
and
the second Steenrod square in $(2)$
is stronger than
the  homotopical Khovanov  homology
as invariants of links in thickened surfaces  with genus$>1$.
That is,
there is a pair of links in a thickened surface  with genus$>1$
with the following properties:
They have different Steenrod squares.
They have different Khovanov-Lipshitz-Sarkar stable homotopy types.
They have the same homotopical Khovanov homology.
\end{mth}

\begin{note}\label{notevirtual}
There are other ways to extend Khovanov homology of links in the 3-sphere
and the second Steenrod square on Khovanov homology of links in the 3-sphere
into thickened surfaces.
Manturov (2006 arXiv)
\cite{Man},
Rushworth
\cite{Ru},
Tubbenhauer
\cite{Tub}, and
Viro (unpublished)
\cite{Viro}
  introduced
 other ways to define
 Khovanov homology for links in thickened surfaces
 by using virtual links.
Dye, Kaestner, and Kauffman
\cite{DKK},
 Nikonov \cite{Igor},
 Kauffman and Ogasa
\cite{KauffmanOgasasq}
wrote alternative definitions of \cite{Man}.
 Kauffman and Ogasa
\cite{KauffmanOgasasq}
 extended
 the second Steenrod square of link in $S^3$
 to thickened surfaces by using virtual links.
See \cite{Kauffman1, Kauffman, Kauffmani} for Virtual links.
 However,
 the readers need not know virtual links or
 the result of any paper in this Note \ref{notevirtual}
 in order to read this paper.
 The techniques for defining Khovanov homology for virtual links are not used in the present paper. Our paper is self-contained, and uses the techniques of Lifshitz and Sarkar.

On the other hand,  for the sake of those who are familiar with virtual links,
we will sometimes %point out differences between theory of `links in thickened surfaces without using virtual links' and that of virtual links, in this paper.
comment on virtual links in this paper.
Of course, the readers may skip the parts
if they do not know virtual links.
\end{note}

\bb
\section{\bf The homotopical  Khovanov
 homology for links in thickened surfaces}\label{KMNhomology}

\h  We review the definition of the homotopical  Khovanov
 homology for links in thickened surfaces introduced in \cite{MN}.
In this paper we define Khovanov homotopy type for
thickened higher genus surfaces, but
  the homotopical  Khovanov
  homology is defined for links in thickened surfaces with any genus.
In this paper, when we say just a surface,
its genus may zero, one, or greater than one.

We translate the definition %\cite[(2.3)-(2.6)]{KMNhomology}
into terminologies %in \cite{LSk}.
%We write many definitions in this section by
 %terminologies and methods of
 in  Lipshitz and Sarkar's paper \cite{LSk},
because we generalize the results about
Khovanov-Lipshitz-Sarkar stable homotopy type there.

\bb
\subsection{Labeled resolution configurations
}\label{lsc}\hskip10mm

%We first define resolution configurations. They are much connected with the theory of links, but we define them without using links.

\begin{defn}\label{2.1}
Let $F$ be a closed oriented surface.
A {\it resolution configuration} $D$ is a pair $(Z(D),A(D))$,
where $Z(D)$ is a set of pairwise-disjoint embedded circles in $F$,
and $A(D)$ is a totally ordered collection of disjoint arcs embedded in $F$,
with $A(D)\cap Z(D)=\partial A(D)$.
We call the number of arcs in $A(D)$ the {\it index} of the resolution configuration $D$,
and denote it by ind$(D)$.
We sometimes abuse notation and write $Z(D)$
to mean $\cup_{Z\in Z(D)} Z$ and $A(D)$ to mean $\cup_{A\in A(D)} A$.
Occasionally, we will describe the total order on $A(D)$
by numbering the arcs: a lower numbered arc precedes a higher numbered one.
\end{defn}

We changed \cite[Definition 2.1]{LSk} into Definition \ref{2.1}
by replacing `$S^2$ in \cite[Definition 2.1]{LSk}' by `$F$ in Definition \ref{2.1}'.

\begin{defn}\label{2.2}
Let $F$ be a surface.
%In this  paper, we abbreviate a closed orineted surface to a surface when it is clear from the context.
%{\bf (A modification of Definition 2.2 in \cite{LSk}.)}
Let $\mathcal L$ be a link in $F\x[-1,1]$.
Let $L$ be a link diagram of  $\mathcal L$.
Note that $L$ is in $F$.
Assume that
 the link diagram $L$ has $n$ crossings,
an ordering of the crossings in $L$,
and a {\it vector} $v\in\{0,1\}^n$.
There is an {\it associated resolution configuration} $D_L(v)$
obtained by taking the resolution of $L$ corresponding to $v$
(that is, taking the 0-resolution at the $i$-th crossing
if $v_i =0$, and the 1-resolution otherwise)
and then placing arcs corresponding to each of the crossings labeled
by 0's in $v$
(that is, at the $i$-th crossing if $v_i=0$). See Figure \ref{r}.\\
\begin{figure}
\includegraphics[width=120mm]{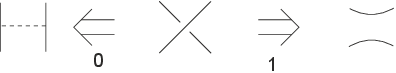}
\caption{{\bf The 0- and 1-resolutions
}\label{r}}
\end{figure}
Therefore, $n-$ind$(D_L(v)) = |v| =\sum v_i$, the (Manhattan) norm of $v$.
(Note that
$\sum v_i\\=\sum (v_i)^2$
in this situation,
since $0^2=0$ and $1^2=1$.)
\end{defn}

Note that
(if $L$ in $F$ denotes a classical link diagram in  $S^2$ or virtual link diagram)
resolution configurations are the same
as what many people often call {\it Kauffman state}s,
which Kauffman first introduced in \cite{Kauffmanstate}.
(The way to draw arcs in \cite{LSk} is different from that in \cite{Kauffmanstate}.)
\\

\h{\bf Note.}
We draw resolution configurations on the plane $\R^2$ in this  paper.
%to explain what kind of them we discuss.
Our way of drawing is similar to that in virtual knot theory
which is  introduced in \cite{Kauffman1, Kauffman, Kauffmani}.
 However,
 in this  paper,  the strict definition of resolution configurations
 is Definition \ref{2.1}.
(On the other hand,
the way of drawing resolution configurations on $\R^2$ (respectively,  $S^2$)
is defined strictly
 in virtual knot theory.)

%After we carry out resolutions on knot diagrams,

\begin{figure}
\includegraphics[width=120mm]{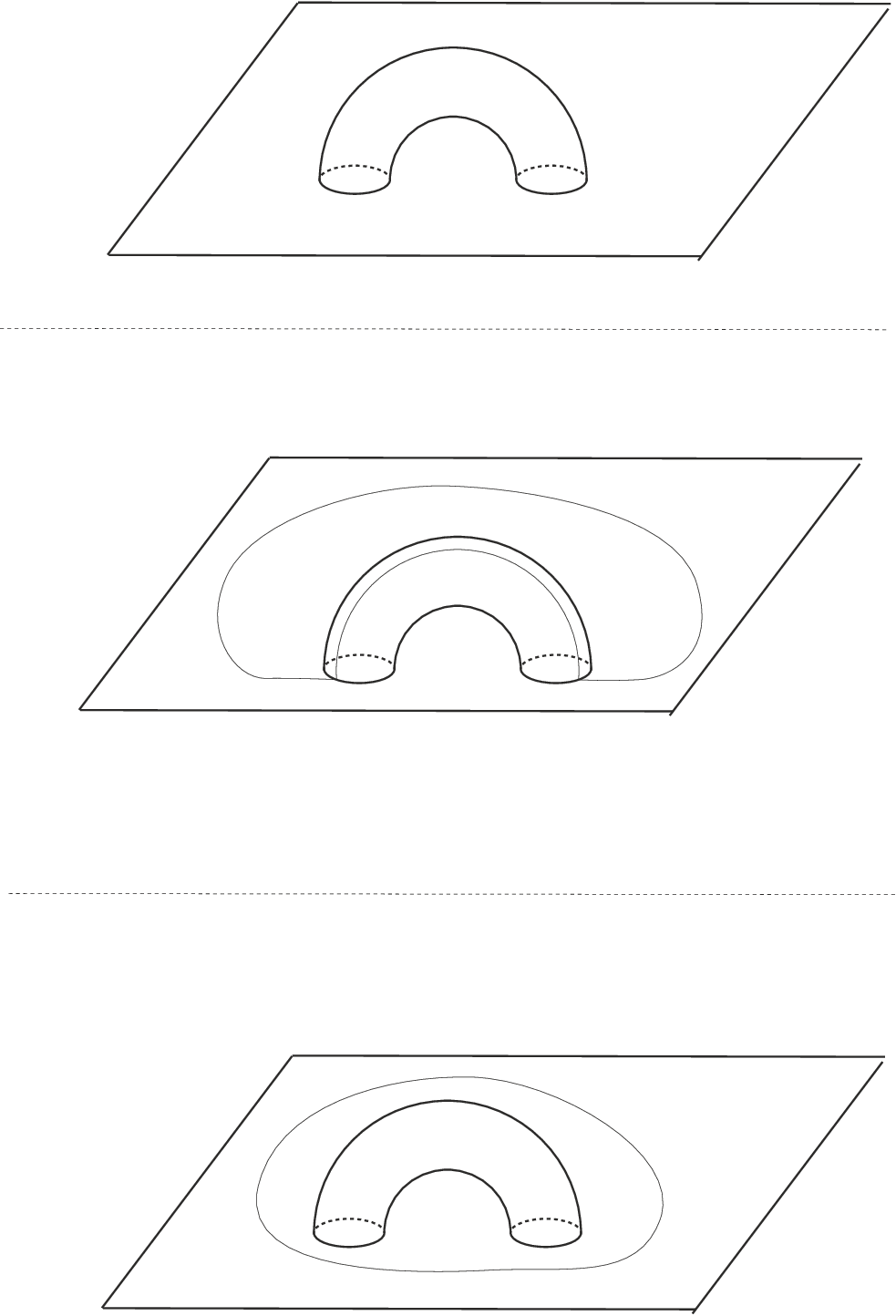}
\caption{{\bf
The upper figure is a part of a surface $F$.
%The middle (respectively, lower) figure is
In each of the middle and lower figures,
there is a non-contractible circle on $F$.
}\label{handle}}
\end{figure}

\begin{figure}
\bb
\includegraphics[width=120mm]{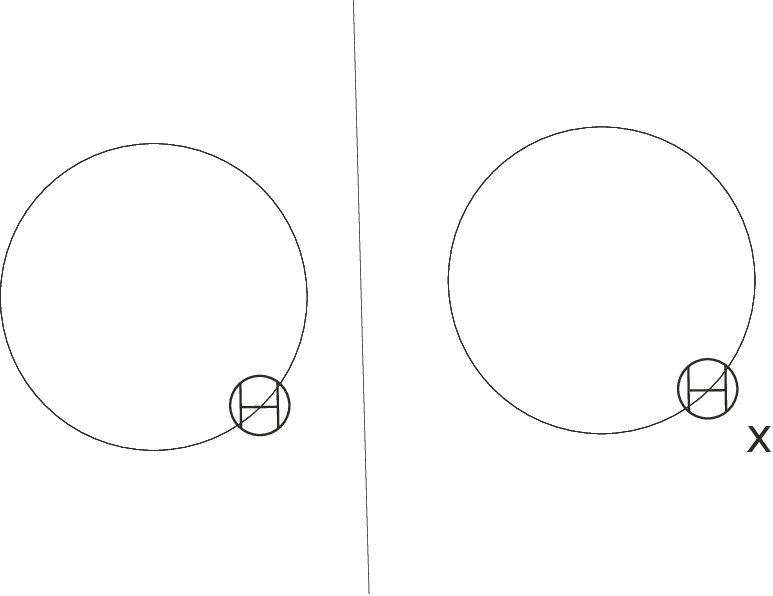}
\caption{{\bf Non-contractible circles}\label{H}}
\bb
\end{figure}

We draw a part of a surface $F$ %as drawn
in the upper figure of Figure \ref{handle}.
We depict a non-contractible circle in $Z(D)$  on $F$
in the middle and lower figures  of Figure \ref{handle}.
These kinds of non-contractible circles are drawn as in the left figure of Figure \ref{H}.
We call this circle a {\it circle with $(H)$}.
We omit drawing a part of the surface $F$ when it is clear from the context.
If we need to explain some property %, e.g. a homotopy class,
of
a non-contractible circle,
we write it in the right lower side where $x$ is written in the right one of Figure \ref{H}.
%(In this example, it is $x$.).

%In a situation, we use a net of surfaces in which link diagrams exist, and draw labeled resolution configurations. See e.g. Figure \ref{quasiT2}.

\begin{figure}
\bb
\includegraphics[width=120mm]{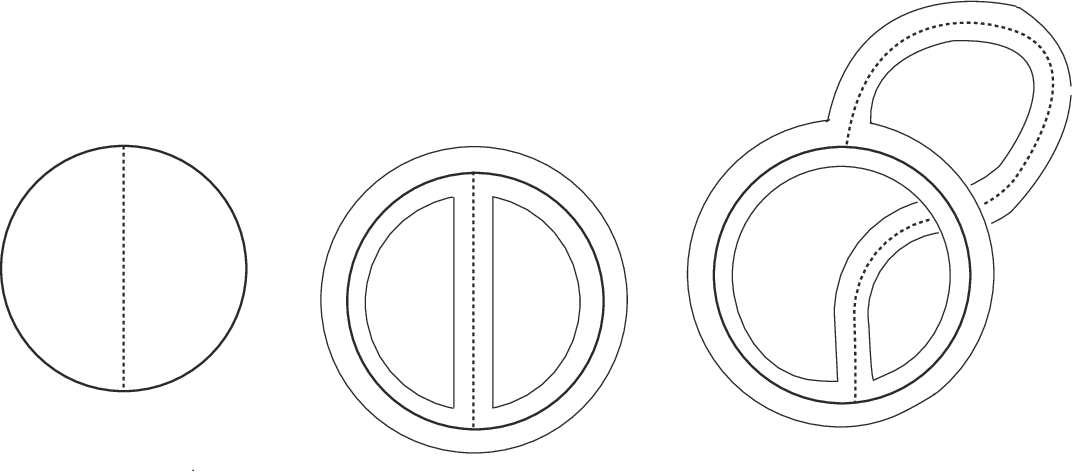}
\caption{{\bf
%Two cases of neighborhoods of the left Kauffman state.
Two Kauffman states with a single circle and a single arc:
The two  are different but make the same abstract graph.
}\label{nbd}}
\bb
\end{figure}

%In another situation,
See Figure \ref{nbd}.
$Z(D)\cup A(D)$ has a neighborhood $N$
such that $N$ is a compact surface and
such that the inclusion map of $Z(D)\cup A(D)$ to $N$
 is a homotopy type equivalence map.
For a given $Z(D)\cup A(D)$, there are many homeomorphism types of $N$ in general.

%Figure \ref{zenbu} We draw only $Z(D)\cup A(D)$ abstractly, like abstract graphs. \\
%the tubular neighborhood of knot diagrams in surfaces. %See \S\ref{KLShomotopy}.
%Note: $Z(D)\cup A(D)$ has a neighborhood $N$
%such that $N$ is a compact surface and
%such that the inclusion map of $Z(D)\cup A(D)$ to $N$  is a homotopy type equivalence map.
%For a given $Z(D)\cup A(D)$, there are many homeomorphism types of $N$ in general.
%See Figure \ref{nbd} for an example.

\begin{defn}\label{2.3}
{\bf (\cite[Definition 2.3]{LSk}.)}
Given resolution configurations $D$ and $E$,
there is a new resolution configuration
$D-E$ defined by
$$Z(D-E)=Z(D)-Z(E)\hskip3mm A(D-E) = \{A\in A(D) |
\forall Z\in Z(E): \partial A\cap Z=\emptyset
\}.$$
Let $D\cap E=D-(D-E)$.
\end{defn}

Note that $Z(D\cap E)=Z(E\cap D)$ and
$A(D\cap E)=A(E\cap D)$; however, the total orders
on $A(D\cap E)$ and $A(E\cap D)$ could be different.

\begin{defn}\label{2.4}
{\bf (\cite[Definition 2.4]{LSk}.)}
The {\it core} $c(D)$ of a resolution configuration $D$ is the resolution configuration
obtained from $D$ by deleting all the circles in $Z(D)$
that are disjoint from all the arcs in $A(D)$.
A resolution configuration $D$ is called {\it basic} if $D = c(D)$, that is,
if every circle in $Z(D)$ intersects an arc in $A(D)$.
\end{defn}

\begin{defn}\label{surger-res-config}
{\bf (\cite[Definition 2.5]{LSk}.)}
  Given a resolution configuration $D$ and a subset $A'\subseteq A(D)$
  there is a new resolution configuration $s_{A'}(D)$, the
  {\it surgery of $D$ along $A'$}, obtained as follows.
  The circles
  $Z(s_{A'}(D))$ of $s_{A'}(D)$ are obtained by performing embedded
  surgery along the arcs in $A'$; in other words, $Z(s_{A'}(D))$ is
  obtained by deleting a neighborhood of $ A'$ from $Z(D)$ and
  then connecting the endpoints of the result using parallel
  translates of $A'$.
  The arcs of $s_{A'}(D)$ are the arcs of $D$ not
  in $A'$, i.e., $A(s_{A'}(D))=A(D)- A'$.
  %See \Figure{res-config-all-c}.

  Let $s(D)=s_{A(D)}(D)$ denote the maximal surgery on $D$.
\end{defn}

\begin{defn}%\label{scs}
Let $D$ be a resolution configuration.
Suppose that, when we carry out a surgery along %a single
one arc of $A(D)$ on circles of $Z(D)$,
the number of the elements of $Z(D)$ is not changed.
Then we call this surgery a {\it single cycle surgery}.
We call this arc a {\it scs arc}.
\end{defn}

There is a scs arc in the right figure of Figure \ref{nbd}.

\begin{figure}
\includegraphics[width=140mm]{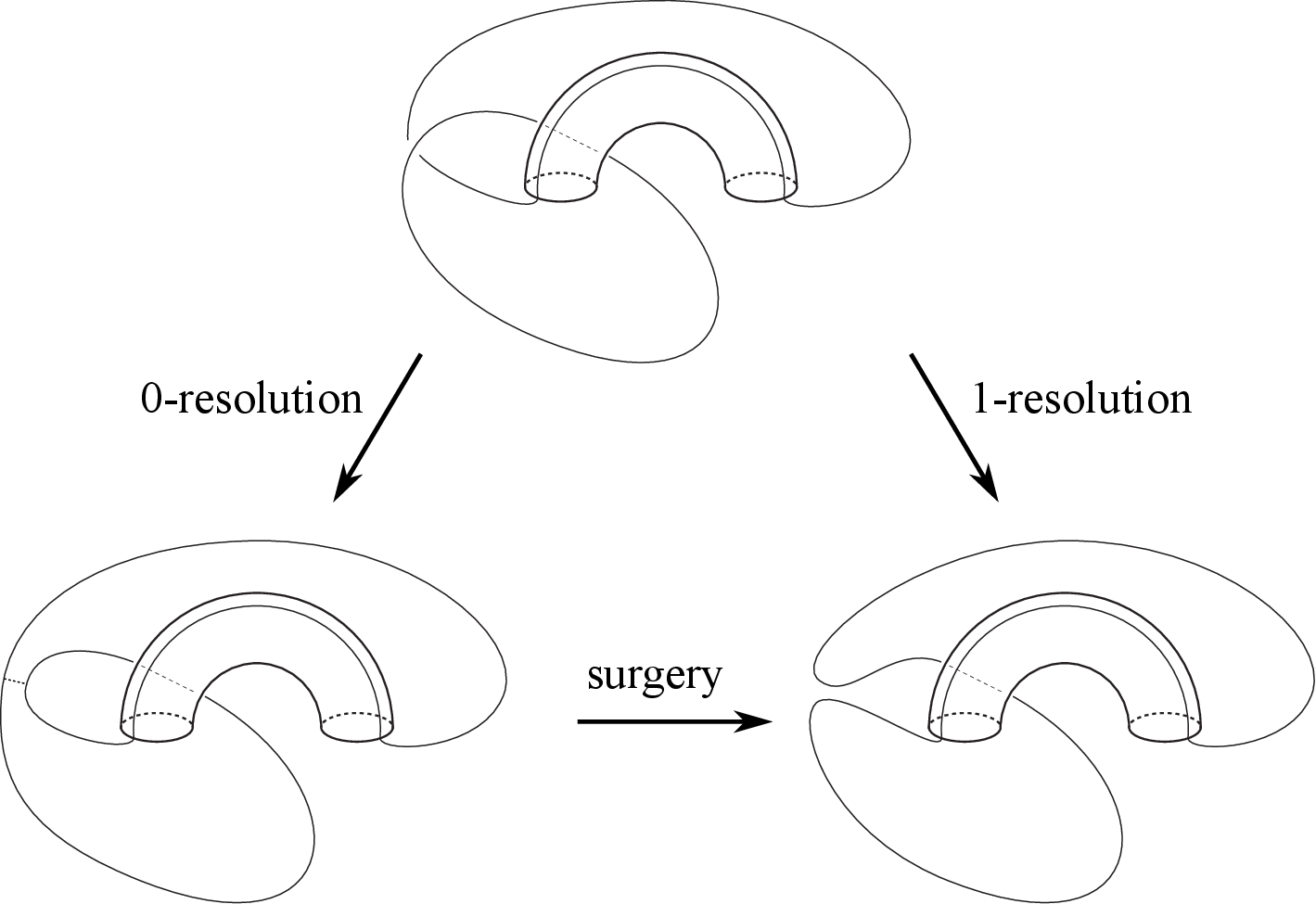}
\caption{{\bf A single cycle surgery.
}\label{vHrei}}
\end{figure}

See another example in Figure \ref{vHrei}:
Each of these three figures are a part of a surface $F$.
The upper is a link diagram in $F$.
%the projection of a link in $F\x[-1,1]$ to $F$.
The middle is obtained from the upper by the 0-resolution.
The lower is obtained from the upper by the 1-resolution.
The lower is obtained from the middle by a surgery along the arc in the middle.
This surgery is a single cycle surgery.

\begin{figure}
\bb
\includegraphics[width=100mm]{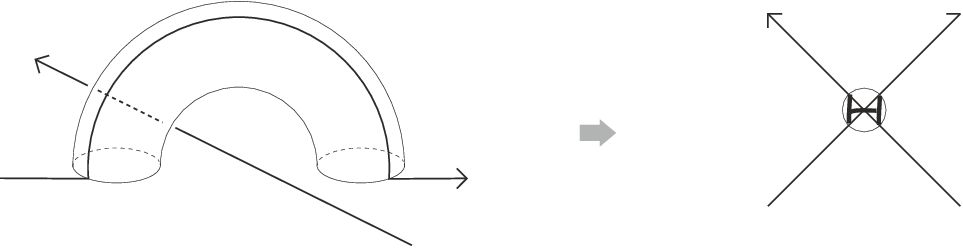}
\caption{{\bf A way of drawing a non-contractible circle in a resolution configuration.
}\label{vH}}
\bb
\end{figure}

%Before the proof, we define some terms.
\bb
See Figure \ref{vH}.
The left figure is a part of a closed oriented surface $F$ %$T^2$
with
a part of a link diagram in $F$. %\x[-1,1]$. %$T^2\x[-1,1]$.
%a part of the projection of a link in $F\x[-1,1]$. %$T^2\x[-1,1]$.
We draw it as the right one for %the
convenience when we discuss a single cycle surgery.
The right one includes $(H)$.
This $(H)$ represents not only the fact in Figure \ref{vH}
but also the fact that the circle is a non-contractible circle as written in Definition \ref{2.2}.

\begin{figure}
\bb
\includegraphics[width=100mm]{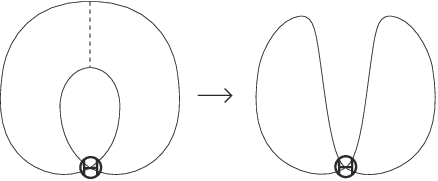}
\caption{{\bf A single cycle surgery
drawn by using  $(H)$.
}\label{scs}}
\bb
\end{figure}

Figure \ref{scs} is an example  of
drawing a single cycle surgery by using $(H)$ of Figure \ref{vH}.

\begin{defn}\label{arcsu}
%If an arc in a labeled resolution configuration produces a single cycle surgery, the arc is called a {\it scs arc}.
If a surgery along an arc increases (respectively, decreases) the number of circles by one,  the surgery is called  a {\it comultiplication} (respectively, {\it multiplication}).
There are just three kinds of surgeries along an arc: a single cycle surgery,
a multiplication, a comultiplication.
If an arc produces a multiplication (respectively, comultiplication),
the arc is called a {\it m-arc} (respectively, {\it c-arc}).
If an arc is an m-arc or a c-arc, that is, it is not a scs arc,
then the arc is called a {\it mc arc}.
\end{defn}

\begin{defn}\label{2.9}  {\bf (\cite[Definition 2.9]{LSk}).}
A {\it labeled resolution configuration} is a pair $(D, x)$ of a resolution configuration
$D$ and a labeling $x$ of each element of $Z(D)$ by either $x_+$ or $x_-$.
\end{defn}

Note that
(if they are associated with classical diagrams in $S^2$ or virtual link diagrams,)
 labeled resolution configurations are the same
as what many people often call
{\it enhanced Kauffman state}s
or
{\it enhanced state}s.
%(The ways to draw arcs associated with resolutions in both papers are different.)
Some people use
$v_+$ (respectively, $v_-$)
for
$x_+$ (respectively, $x_-$).
\\

Let $\{A_i\}_{i\in\Lambda}$ be the set of all labeled resolution configurations made from
an arbitrary link diagram $L$ in a surface.
Note that $\Lambda$ is a finite set.
$\{A_i\}_{i\in\Lambda}$ composes
a basis of the Khovanov homology for $L$ %link diagrams in surfaces
as in \cite{APS, B, DKK, KauffmanOgasasq, K, LSk}.
We call $\{A_i\}_{i\in\Lambda}$  the {\it Khovanov basis}.
We call each  $A_i$ a {\it Khovanov basis element}.
\\

We will define a partial order on $\{A_i\}_{i\in\Lambda}$.
After that, by using the partial order,
we will define the differential acting on each $A_i$,  %labeled resolution configurations,
and
introduce the Khovanov homology  for $L$
as in \cite{LSk, KauffmanOgasasq}.
See the definitions in the following subsections % \S\ref{pos}
for the detail.
In order to define moduli spaces and to construct
Khovanov-Lipshitz-Sarkar stable homotopy type,
we need the partial ordered set.
\\

\h{\bf Review.}
%The following two  ways are the same.
Comparing \cite{B,K} with \cite{LSk},
a differential of the Khovanov complex and a partial order on the set of labeled resolution configurations are essentially the same thing. We explain it below.
\bs

Let  $\{X_i\}_{i\in\Theta}$ be a set of  all labeled resolution configurations
 of an arbitrary link diagram in $S^2$.
We want to let the discussion in this Review separated from the one above here,
so we use different notation  $\{X_i\}_{i\in\Theta}$.
\bs

\h(1)
The Khovanov differential acting on  $\{X_i\}_{i\in\Theta}$
in \cite{K, B} induces a partial order on  $\{X_i\}_{i\in\Theta}$
as below.
We use the following notation

\begin{equation}\label{keisu}
\delta X_i = \displaystyle\sum_{j\in\Theta}
[X_i;X_j]\cdot X_j,
%c[X_i;X_j]\cdot X_j,
\end{equation}

\h
and
$[X_i;X_j]$
%$c[X_i;X_j]$
is an integer coefficient.
Recall that
$[X_i;X_j]\in\{-1,0,1\}$.
%$c[X_i;X_j]\in\{-1,0,1\}$.
The partial order $\prec$ on
 $\{X_i\}_{i\in\Theta}$ is induced as follows:
If
$[X_i;X_j]=\pm1$,
%$c[X_i;X_j]=\pm1$,
then $X_i\prec X_j$.
It follows from this definition that
if $u \prec v$ and $v \prec w$ then $u \prec w$.
\\

\h(2)
In \cite[Definition 2.10]{LSk}, Lipshitz and Sarkar introduced
a partial order $\prec$ on the set $\{X_i\}_{i\in\Theta}$
%of all labeled resolution configurations of an arbitrary link diagram in $S^2$. in $S^2$
before they define the Khovanov differential acting on each $X_i$. %$\{X_i\}_{i\in\Theta}$.

Furthermore, they use the following method:
Let $Y, Z\in\{X_i\}_{i\in\Theta}$. Let $Y\prec Z$.
Assume that $Y\prec X\prec Z$ does not hold for any $X\in\{X_i\}_{i\in\Theta}-\{Y,Z\}$.
 They define an  explicit way to
assign  $+1$ or $-1$
to the pair,  $Y$ and $Z$.
% (This method is also used in \cite{B,K}.)

By using this partial order and this method to give a sign,
they induce the Khovanov differential,
acting on
 $\{X_i\}_{i\in\Theta}$
 (\cite[Definition 2.15]{LSk}).
\\

We follow Lipshitz and Sarkar's method above, %in \cite[Definition 2.10]{LSk}
in the following subsections.  %\S\ref{pos}.

Lipshitz and Sarkar used the partial ordered set, defined moduli spaces, and constructed
Khovanov-Lipshitz-Sarkar stable homotopy type.

\bb
\subsection{
A partial order on the set of labeled resolution configurations
}\label{pos}\hskip10mm

\begin{defn}\label{induce}
Let  $(D, x)$ and $(E, y)$ be labeled resolution configurations.
Then there is a natural labeling $x|$ (respectively, $y|$) on $D-E$ (respectively, $E-D$),
say the {\it induced labeling}.
We often consider  induced labelings in the following case:
if we restrict each of the labelings $x$ and $y$ to $D\cap E = E\cap D$,
we obtain the same labeling from $x$ and $y$.
\end{defn}

We define a partial order on $\{A_i\}_{i\in\Lambda}$ in the following definition.
By using the partial order,
we will induce the differential (Definition \ref{korekos}, cited from \cite[\S2]{MN}).
One can induce this partial order from the differential: Use
\cite[(2.1), (2.3)-(2.6)]{MN}.
(One can do it in a similar way in the Review in the last part of \S\ref{lsc}.)

\begin{figure}
%\vskip-10mm
\includegraphics[width=170mm]{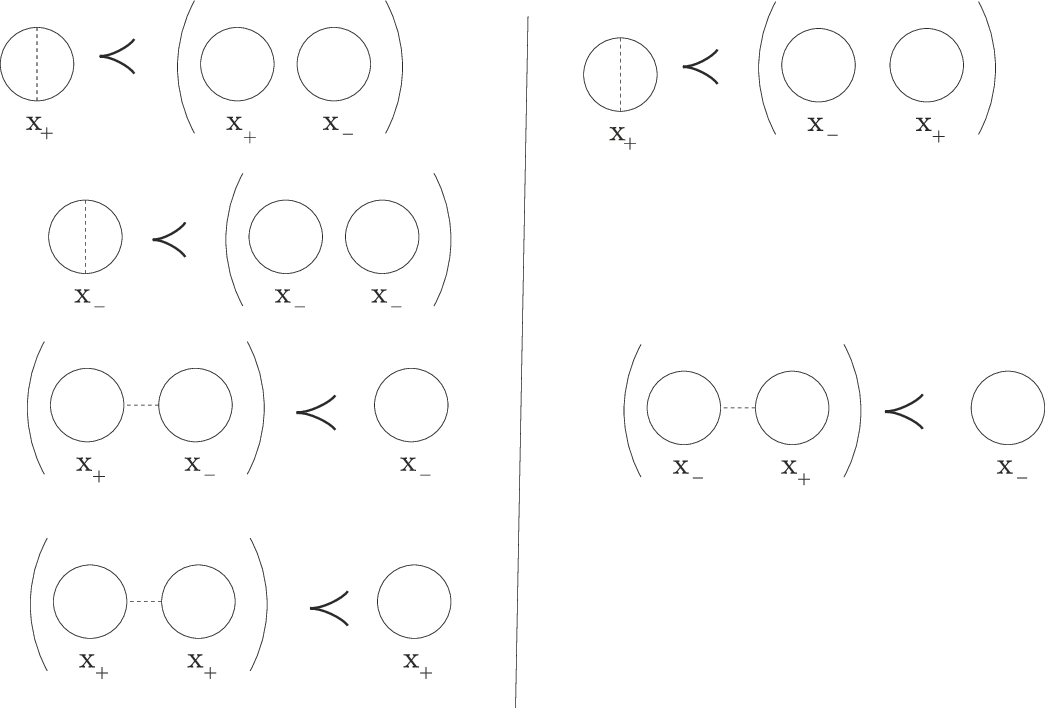}
\caption{{\bf
The partial order of the set of labeled resolution configurations
in the upper three of Table \ref{PtoQR} and in the upper three of Table \ref{PQtoR}.
%is made by the rule in this figure and Figure \ref{resolNC}.
%In this Figure, we draw the case where all circles are contractible.
See Figures \ref{resolreiC} and \ref{resolreiNC}.
}\label{resolC}}
\end{figure}

\begin{figure}
\vskip-10mm
\includegraphics[width=130mm]{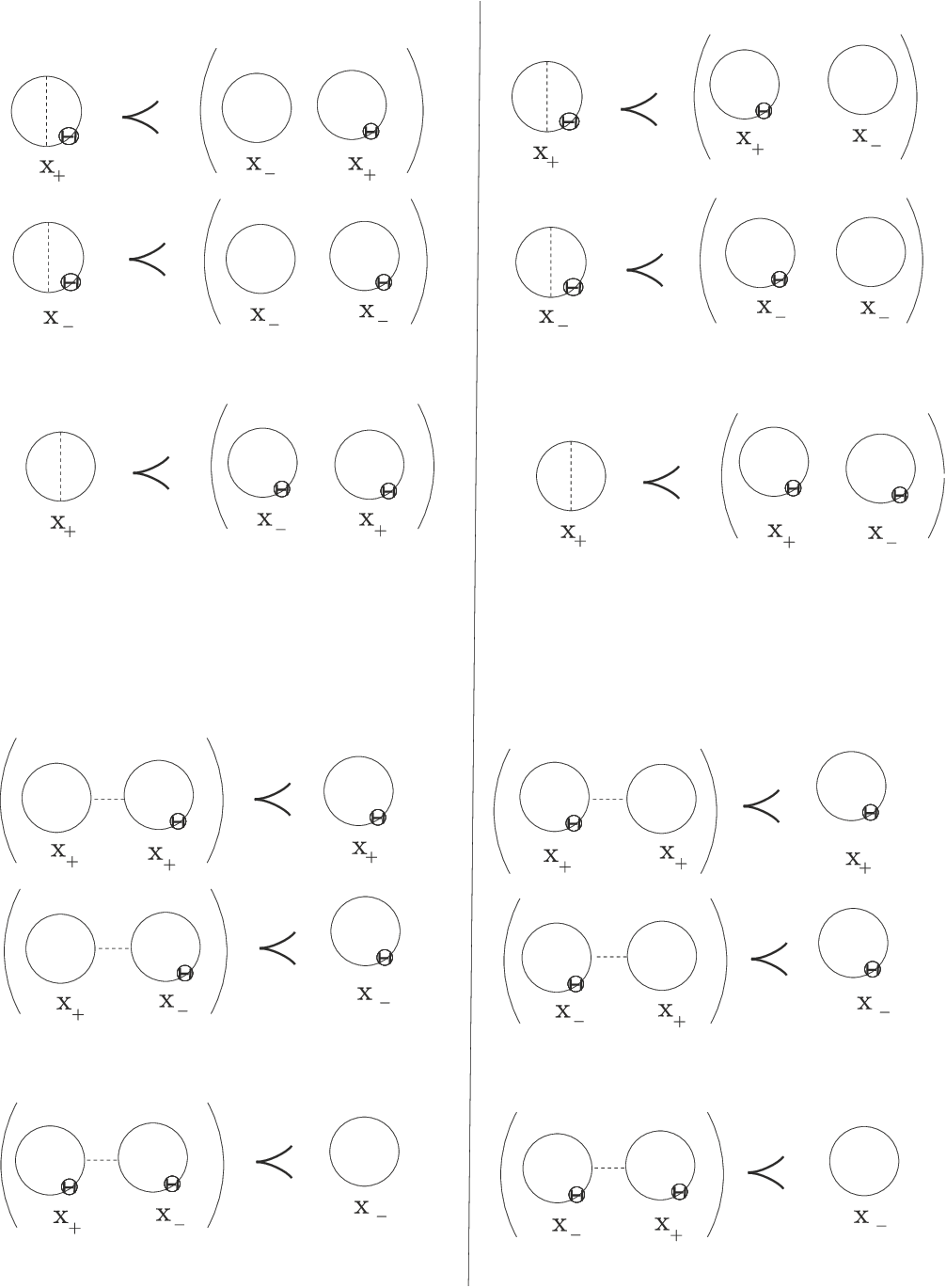}
\caption{{\bf
The partial order of the set of labeled resolution configurations
in the lower six of Table \ref{PtoQR} and in the lower six of Table \ref{PQtoR}:
%is made by the rule in this figure and Figure \ref{resolC}.
%In this Figure, we draw the case where there are non-contractible circles.
There appear non-contractible circles.
 Circles with (H) denote non-contractible circles.
See Figures \ref{resolreiC} and \ref{resolreiNC}.
}
\label{resolNC}}
\end{figure}

\begin{defn}\label{2.10}  %{\bf (\cite[Definition 2.10]{LSk}).}
There is a partial order $\prec$ on labeled resolution configurations defined as
follows.
Note that  $\alpha\prec\alpha$ holds.
Then we declare that $(E, y)\prec(D, x)$ if:

\bs
\begin{enumerate}
\item[(1)]
 The labelings $x$ and $y$ induce the same labeling on $D\cap E = E\cap D.$

\bs
\item[(2)]
$D$ is obtained from $E$ by surgering along a single arc of $A(E)$.

\h
$Z(E-D)$ (respectively,  $Z(D-E)$)
has an induced labelling
$y|$ (respectively, $x|$)
from
$y$ (respectively, $x$).
There are two sub-cases (i) and (ii).

\bs
\h
(i) $(E-D,y|)$ has just one circle $P$, and
 $(D-E,x|)$ has just two circles, $Q$ and $R$.
 %They
 These two labeled resolution configurations
 satisfy  the conditions in Table \ref{PtoQR}.
Here, c means a contractible circle,
and n means a non-contractible circle.
Examples are
the upper three figures in Figure \ref{resolC}
and the upper six figures in Figure \ref{resolNC}.

\bb

\begin{table}[htb]

  \begin{tabular}{ccccc}
    $P$  & $\prec$ & $Q$& $\&$ & $R$ \\ \hline
    c, $x_+$&  $\prec$ & c, $x_+$& $\&$& c, $x_-$\\
    c, $x_+$&  $\prec$ & c, $x_-$& $\&$& c, $x_+$\\
    c, $x_-$&  $\prec$ & c, $x_-$& $\&$& c, $x_-$\\
              &              &            &        & \\
  n, $x_+$ &  $\prec$  & c, $x_-$ & \&& n, $x_+$\\
  n, $x_+$ &  $\prec$  & n, $x_+$ &\& & c, $x_-$\\
   n, $x_-$&  $\prec$  & c, $x_-$ & \&& n, $x_-$\\
   n, $x_-$&   $\prec$ & n, $x_-$ & \&& c, $x_-$\\
  c, $x_+$ & $\prec$  & n, $x_-$ & \&& n, $x_+$\\
  c, $x_+$ & $\prec$  & n, $x_+$ & \&& n, $x_-$\\
    \end{tabular}
    \bb
\caption{{\bf $P$, $Q$, $R$ in Definition \ref{2.10}.(2).(i).
}\label{PtoQR}}
\end{table}

\bb
\h
(ii) $(E-D,y|)$ has just two circles, $P$ and $Q$, and
 $(D-E,x|)$ has just one circle $R$.
 %They
  These two labeled resolution configurations
satisfy  the conditions in Table \ref{PQtoR}.
%Here, c means a contractible circle, and n means a non-contractible circle.
Examples are
the lower three figures in Figure \ref{resolC}
and the lower six figures in Figure \ref{resolNC}.

\bb

\begin{table}[htb]
\begin{tabular}{ccccc}
          $P$&\&  &   $Q$& $\prec$&$R$  \\ \hline
    c, $x_+$& \& &c, $x_-$ &$\prec$ &c, $x_-$  \\
    c, $x_-$&\&  &c, $x_+$ &$\prec$ &c, $x_-$  \\
    c, $x_+$&\&  &c, $x_+$ &$\prec$ &c, $x_+$  \\
    & &   & & \\
    c, $x_+$&\&  &n, $x_+$ &$\prec$ &n, $x_+$  \\
     n, $x_+$&\&  &c, $x_+$& $\prec$ &n, $x_+$  \\
     c, $x_+$&\&  &n, $x_-$& $\prec$ &n, $x_-$  \\
     n, $x_-$&\&  &c, $x_+$& $\prec$ &n, $x_-$  \\
     n, $x_+$&\&  &n, $x_-$& $\prec$ &c, $x_-$  \\
     n, $x_-$&\&  &n, $x_+$& $\prec$ &c, $x_-$  \\
 \end{tabular}
\bb
\caption{{\bf $P$, $Q$, $R$ in Definition \ref{2.10}.(2).(ii).
}\label{PQtoR}}
\end{table}

\end{enumerate}

Now, $\prec$ is defined to be the transitive closure of this relation.
\end{defn}

\h{\bf Note.}
(1)
The upper three relations in
each of Tables \ref{PtoQR} and \ref{PQtoR}
are the same as those in  \cite[Definition 2.10]{LSk}.
The lower six in each is associated with the case where
%contractible and
non-contractible circles appear.
These cases conform to the conditions for defining the boundary in
\cite[(2.1), (2.3)-(2.6)]{MN}.

\bs

\h(2) In Definition 2.9 the relations for homotopy classes of  circles are determined naturally.
(When we consider a homotopy class of the circles,
we
%give them orientations
let them orient %ations
appropriately.)
\bs

\h(3) No surgery in Definition \ref{2.10} is a single cycle surgery although we now consider links in thickened %oriented closed
surfaces.
\bs

\begin{figure}%[t]
\bb
\includegraphics[width=70mm]{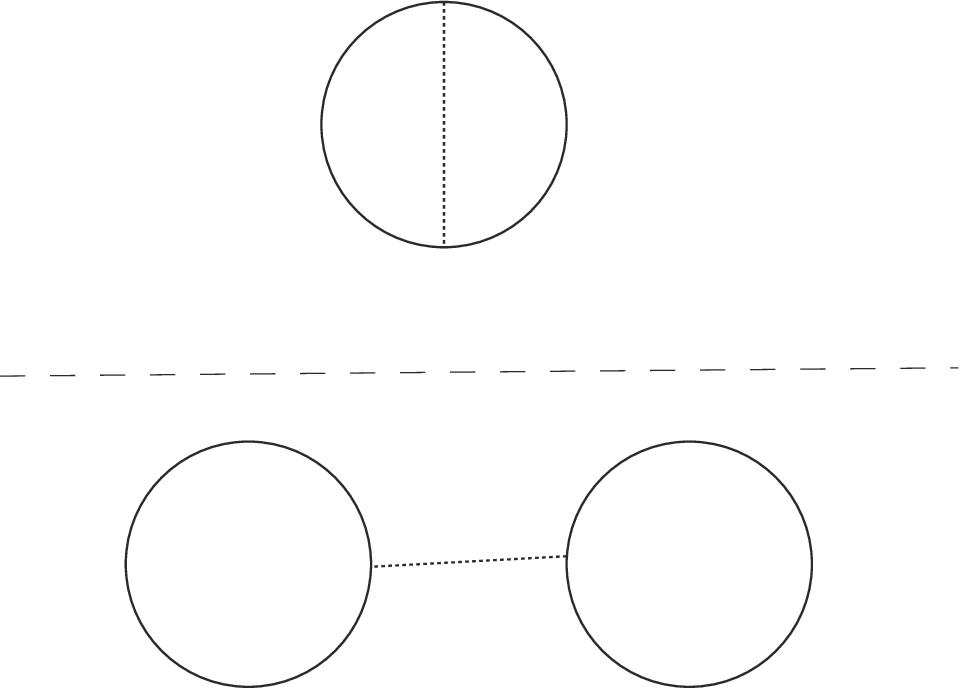}
\caption{{\bf
Graphs of the basic resolution configurations of index 1.
%We draw abstractly $Z(E-D)\cup A(E-D)$ in the case where  $A(E-D)$ has a single arc. There are only two cases.
}\label{1arc}}
\bb
\end{figure}

\h(4) In Figure \ref{1arc},
we draw abstractly $Z(E-D)\cup A(E-D)$
in the case where  $A(E-D)$ has a single arc. There are only two cases.
\bs

\begin{defn}\label{2.11}   {\bf (\cite[Definition 2.11]{LSk}.)}
A {\it decorated resolution configuration} is a triple $(D, x, y)$
where $D$ is a resolution configuration and
$x$ (respectively,  $y$) is a labeling of each component of $Z(s(D))$
(respectively, $Z(D)$) by an element of $\{x_+, x_-\}$.
%It means a partial ordered set
Associated to a decorated resolution configuration $(D, x, y)$ is the poset $P(D, x, y)$
 consisting
of all labeled resolution configurations
$(E, z)$ with
$(D, y)\prec(E, z)\prec(s(D), x)$.
We call $P(D, x, y)$ the poset for $(D, x, y)$.

If we never have $(D, y) \prec (s(D), x)$,
we say that $(D, x, y)$ is empty. % and write $(D,x,y)\\=\emptyset$.
%the empty set.
\end{defn}

Note:
In \cite[Definition 2.11]{LSk},
%if we write $(D, x, y)$, then $(D, x, y)\neq\emptyset$.
%
% Lipshitz and Sarkar let $(D, x, y)$ just represent $(D, y)$ and $(s(D), x)$ above, and call the poset $P(D,x,y)$.   On the other hand, in this paper,   $(D, x, y)$ also means $P(D,x,y)$.
   Lipshitz and Sarkar
%Furthermore they
dealt with only non-empty decorated resolution configurations.
It means that, in their paper, they only consider
 labelings $x$ on  $s(D)$ such that
 $(D, y) \prec (s(D), x)$.
 On the other hand, in this paper, we define the case
where $(D, x, y)$ is empty,
 %$(D, x, y)=\emptyset$
 for convenience.

\begin{defn}\label{2.15gr}
{\bf(\cite[A part of Definition 2.15]{LSk}, and \cite[Definition 2.2]{MN}.)}
Let $L$ be an oriented link diagram in a surface $F$.
Let $n$ (respectively, $n_+$, $n_-$) be
the number of crossing  (respectively, positive crossing, negative crossing)
points of $L$.

For labeled resolution configurations,
{\it homological grading $\text{gr}_h$},
{\it a quantum grading $\text{gr}_q$}, and
{\it a homotopical grading $\text{gr}_\mathfrak H$}
are defined as follows:

%\vskip2mm\hskip7mm
\begin{equation}
\text{gr}_h((D_L(u), x)) = -n_- + |u|,
\end{equation}

\begin{equation}
\text{gr}_q((D_L(u), x)) = n_+ - 2n_- + |u|+
\sharp\{Z\in Z(D_L(u)) | x(Z) = x_+\}
-\sharp\{Z\in Z(D_L(u)) | x(Z) = x_-\}.
\end{equation}

%\vskip2mm\hskip7mm
%$\text{gr}_q((D_L(u), x)) = n_+ - 2n_- + |u|$

%\hskip40mm
%$+
%\sharp\{Z\in Z(D_L(u)) | x(Z) = x_+\}
%-\sharp\{Z\in Z(D_L(u)) | x(Z) = x_-\}.$
%\vskip2mm

%Here $n_+$ denotes the number of positive crossings in $L$; and $n_- = n-n_+$ denotes the number of negative crossings.\\

%Let $L$ be a $S$ be a closed oriented surface.
We consider the set $\mathfrak L = [S^1; F]$
of all the homotopy classes of free oriented loops in $F$.  %$S$.
Let $\bigcirc\in \mathfrak L$ be the homotopy class of contractible loops.
For any closed curve $\gamma$,
one can consider the curve $-\gamma$ obtained from $\gamma$
by the orientation change.
Let $\mathfrak H$ be the quotient group of the free abelian group
with generator set $\mathfrak L$ modulo the relations
 $\bigcirc= 0$ and $[\gamma]=[-\gamma]$ for all free loops $\gamma$.
%
%
%We consider a new homotopical grading in the complex $[[L]]$ valued
%in the group
Let $Z(D_L(u))=\{C_1, C_2,...,C_\nu\}$.  %\\  %{\gamma(s)}\}$.
%Let $x_i$ be as follows: $C_i$ equips $x_i$,
Let $C_i$ equip $x_i$,
where  $x_i \in\{x_+, x_-\}$ and $i\in\{1,...,\nu\}$. %\\ %\gamma(s)\}$. \\
%Let $s$ be a state and $C_1, C_2,...,C_{\gamma(s)}$ be the circles of the resolution $L_s$.
%For each element $x=x_1\otimes x_2 \otimes...\otimes x_{\gamma(s)}\in V(s)$,  where $x_i = v_{\pm}$ corresponds to the circle $C_i$, $i = 1,...,\gamma(s)$,
%Let deg$(x_i)=$$+1$(respectively $-1$) if $x_i=x_+$ (respectively, $x_i=x_-$).
Define

\begin{equation}
\text{gr}_{\mathfrak H}((D_L(u), x))=
%\displaystyle\mathfrak h(x) =
\sum^{\nu}_{i=1} %{\gamma(s)}_{i=1}
{\text{deg}}(x_i)\cdot[C_i] \in\mathfrak H,
\end{equation}

\h
where $\deg(x_\pm)=\pm 1$.
\end{defn}

\begin{note}\label{muki}
If $L$ has only one component,
%when we define
then
$n$, $n_+$, $n_-$, $\text{gr}_h$, $\text{gr}_q$, and $\text{gr}_\mathfrak H$,
we do not %need to use
depend on
the orientation of $L$.
If $L$ has greater than one component,
we use that of $L$ when we define $n_+$ and $n_-$.
Note that, if we change the orientation of $L$ into
the opposite  one,
then
neither $n_+$ nor $n_-$
%$n$, $n_+$ $n_-$, $\text{gr}_h$, $\text{gr}_q$, nor $\text{gr}_\mathfrak H$
changes.

We suppose $[\gamma]=[-\gamma]$ when we define $\mathfrak H$.
Hence, whichever orientation we give to a disjoint union of circles $Z(D)$ in
a labeled resolution configuration $(D,x)$, $\text{gr}_\mathfrak H(D,x)$ is the same.

Note the following facts.
Let $\gamma$ and $\gamma'$ be circles embedded in a surface $F$,
and $\gamma\cap\gamma'=\emptyset$.
Let $A$ be an arc which connects  $\gamma$ and $\gamma'$.
Assume that we do not consider the orientation of $\gamma$ nor that of $\gamma'$.
Let $\sigma$ be an embedded circle in $F$
which is obtained
from  $\gamma$ and $\gamma'$
by the surgery along $A$.
Suppose that we do not consider the orientation of $\sigma$.
Once we are given $A$, $\gamma$ and $\gamma'$,
then $\sigma$ is determined, and henceforth
$[\sigma]\in\mathfrak H$ is determined.

Let $\gamma$ be a circle embedded in a surface $F$.
Let $A$ be a c-arc both of whose endpoints are in $\gamma$.
Assume that we do not consider the orientation of $\gamma$.
Let $\mathfrak s$
%$\sigma\amalg\sigma'$
be a disjoint union of two embedded circles in $F$
which is obtained from $\gamma$
by the surgery along $A$.
We do not consider the orientation of $\mathfrak s$.
Note that $\mathfrak s$
determines one or two elements in $\mathfrak H$.
Once we are given $A$ and $\gamma$,
then $\mathfrak s$ is determined.
Furthermore we can know
how many elements of $\mathfrak H$ are determined,
and what the elements are (respectively, the element is).

Let $\gamma$ be a circle embedded in a surface $F$.
Let $A$ be a scs arc both of whose endpoints are in $\gamma$.
%Assume that we do not consider the orientation of $\gamma$.  Let $\sigma$ be an embedded circle in $F$ which is obtained from $\gamma$ by the surgery along $A$. Suppose that we do not consider the orientation of $\sigma$. Then $[\gamma]=[\sigma]\in\mathfrak H$. Furthermore, n
Note that
we do not use a scs arc when we define
the partial order in Definition \ref{2.10}.
%(See also Proposition \ref{ichigo}.)
\\

Each of
$\text{gr}_h(D,x)$,
$\text{gr}_q(D,x)$, and
 $\text{gr}_\mathfrak H(D,x)$
is independent of
which orientation we give $Z(D)$,
%the disjoint union of circles in $(D,x)$,
and is independent of which we choose $L$ or $-L$.
Here,
%let $-L$ be
 $-L$ is the link made from $L$ by reversing the orientation of $L$.

\end{note}
\bb

By Definition \ref{2.15gr}, we have the following.

\begin{fact}\label{hode}
Let $(E, y)$ and $(D, x)$ be labelled resolution configurations of a link diagram %$L$
in a %thickened %oriented closed
surface.
Let $(E, y)\prec(D, x)$.
Then $(E, y)$ and $(D, x)$ have the same quantum grading and
the  same homotopical grading.
\end{fact}

We have the following.

\begin{fact}\label{qde}
Assume that a single cycle surgery changes
a $($non-labeled$)$ resolution configuration $D_L(u)$ into $D_L(v)$.
Let $A_i$ $($respectively, $A_j$$)$ be a labeled resolution configuration defined on
 $D_L(u)$ $($respectively, $D_L(v)$$)$.
 Then $A_i$ and $A_j$ have different quantum gradings.
\end{fact}

\h{\bf Proof of Fact \ref{qde}.}
Recall the definition of quantum gradings, $\text{gr}_q((D_L(u), x))$, above.
Since $n_+$ and $n_-$ are determined by a given  link diagram,
a single cycle surgery does not change $n_+ - 2n_-$.
By the definition of a (single cycle) surgery and that of $|\hskip2mm|$,
 a single cycle surgery changes the parity of $|u|$.
Since a single cycle surgery does not change the number of circles in
a labelled resolution configuration,
 a single cycle surgery does not change the parity of
$ \sharp\{Z\in Z(D_L(u)) | x(Z) = x_+\}
-\sharp\{Z\in Z(D_L(u)) | x(Z) = x_-\}.$
Therefore a single cycle surgery always changes
the quantum grading $\text{gr}_q((D_L(u), x))$.
\qed\bb

We use the partial order in Definition \ref{2.15gr},
and define
the integral ($\Z$-coefficient) Khovanov chain complex for
a link diagram %$L$
in a thickened
%closed oriented
surface in the following subsection.

\bb
\subsection{
The differential,
and the homotopical Khovanov homology for link diagrams in
surfaces
}\label{dif}\hskip10mm

Our definition is the same as that in \cite{MN}.

\begin{defn}\label{korekos}
 {\bf The differential $\delta$.}
Given an oriented link diagram $L$ with $n$ crossings and an ordering of the crossings in $L$,
the {\it Khovanov chain} complex is defined as follows.
The Khovanov chain group $KC(L)$ is the $\Z$-module freely generated
by labeled resolution configurations of the form $(D_L(u), x)$ for $u\in\{0, 1\}^n$.
(The set of all labeled resolution configurations of $L$ is
a basis of $KC(L)$.)
The differential preserves the quantum grading and the homotopical grading,
increases the homological grading by 1,
and is defined as

\begin{equation}\label{bibun}
{\displaystyle
\delta(D_L(v),y)
=
\sum_{{\text{For all}} (D_L(u),x),{\text{such that}} |u|=|v|+1,
{\text{and such that}} (D_L(v),y)\prec (D_L(u),x)}
(-1)^{s_0(\mathcal C_{u,v}) %+\zeta((D_L(u),x), (D_L(v),y))
}}(D_L(u),x),
\end{equation}

\noindent
where $s_0(C_{u,v})\in\Z_2$ is defined as follows:
if $u$$= (\epsilon_1,..., \epsilon_{i-1}, 1, \epsilon_{i+1}, . . . , \epsilon_n)$ and
$v$\\$=(\epsilon_1,..., \epsilon_{i-1}, 0,
\epsilon_{i+1}, . . . , \epsilon_n)$,
then $s_0(C_{u,v}) = (\epsilon_1+\cdot\cdot\cdot+ \epsilon_{i-1})$.
%; see also Definition \ref{4.5}.

%$\zeta((D_L(u),x), (D_L(v),y))$ is defined in  Definition \ref{korekos}.(3) by using the number $\xi$ defined in Definition \ref{korekos}.(2).
\end{defn}

%\h{\bf Note.}  See \S\ref{konnyaku}  There is a difference between the definition of the sign in the differential of Khovanov homology for links in thickened surfaces, defined above,  and that for virtual links in \cite{Man}.

%By Definition \ref{korekos}, we have the following. \begin{prop}\label{rerenore} (D_L(v),y)\prec (D_L(u),x) Let $(D_L(v),y)$ and $(D_L(u),x)$ be labeled resolution configurations. Let $y\vert$ is a natural labeling on  $D_L(v)-D_L(u)$. Let $x\vert$ is a natural labeling on  $D_L(u)-D_L(v)$. Assume that $(D_L(v)-D_L(u), y|)$ is one in $(\quad)$ of the left hand side of one of the identity in Figure  \ref{resolold}, and $(D_L(u)-D_L(v), x|)$ is that of the right hand side.

\begin{exa}\label{rerenore}
In Figures \ref{resolreiC} and \ref{resolreiNC},
examples of how the differential acts.
Note: Let $\bf x$ be a labeled resolution configuration
such that $\bf x$ has two circles and one arc which connects the two circles.
If $\bf x$ is different from the left hand side of all identities
in this figure and Figure \ref{resolreiNC},
then $\delta \bf x=0.$
Let $\bf y$ be a labeled resolution configuration
such that $\bf y$ has one circle and one arc.
If $\bf y$ is different from the left hand side of all identities
in this figure and Figure \ref{resolreiNC},
then $\delta \bf y=0.$
%The labeled resolution configuration out $(\hskip3mm)$ in the right hand side is the same as that in the right hand side for each identity. $\delta'=(-1)^{s_0(\mathcal C_{u,v})}\delta$ in the right hand side of the identity in Definition \ref{korekos}.
%In fact, by using only these rules, we can describe all cases of Definition \ref{korekos}.
%
%
%Note:  Let $\bf x$ be a labeled resolution configuration such that $\bf x$ has two circles and one arc which connects the two circles. If $\bf x$ is different from the left hand side of all identities in Figures \ref{resolreiC} and \ref{resolreiNC} %Figure  \ref{resolold},   then $\delta \bf x=0.$  That is, the case is that both circles equip $x_-$.

%Let $\bf y$ be a labeled resolution configuration such that $\bf y$ has one circle and one arc.  If $\bf y$ is different from the left hand side of all identities in Figures \ref{resolreiC} and \ref{resolreiNC} %Figure  \ref{resolold},   then $\delta \bf y=0.$

\begin{figure}
\vskip-10mm
\includegraphics[width=160mm]{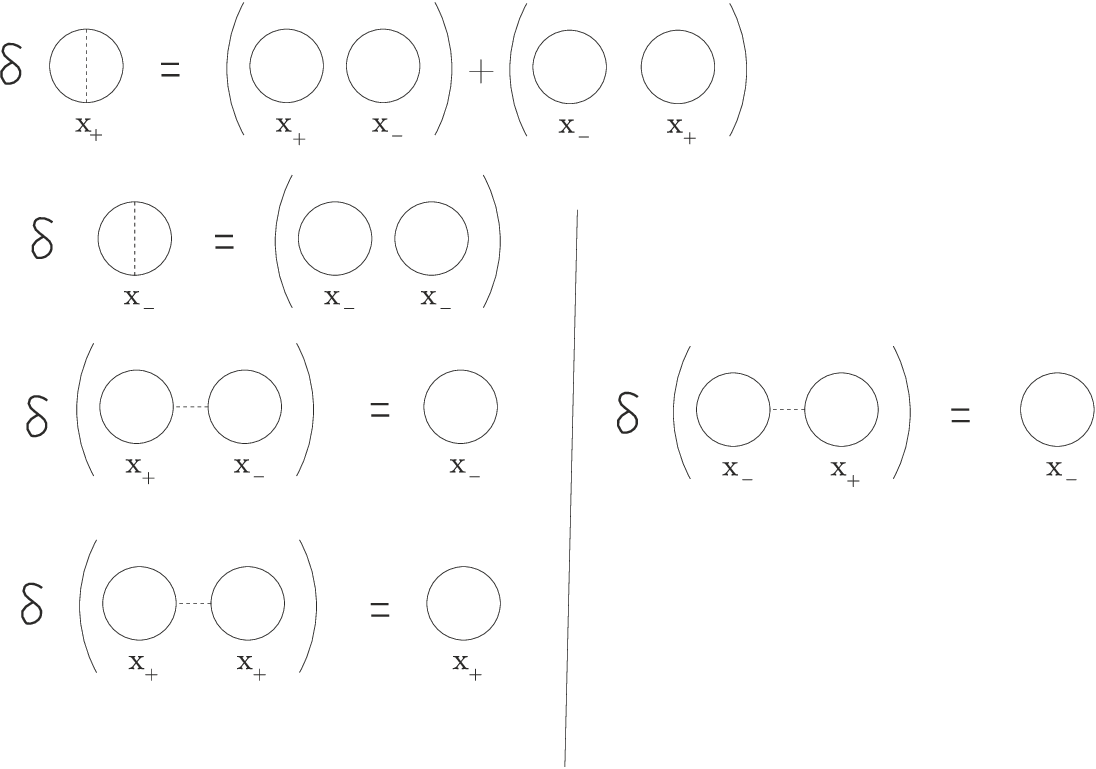}
\caption{{\bf
Examples of how the differential acts. See other examples in Figure \ref{resolreiNC}.
%The labeled resolution configuration out  $(\quad)$  in the right hand side is the same as that in the right hand side for each identity. $\delta'=(-1)^{s_0(\mathcal C_{u,v})}\delta$.  $(-1)^{s_0(\mathcal C_{u,v})}$ is in the right hand side of the identity in Definition \ref{korekos}.
%In fact, by using only these rules, we can describe all cases of Definition \ref{korekos}.
}\label{resolreiC}}
\end{figure}

\begin{figure}
\vskip-10mm
\includegraphics[width=160mm]{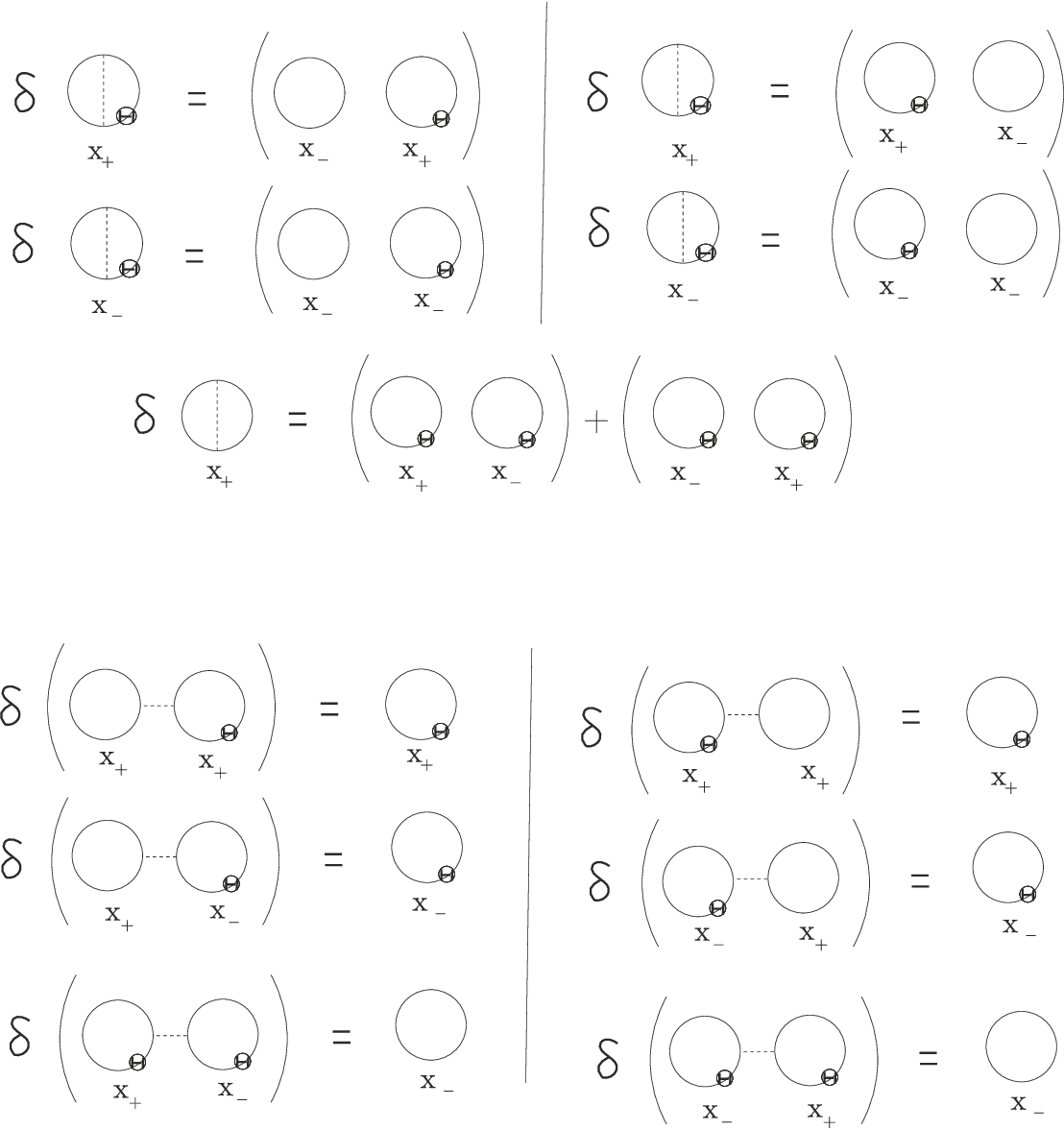}
\caption{{\bf
Examples of how the differential acts. See other examples in Figure \ref{resolreiC}.
%Read the note in the caption of Figure \ref{resolreiC}.
}\label{resolreiNC}}
\end{figure}
\end{exa}

%%%%%

%\bb  We prepare a notation $c[A_i;A_j]$. A standard basis of the Khovanov chain complex  is the set of all labeled resolution configurations made from a link diagram $L$ in a thickened surface.  Let $\{A_i\}_{i\in\Lambda}$ be the set of all labeled resolution configurations made from $L$. Note that $\Lambda$ is a finite set. We use the following notaion \begin{equation}\label{keisu} \delta A_i = \displaystyle\sum_{j\in\Lambda} c[A_i;A_j]\cdot A_j, \end{equation}

%\h and $c[A_i;A_j]$ is an integer coefficient.  \\

\begin{note}\label{daidai}
Recall that
 $\{A_i\}_{i\in\Lambda}$ denotes the set of all labeled resolution configurations made from
an arbitrary link diagram $L$ in a surface.
We use the following  notation
 \begin{equation}\label{keisu2}
\delta A_i = \displaystyle\sum_{j\in\Lambda}
[A_i;A_j]\cdot A_j,
%c[A_i;A_j]\cdot A_j,
\end{equation}

\h like the identity (\ref{keisu}). Here,
$[A_i;A_j]$
%$c[A_i;A_j]$
is an integral coefficient.
Note that,
by Definition \ref{korekos},
$[A_i;A_j]\in\{-1,0,1\}$.
%$c[A_i;A_j]\in\{-1,0,1\}$.

By Definition \ref{korekos},
 the following facts hold.

 %Recall the notation $c[A_i;A_j]$ in the identity (\ref{keisu}). \\

\noindent(1)
If $A_i$ and $A_j$ have different quantum gradings, then
$[A_i;A_j]=0$.
%$c[A_i;A_j]=0$.
(This is related to Fact \ref{qde}.)
Note: This condition holds in the case of the Khovanov homology for  links in $S^3$.
%If $A_i$ and $A_j$ are given as in
%Fact \ref{qde} above, then $c[A_i;A_j]=0$.
% because of Fact \ref{qde}.
%We will explain why we want this condition in Note \ref{point}.
\\

\noindent(2)
If
%(the homological grading $A_i)+1\neq$ (that of $A_j$),
$\text{gr}_h(A_i)+1\neq\text{gr}_h(A_j)$,
then
$[A_i;A_j]=0$.
%$c[A_i;A_j]=0$.
\\

\noindent(3)
$[A_i;A_j]\neq0$
%$c[A_i;A_j]\neq0$
only if
 %(the quantum grading $A_i)=$ (that of $A_j$),
$\text{gr}_q(A_i)=\text{gr}_q(A_j)$,
%  (the homotopical grading $A_i)=$ (that of $A_j$), and
$\text{gr}_{\mathfrak H}(A_i)=\text{gr}_{\mathfrak H}(A_j)$, and

%(the homological grading $A_i)+1=$ (that of $A_j$).

  $\text{gr}_h(A_i)+1=\text{gr}_h(A_j)$.

\bs
 Note that the above facts (1) and (2) follows from this fact (3).
\\

\noindent(4)
If $\delta A_i=0$,  then $A_i$ is a maximal element in  $\{A_i\}_{i\in\Lambda}$.

Therefore we have the following.
Let $A$ (respectively $A'$) be a labeled resolution configuration.
Assume that
a (non-labeled) resolution configuration under $A$ is
changed into that under $A'$ by just one surgery.
Suppose that $\delta A=0$.
Then $A$ and $A'$ are not related by $\prec$.

\end{note}
\bb

By Definition \ref{korekos}, we have the following.
We use notaions in the identities (\ref{bibun}) and (\ref{keisu2}).

\begin{pr}\label{ittoku}
Let $(D_L(v),y)$ and $(D_L(u),x)$ be labelled resolution configurations.

\bs\h$(1)$
We have $[(D_L(v),y); (D_L(u),x)]\neq0$
only if we have the following:
$(D_L(v)-D_L(u),y|)$
is the left hand side of one of the relations by $\prec$ in the tables in Definition \ref{2.10}.  \\
$(D_L(u)-D_L(v),x|)$
is the right hand side of the above relation by $\prec$.
%$($See Definition $\ref{2.10}.(2)$ for the definition of  $x|$ and $y|.)$

\bs\h$(2)$
If $[(D_L(v),y); (D_L(u),x)]\neq0$,
then
 $[(D_L(v),y); (D_L(u),x)]=(-1)^{s_0(\mathcal C_{u,v})}$.
\end{pr}
%Each chain has three kinds of degree: the homological degree, the quantum degree, and the homotopical degree.The differential changes only homological degree: \\
%(The homological degree of $\delta x)=($the homological degree of $ x)+1$.

%Our rule of resolutions is drawn in Figure \ref{resol}. Note that, in each of the lower eight resolutions in Figure \ref{resol}, the labeled resolution configuration with the notation $(H)$ in the left hand side of the arrow  that in the right hand side have the same quantum and homotopical degree.

%The Khovanov-Manturov-Nikonov chain complex of links in thickened surfaces
%in \cite{MN} is defined by the $\Z$ coefficient.

By
Fact \ref{qde} and Note \ref{daidai}.(1), we have the following.

\begin{pr}\label{ichigo}
Assume that a single cycle surgery changes
a $($non-labeled$)$ resolution configuration $D_L(u)$ into $D_L(v)$.
Let $A_i$ $($respectively, $A_j$$)$ be a labeled resolution configuration defined on
 $D_L(u)$ $($respectively, $D_L(v)$$)$.
 Then
 $[A_i;A_j]=0$.
%$c[A_i;A_j]=0$.
\end{pr}

%\bb
\subsection{
The well-definedness of the homotopical Khovanov homology for links in
thickened surfaces
}%\label{dif}
\hskip10mm

%\h We have the following.

\begin{thm}\label{sikon}{\bf (\cite{MN}.)}
Let $L$ be a link diagram of a link $\mathcal L$ in a thickened surface.
For $\delta$ in Definition $\ref{korekos}$,
\begin{equation}
 \delta^2=0.
\end{equation}

%\h That is, the homotopical  Khovanov
 %homology of $L$ is a $\Z$ coefficient homology.
\end{thm}

\h{\bf Note.}
(1) %The coefficient of the above homology is the integer.
The above homology uses integer coefficients.

\bs
\h
(2) By Note \ref{muki},
this homology of $L$ is the same as that of $-L$.
This homology is independent of
which orientation on the disjoint union of all circles in each labeled resolution configuration
we choose.
\\

%Make $G-G'$. Suppose that, in  $G-G'$, there is an arc made from $R$ in $G$. Then this arc in   $G-G'$ is called $R$ again  when it is clear from the context. Assume that we carry out a surgery on $G$ along an arc different from $R$ and we get $G''$. In $G''$ there is an arc made from $R$ in $G$.  This arc in $G''$ is called  $R$ again. % when it is clear from the context.
%
%
Let $G$ and $G'$ be resolution configurations.
Let $R$ be an arc in $G$.
 By definition,
 $A(G-G')$ and $A(s_A(G))$ are subsets of $A(G)$.
Assume that
$A(G-G')$ or $A(s_A(G))$ includes an arc which was $R$ in $A(G)$.
 Then we also call this arc in $A(G-G')$ or $A(s_A(G))$, $R$.
% If $R$ belongs to these subsets,
 %then %$R$ does not change.
% Thus, we can omit this remark.

\bb
\h{\bf Proof of Theorem \ref{sikon}.}
Let $(D,x)$ be any labeled resolution configuration of  $L$.
Take two arbitrary arcs, $A$ and $A'$, in $D$.
Let $x|$ be labelings on $D-s_{A,A'}(D)$ induced by $x$.
If $\delta^2(D-s_{A,A'}(D), x|)=0$,
then  Theorem \ref{sikon} holds.
%Note that $(D-s_{A,A'}(D), x|)$ has only two arcs.
%Note that  $D-s_{A,A'}(D)$ is basic.
Let $(E,z)=(D-s_{A,A'}(D), x|)$.

%We call the arc in $E$ which is derived form $A$ (respectively, $A'$) in $D$,  $A$ (respectively, $A'$) again.
Note that $E$ has only two arcs $A$ and $A'$.
We check all cases of $(E,z)$ and prove $\delta\cdot\delta(E,z)=0$ below.

The union of all arcs  $A$ and $A'$ and all circles in $E$ is a topological space, say $P$.
If $A$ and $A'$ are included in different connected components of $P$,
then the proof is easy.
Note that $E$ is a basic resolution configuration.
Assume that
$P$ above %the union of all arcs and circles in  $E$
is connected.
\bb

\begin{figure}%[t]
\bb
\includegraphics[width=70mm]{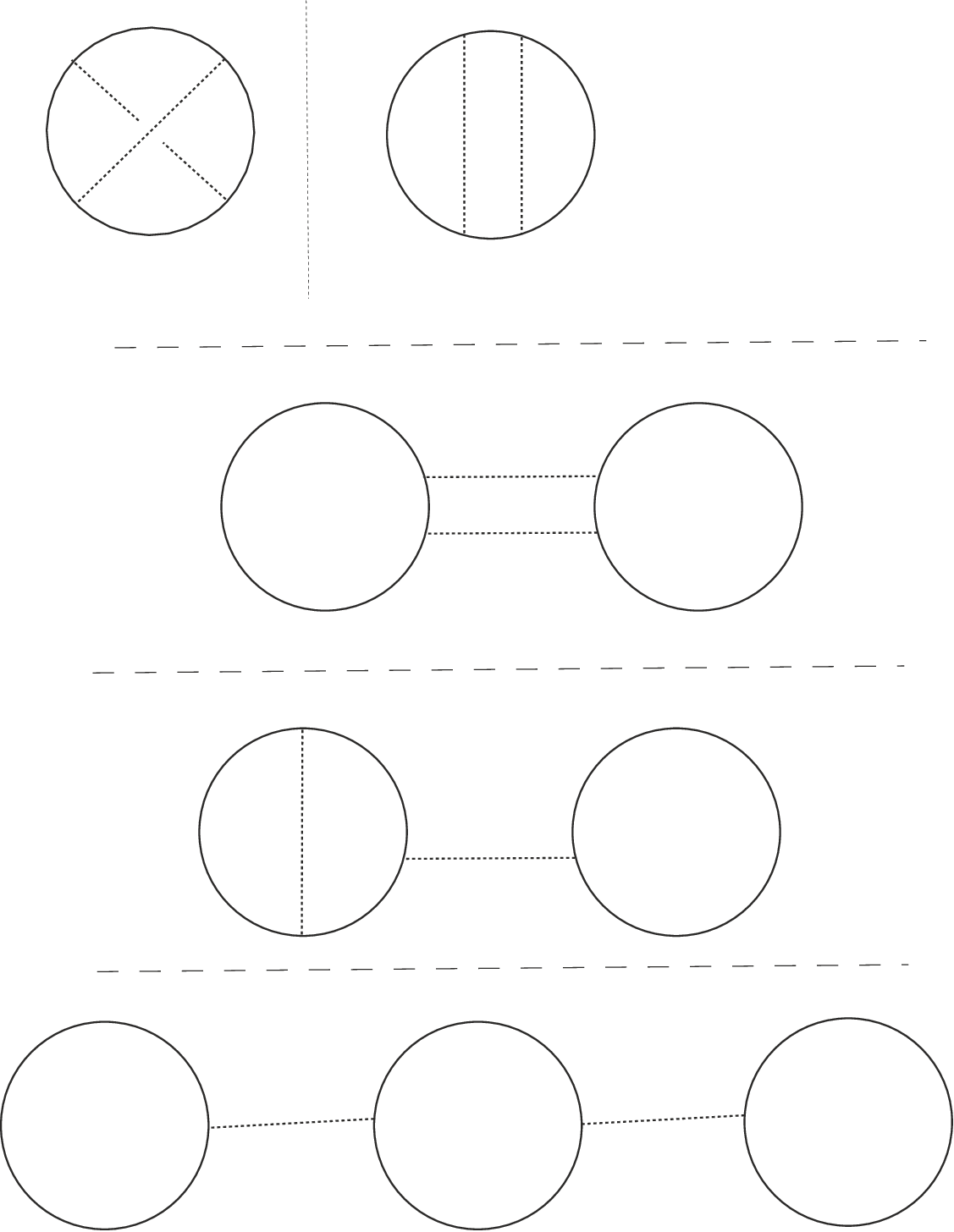}
\caption{{\bf
Connected graphs of the resolution configurations of index 2.
%We draw all cases of $Z(E)\cup A(E)$ in Proof of of Theorem \ref{sikon} abstractly.
}\label{2arc}}
\bb
\end{figure}

In Figure \ref{2arc},  we draw all cases of $Z(E)\cup A(E)$ abstractly.

We check all $(E,z)$ as follows.

\bb
\h{\bf Case 1. }  %これを最後に持ってって ladybugとかをくわしくするか
If only two arcs, $A$ and $A'$, in $E$,
only one arc $A'$ in $s_{A}(E)$, and
only one arc $A$ in $s_{A'}(E)$ are mc arcs,
we have
$\delta^2(E,z)=0$ by applying the rule in
Figures $\ref{resolC}$ and $\ref{resolNC}$ of Definition \ref{2.10}.

\bb
\h{\bf Note.}
Let $y$ be a labeling on $s(E)$.
%Assume $(E,y,z)\neq\emptyset.$
%$\{p|p$ is a labeled resolution configuration. \\
%\hskip20mm $(E,z)\prec p\prec(s(E),y)$, $p\neq(E,z), p\neq(s(E),y)\}\neq\emptyset.$
%
%
%$\{p|p$ is a labeled resolution configuration. $p$ is obtained form $(E,z)$ by one surgery. $\}\neq\emptyset.$
%Then the number of the elements of this set is 4 or 6.
There is a case where
the number of the elements in $P(E,y,z)$ %$(E,y,z)$
is six.
%
%
%The case of 4 includes the ladybug configuration in \cite[\S 5.4.2]{LSk}. It is important that the case of 4 in the case of links in thickened surfaces includes a different feature from the ladybug configuration. See  Note \ref{insight}.
This case is important. We will discuss it in \S\ref{tamago} and \S\ref{hiyoko}.
\bb

\h{\bf Case 2. }
If only two arcs, $A$ and $A'$, in $E$ are scs arcs,
%each of $A$ and $A'$ induces a single cycle surgery,
$\delta(E,z)=0$ %(D-s_{A,A'}(D), x|)=0$
by Proposition \ref{ichigo}. %Note \ref{daidai}.(1).
Hence $\delta^2(E,z)%(D-s_{A,A'}(D), x|)
=0$.

\bb

\begin{figure}
\includegraphics[width=90mm]{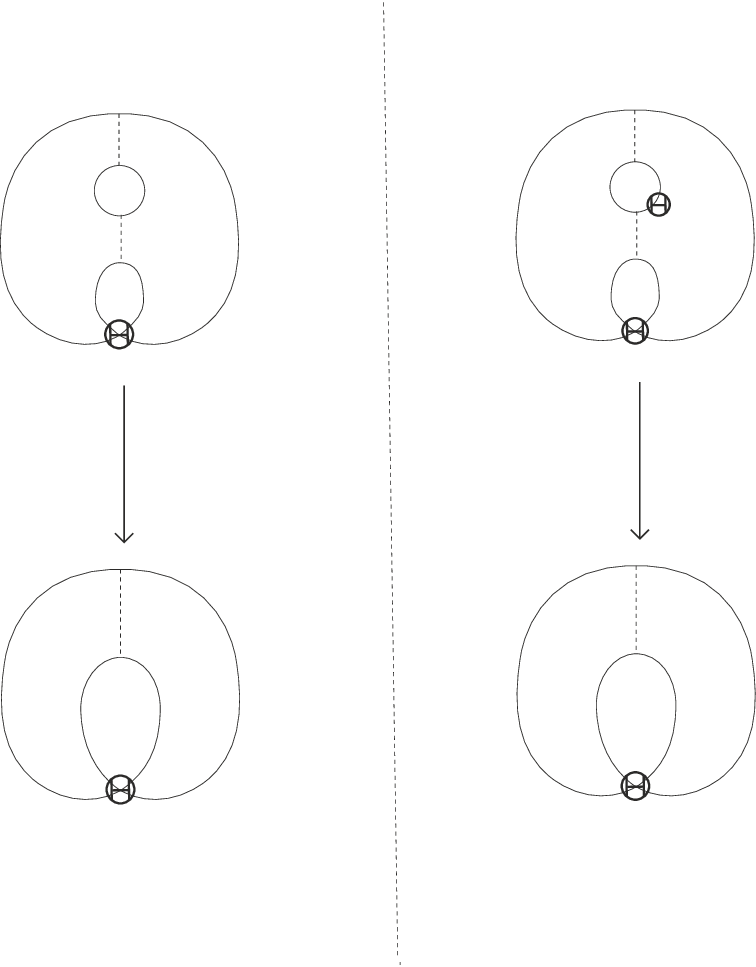}
\caption{{\bf
Surgery along $A$ in the case when $A$ and $A'$ in $E$ are mc arcs,
and $A'$ in $s_A(E)$ is a scs arc.
%
%Each side is an example of surgery: Only two arcs, $A$ and $A'$, in $E$ are mc arcs and only one arc in $s_{A}(E)$ is a scs arc.
}\label{nanaba}}   \bb\end{figure}

\bb
\h{\bf Case 3. }
Suppose that only two arcs, $A$ and $A'$, in $E$ are mc arcs, and
that only one arc $A'$ in $s_{A}(E)$ is a scs arc.
%Note that there is a case such that both $A$ and $A'$ are mc arcs and such that $s_{A}(E)$ has a scs arc.
(See two examples in Figure \ref{nanaba}.)
Then we have the following:
The surgery from $s_{A}(E)$ to $s_{A,A'}(E)$ is a  single cycle surgery.
%The difference of the number of circles in $E$ and that in $s_{A,A'}(E)$ is one.
Hence the number of circles in $E$ has a different parity of that in $s_{A,A'}(E)$.
Therefore
%one of
the surgery %from $E$ to $s_{A'}(E)$ and that
from $s_{A'}(E)$ to $s_{A,A'}(E)$ is a  single cycle surgery.
By Proposition \ref{ichigo}, %Note \ref{daidai}.(1),
we have $\delta^2(E,z)=0$.
%(We can prove this although we do not know the following:  Since the surgery from $E$ to $s_{A'}(E)$ is not a single cycle surgery, that from $s_{A'}(E)$ to $s_{A,A'}(E)$ is a  single cycle surgery.)

Assume that only two arcs, $A$ and $A'$, in $E$ are mc arcs, and that
only one arc $A$ in $s_{A'}(E)$ is a scs arc. We can prove $\delta^2(E,z)=0$ by the same method as the above one.

\bb
\h{\bf Case 4. }
Assume that one of only two arcs, $A$ and $A'$, in $E$ is a mc arc, and
the other a scs arc.
We can let $A$ be a mc arc and $A'$ a scs arc without loss of generality.

%We have the following. There are orientations of circles in $E$ with the following properties: Take a small disc which includes $A'$. Then the orientations of the intersection of the circles and the disc are as drawn in Figure \ref{osa}.O.
%
%By the definition of scs arcs, we have the following. Take a small disc which includes $A$. Then the orientations of the intersection of the circles and the disc is as drawn in Figure \ref{osa}.S or   \ref{osa}.S$'$.
%
%Carry out a surgery along $A'$. After that, let  the orientations of the intersection of the circles and the disc be drawn as in Figure \ref{osa}.G. Do not change  the orientations of circles out the disc. Hence there are orientations of circles in $s_{A}(E)$ with the following properties: the orientations of the intersection of the circles and the disc is as drawn in Figure \ref{osa}.G. (See the relation between Figures \ref{osa}.O and \ref{osa}.G.). Then the orientations of the intersection of the circles and the disc with $A'$is not changed. (They are drawn as in Figure \ref{osa}.S or   \ref{osa}.S$'$.)

%Therefore only one arc $A'$ in $s_{A}(E)$ is a scs arc after the surgery.

%By Proposition \ref{ichigo}, we have $\delta^2(E,z)=0$.

%\begin{figure} \includegraphics[width=120mm]{osa.eps} \caption{{\bf Arcs and small discs which include the arcs }\label{osa}}   \end{figure}

%\bb  \h{\bf Case 5.} Assume that only two arcs, $A$ and $A'$, in $E$ are scs arcs.

Note that, under this setting,  there is a case such that
%only two arcs, $A$ and $A'$, in $E$ are scs arcs and such that
only one arc $A'$ in $s_{A}(E)$ is a mc arc.
The other possibility,
when $A'$ in $s_A(E)$ is a scs arc, % and $A$ in $s_{A'}(E)$ is a mc arc,
is easy.

\begin{figure}
\includegraphics[width=55mm]{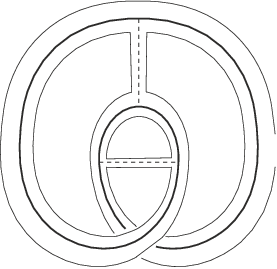}
\caption{{\bf
A shell configuration and its neighborhood
}\label{shell}}
\end{figure}

In order to show such examples, we need to introduce a terminology.

A {\it shell configuration} or {\it shell Kauffman state}
is a resolution configuration of a link diagram in a surface $F$
drawn in Figure \ref{shell}: We draw only a neighborhood $N$ of the shell configuration in $F$.
%We do not draw all $F$.
We assume that $N$ is a compact surface and
 that the inclusion map of the shell configuration to $N$ is a homotopy type equivalence map.

%The left upper resolution configuration $Q$ in each of Figures \ref{Hp} and \ref{Hq} is a shell configuration.

\begin{figure}
\includegraphics[width=90mm]{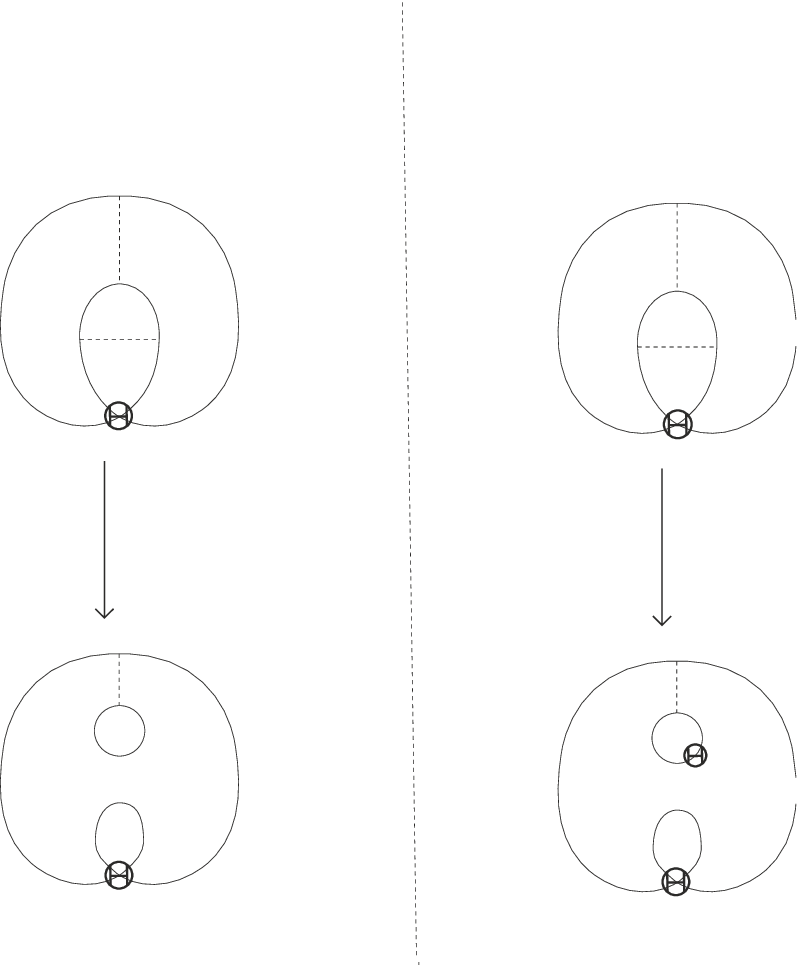}
\caption{{\bf
Each side is an example of surgery:
Only two arcs, $A$ and $A'$, in $E$ are mc arcs
and only one arc in $s_{A}(E)$ is a scs arc.
}\label{banana}}
\end{figure}

Two examples of the above case are drawn in Figure \ref{banana}:
The upper of each side is a shell configuration.

\begin{figure}
\includegraphics[width=100mm]{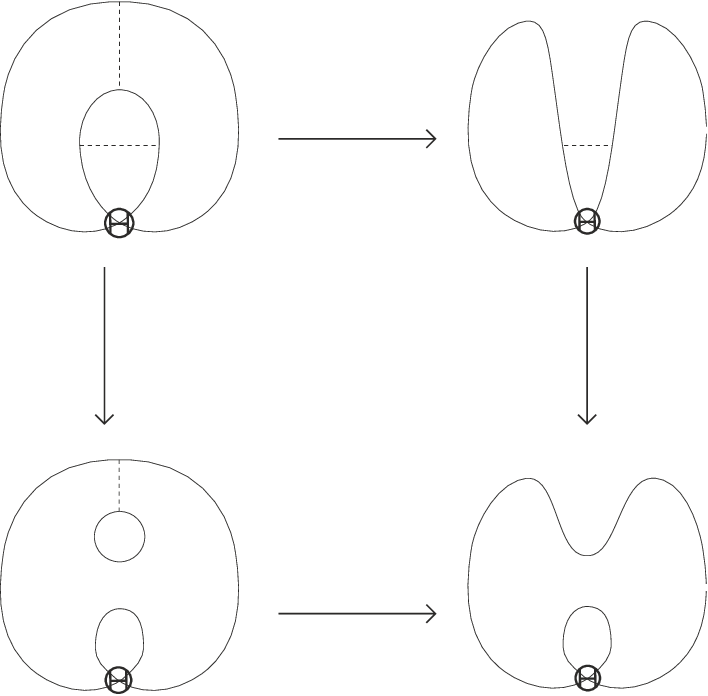}
\caption{{\bf A link diagram $L$ in a thickened
%closed oriented
surface, and
the relation,  which is made by surgeries,
among all (non-labeled) resolution configurations made from $L$.
One circle in the left lower labeled resolution configurations
is non-contractible, and the other contractible.
}\label{Hp}}
\end{figure}

\begin{figure}
\includegraphics[width=100mm]{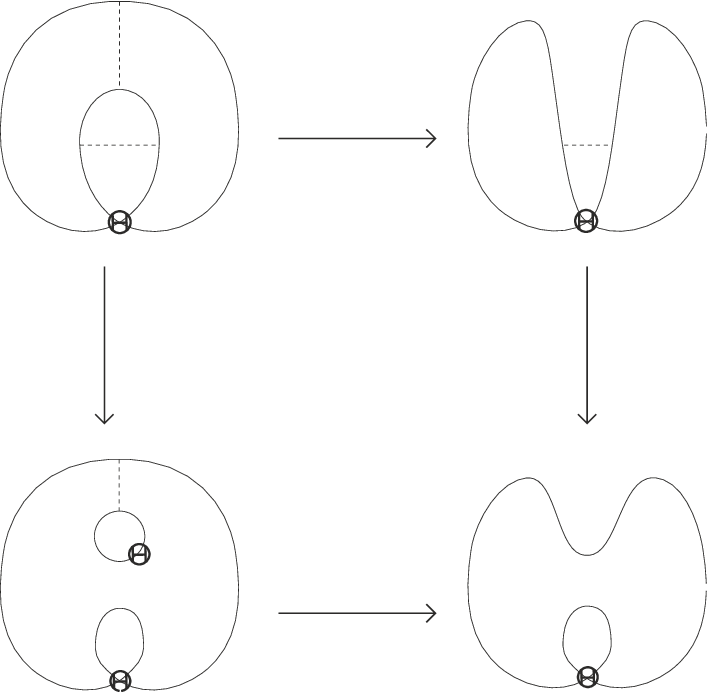}
\caption{{\bf A link diagram $L$ in a thickened surface, and
the relation,  which is made by surgeries,
among all (non-labeled) resolution configurations made from $L$.
The two circles in the left lower labeled resolution configurations
are non-contractible.
}}\label{Hq}
\end{figure}

We must take care of the cases in Figures \ref{Hp} and \ref{Hq},
which are associated with Figure \ref{banana}.
 The other cases are easy. %, or are proved by the same method as that in these cases.

%We may take care of the cases in Figure \ref{nanaba} % \ref{S1X} but it is not so difficult.

\begin{figure}
\includegraphics[width=90mm]{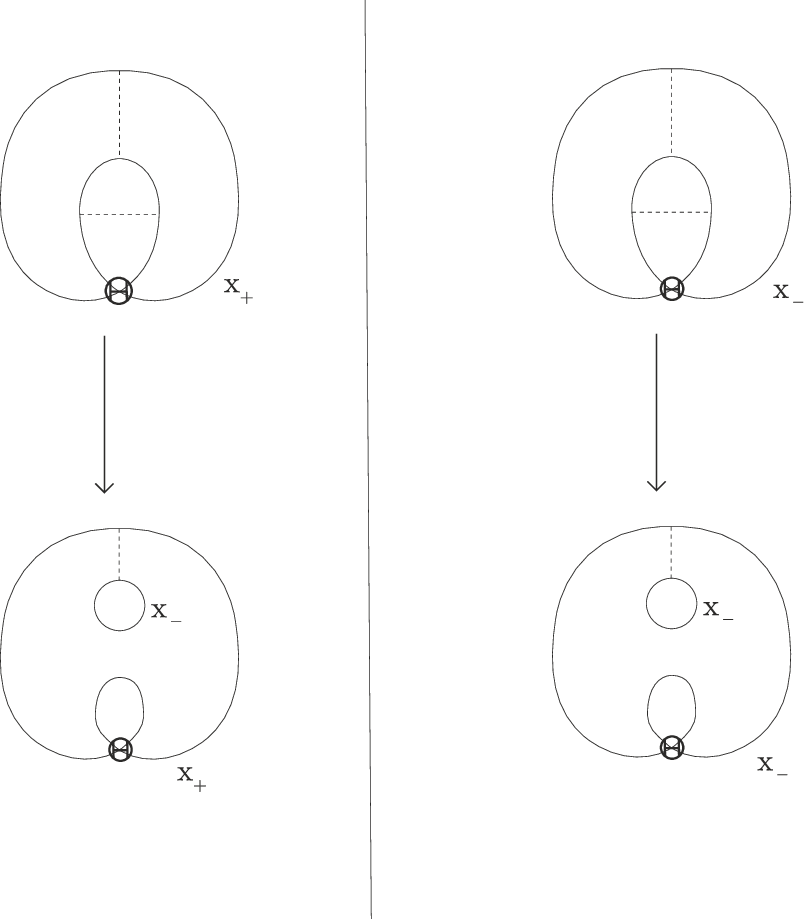}
\caption{{\bf
We put a labbeling $x_+$ (respectively, $x_-$)
on the left (respectively, right) upper resolution configuration
 in Figure \ref{Hp}.
 We check what labelled resolution configurations are made from it
 by the differential.
}\label{H+}}  %
\end{figure}

In the case in Figure \ref{Hp}, we have $\delta\cdot\delta=0$
because of Figure \ref{H+}.

\begin{figure}
\includegraphics[width=90mm]{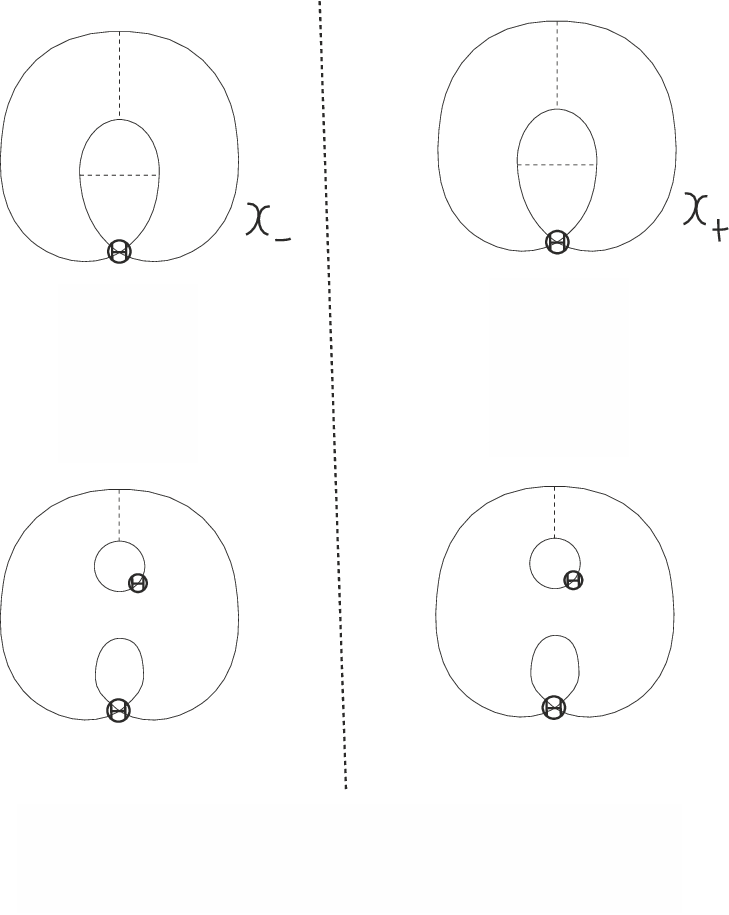}
\vskip-12mm
\caption{{\bf
In each side the lower (non-labelled)  resolution configuration
is obtained from that under the above labeled resolution configuration
 by a surgery along a m-arc.
By Note \ref{daidai}.(4),
we have the following:
Whichever of $x_+$ and $x_-$ we give the lower resolution configuration,
the upper and the lower are not related by $\prec$.
}\label{HqZ}}
\end{figure}

In the case in Figure \ref{Hq}, we have $\delta\cdot\delta=0$
because of Figure \ref{HqZ}.

\bb
Therefore $\delta^2(E,z)
=0$ in all cases.
This completes the proof of Theorem \ref{sikon}. \qed
\bb

\h{\bf Note.}
When we consider this case, we should take care of  the following facts.

Let $(E,z)$ be a labeled resolution configuration
which is obtained by giving a labeling to  the left upper
shell configuration in Figure \ref{Hp} (respectively, \ref{Hq}).
Then there is only one element in the set
$\{p|p$ is a labeled resolution configuration. $(E,z)\prec p,  (E,z)\neq p\}.$

   In the identities in  Figures \ref{resolreiC} and \ref{resolreiNC},
only
the most upper identity in Figure \ref{resolreiC}
and
the third identity in Figure \ref{resolreiNC}
have two labeled resolution configurations in the right hand side.
In both cases, the left hand side has
only one labeled resolution configuration and
the circle in it is a contractible circle.

Shell configurations,
Figures \ref{Hp} and \ref{Hq} are important. %See \S\ref{konnyaku}.
\\

%\begin{figure} \includegraphics[width=80mm]{S1X.eps}\caption{{\bf }\label{S1X}}   \end{figure}

\bb
By using the same method in \cite{B,K},
we have the following.

\begin{thm}\label{daiji}{\bf (\cite{MN}.)}
Let $L$ and $L'$ be link diagrams of a link $\mathcal L$ in a thickened surface.
Then the homotopical  Khovanov %-Manturov-Nikonov
 homology of $L$ and
that of $L'$ are equivalent.
\end{thm}

By Theorem \ref{daiji}, the following definition  is well-defined.

\begin{defn}\label{budo}
Let $\mathcal L$ be a link in a thickened surface.
Let $L$ be a link diagram of a link $\mathcal L$.
We define the homotopical  Khovanov %-Manturov-Nikonov
 homology of  $\mathcal L$   to be that of $L$.
 \end{defn}

\section{\bf The ladybug configuration for link diagrams in $S^2$}\label{tamago}

\noindent
We review the ladybug configuration for link diagrams  in $S^2$,
which is introduced in \cite[section 5.4]{LSk}.
Lipshitz and Sarkar introduced it to define a CW complex for any link diagram  in $S^2$.
We cite the definition of it,
that of the right pair, and that of the left pair associated with it
from \cite[section 5.4.2]{LSk}. \\

\begin{figure}
\includegraphics[width=40mm]{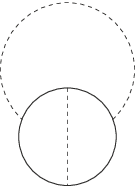}
\caption{{\bf The ladybug configuration
}\label{tento}}
\end{figure}

\begin{defn}\label{teten}  {\bf (\cite[Definition 5.6]{LSk}).}
An index 2 basic resolution configuration $D$ in $S^2$
is said to be a ladybug configuration
if the following conditions are satisfied (See Figure \ref{tento}.).

$\bullet$ $Z(D)$ consists of a single circle, which we will abbreviate as $Z$;

$\bullet$ The endpoints of the two arcs in $A(D)$, say $A_1$ and $A_2$,
alternate around $Z$

\hskip3mm (that is, $\partial A_1$ and $\partial A_2$ are linked in $Z$).\\

\end{defn}

\begin{defn}\label{rl}  {\bf (\cite[section 5.4.2]{LSk}).}
Let $D$ be as above.
Let $Z$ denote the unique circle in $Z(D)$.
The surgery $s_{A_1}(D)$ (respectively,  $s_{A_2}(D)$) consists of two circles;
denote these $Z_{1,1}$ and $Z_{1,2}$ (respectively, $Z_{2,1}$ and $Z_{2,2}$);
that is, $Z(s_{A_i}(D)) = \{Z_{i,1}, Z_{i,2}\}$.
Our main goal is to find a bijection between
$\{Z_{1,1}, Z_{1,2}\}$  and $\{Z_{2,1}, Z_{2,2}\}$;
this bijection will then tell us which points in
$\partial_{\rm exp}\mathcal M(x, y)$ to identify.
See \cite[(RM-2) in section 5.1]{LSk} for the notation $\partial_{\rm exp}$.

As an intermediate step, we distinguish two of the four arcs in
$Z - (\partial A_1\cup \partial A_2)$.
Assume that the point $\infty\in S^2$ is not in $D$,
and view $D$ as lying in the plane $S^2-\{\infty\}\cong\R^2$.
Then one of $A_1$ or $A_2$ lies outside $Z$ (in the plane)
while the other lies inside $Z$.
Let $A_i$ be the inside arc and $A_o$ the outside arc.
The circle $Z$ inherits an orientation from the disk it bounds in $\R^2$.
With respect to this orientation, each component of
$Z - (\partial A_1\cup\partial A_2)$
either runs from the outside arc $A_o$ to an inside arc $A_i$ or vice-versa.
The {\it right pair} is the pair of components of $Z-(\partial A_1\cup\partial A_2)$
which run from the outside arc $A_o$ to the inside arc $A_i$.
The other pair of components is the {\it left pair}. See \cite[Figure 5.1]{LSk}.
\end{defn}

We explain why the ladybug configuration is important:
By each ladybug configuration,
there are many possibilities of homotopy type of Khovanov CW complex.
So we must check whether such homotopy types are  stable homotopy type equivalent or not.
See below and Fact \ref{konnya}.

\begin{pr}\label{4}
Let ${\bf x}$ $($respectively, ${\bf y})$ be a labelled resolution configuration in $S^2$
of
homological grading $n$ $($respectively, $n+2).$
 Then the cardinality of the set

\hskip3mm $\{p| p$ is a labelled resolution configuration.
${\bf x}\prec p, p\prec{\bf y},$ $p\neq{\bf x}$, $p\neq{\bf y}\}$

\noindent
is 0, 2, or 4, where $\prec$ is defined in Definition \ref{2.10}.
  \end{pr}

\begin{figure}
\includegraphics[width=155mm]{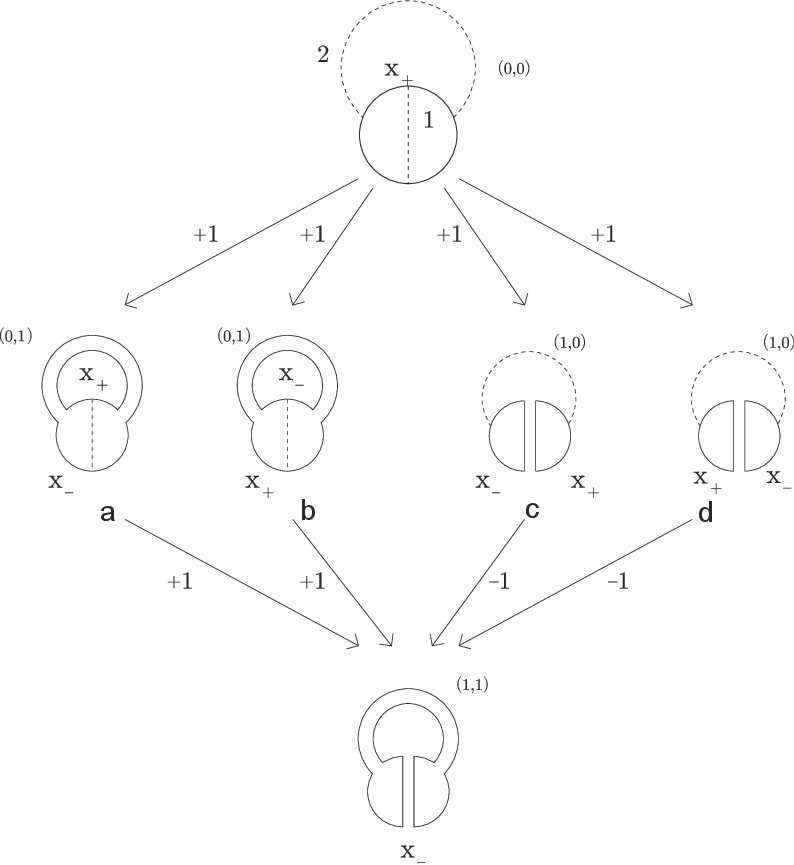}
\caption{
                 {\bf  The poset for
                 the decorated resolution configuration associated with
                                       a ladybug  configuration in $S^2$
                  }\label{uenp}
             }
\end{figure}

Let $D$ be the ladybug configuration  in $S^2$.
Since each of $D$ and $s(D)$ has only one circle,
we can let $x_+$ or $x_-$ denote a labeling on it.
Give $D$ (respectively, $s(D)$) a labeling $x_+$ (respectively, $x_-$).
We call the resultant labeled resolution configuration
 $(D,x_+)$ (respectively, $(s(D),x_-)$).
We obtain a decorated resolution configuration $(D,x_-,x_+)$
as drawn in Figure \ref{uenp}.
%This decorated resolution configuration is called the {\it decorated resolution configuration associated with the ladybug configuration $D$}.

\begin{fact}\label{tenten}
The case of 4 in Proposition \ref{4} occurs in the above case $(D,x_-,x_+)$.
%when we have the decorated resolution configuration associated with the ladybug configuration.
\end{fact}

Fact \ref{tenten}  is also explained in \cite[section 5.4]{LSk}.

\bigbreak
\section{\bf  Khovanov-Lipshitz-Sarkar stable homotopy type
and Steenrod square for  links  in $S^3$}\label{Khpi}

\noindent
The natural dual of the Khovanov chain group is called the {\it dual Khovanov chain group}.
The natural dual of the Khovanov basis % $($respectively, a Khovanov basis element$)$
is called the {\it dual Khovanov basis}. %$($respectively, \text{a} Khovanov basis element$)$}.
See \cite[\S3]{KauffmanOgasasq} for the explanation of this kind of dual. % for detail.
\bs

Let $\mathcal L$ be a link  in $S^3$.
Let $L$ be a link diagram  in $S^2$ which represents $\mathcal L$.
In \cite{LSk} Lipshitz and Sarkar
made a consistent method to construct a CW complex for $L$
whose singular homology is the homology of the dual Khovanov chain complex of $\mathcal L$.

Note that, by the definition of the Khovanov homology,
the cohomology, not homology, of Khovanov homotopy type is
the Khovanov homology.

Let $L'$ be a link diagram  in $S^2$ which represents $\mathcal L$.
They also proved that the stable homotopy type of the CW complex for $L$ and that for $L'$ are the same.
We call this stable homotopy type,
Khovanov-Lipshitz-Sarkar stable homotopy type for $\mathcal L$.
Thus this stable homotopy type has Steenrod squares.
In \cite{LSs} Lipshitz and Sarkar found how to calculate the second Steenrod square
associated with Khovanov-Lipshitz-Sarkar stable homotopy type
by using labeled resolution configurations.
In \cite{Seed}
Seed proved the following fact by making a computer program
of Lipshitz and Sarkar's calculation of the second Steenrod square:
There is a pair of links  in $S^3$ with the following properties.
Their second Steenrod squares and
their Khovanov-Lipshitz-Sarkar stable homotopy types
are different.
Their Khovanov homologies are the same.

We review Lipshitz and Sarkar's method below.

\bigbreak

Let $L$ in $S^2$ and $\mathcal L$ in $S^3$ be as above.
In \cite[\S5 and \S6, in particular, Definition 5.3]{LSk},
the Khovanov-Lipshitz-Sarkar stable homotopy type of $L$ and
that of $\mathcal L$
is defined as follows.

Let $(D, x, y)$ be any index $n$ basic decorated resolution configuration.
In \cite[\S5 and \S6]{LSk} Lipshitz and Sarkar
associate to
%each index $n$ basic decorated resolution configuration
$(D, x, y)$
an $(n-1)$-dimensional $<n-1>$-manifold $\mathcal M(D, x, y)$
together with an $(n-1)$-map
$$\mathcal F :\mathcal  M(D, x, y)
\to\mathcal M_{\mathscr C(n)}(\overline1, \overline0).$$

See \cite[\S3.1]{LSk} for $<m>$-manifolds and $m$-maps, where $m$ is an integer.

See \cite[Definition 3.3]{LSk} for the definition of $\bar 0,\bar 1$.

Note that  we regard
the poset for
each index $n$ basic decorated resolution configuration
as a flow category.
See \cite[Definition 3.12]{LSk} for the definition of flow category.
\\

The {\it Khovanov flow category} $\mathcal C_K(L)$ has one object
for each %of the
Khovanov basis element. %s.
%the standard basis elements %generators
%of Khovanov homology, cf. \cite[Definition 2.15]{LSk}.
That is, an object of $\mathcal C_K(L)$ is a labeled resolution configuration of
the form $\mathbf{x}=(D_L(u), x)$  with $u\in\{0, 1\}^n$.
The grading on the objects is
the homological grading gr$_h$; % from \cite[Definition 2.15]{LSk};
the quantum grading gr$_q$ is an additional grading on the objects.
We need the orientation of $L$ in order to define these
gradings, but the rest of the construction of $\mathcal C_K(L)$ is independent of the orientation.
Consider objects $\mathbf{x}=(D_L(u), x)$ and
$\mathbf{y}=(D_L(v), y)$ of $\mathcal C_K(L)$.
The space $\mathcal M_{\mathcal C_K(L)}(\mathbf{x},\mathbf{y})$
is defined to be empty unless $y\prec x$
with respect to the partial order from Definition \ref{2.10}.
So,
assume that $y\prec x$.
Let $x|$ denote the restriction of $x$ to
$s(D_L(v)-D_L(u))=D_L(u)-D_L(v)$
and
let $y|$ denote the restriction of $y$ to $D_L(v)-D_L(u)$.
Therefore, $(D_L(v)-D_L(u), x|, y|)$  is
a basic decorated resolution configuration.
%See \cite[Definition 2.4]{LSk} for the term, ``basic''.
We define $\mathcal M(D_L(v)-D_L(u),x|, y|)$ as above.
Use it, and define
$$\mathcal M_{\mathcal C_K(L)}(\mathbf{x},\mathbf{y})\\=\mathcal M(D_L(v)-D_L(u),x|, y|),$$
as smooth manifolds with corners.
\\

In \cite[\S5]{LSk}, it is proved that,
if gr$_h\mathbf{x}-$gr$_h\mathbf{y}=n$,
$\mathcal M_{\mathcal C_K(L)}(\mathbf{x},\mathbf{y})$
is a disjjoint union of some copies of
the $n$-dimensional cube moduli
$\mathcal M_{\mathscr C(n)}(\overline1, \overline0)$.

\bb

\begin{defn}\label{poy}
%Let $\mathcal L$ be a link  in $S^3$.
%Let $L$ be a link diagram  in $S^2$ which represents $\mathcal L$.
Let $L$ in $S^2$ and $\mathcal L$ in $S^3$ be as above.
We define a CW complex $Y(L)$ for $L$ below.
\bs

Let $\{a_p\}_{p\in\Lambda}$ be the dual Khovanov basis for $L$.
Note that $\Lambda$ is a finite set.
\bs

Fix $n\in\Z$. Let $g^n_i$ be
all dual Khovanov basis elements
whose homological grading is $n$ in $\{a_p\}_{p\in\Lambda}$.
We assign
to
$g^n_i$
a $(n+N)$-cell   $e^{n+N}_i$,
where $N$ is a large integer.  \bs

We attach the cells, $e^{n+N}_i$, for all $n$:
We use the moduli spaces defined above,
with an arbitrary set of framings which satisfy \cite[Definition 3.20]{LSk},
according to the method in \cite[Definition 3.23]{LSk}
(note Proposition \ref{kya} below).
The result is $Y(L)$. \bs

The stable homotopy type of the CW complex $Y(L)$
is called the {\it Khovanov-Lipshitz-Sarkar stable homotopy type}
for the link diagram $L$.
More precisely, since we use an arbitrary large integer $N$ to construct $Y(L)$,
we must handle $N$ as follows.
$N$ times of the formal desuspension of $Y(L)$ is called
the  {\it Khovanov-Lipshitz-Sarkar spectrum} for the link diagram $L$.
\end{defn}

\h
Note:
Following Lipshitz and Sarkar \cite[Definition 5.5]{LSk},
 the Khovanov homology is the reduced cohomology of
 the Khovanov space
 shifted by $(-C)$ for some positive integer $C$.
 The Khovanov spectrum $\chi_{Kh}(L)$ is the suspension spectrum of the Khovanov space, de-suspended $C$ times. Here we can take $C = N$.
\\

\noindent{\bf Note.}
The dual Khovanov chain complex is made from Khovanov chain complex uniquely, and vice versa.
So the following two sentences (1) and (2) have the same meaning.

\bs\h
(1) We associate the framed flow category $\mathcal C$ to Khovanov chain complex.

\bs\h
(2) We associate the framed flow category $\mathcal C$ to
the dual Khovanov chain complex.

\bs\h
(Here, suppose that the above Khovanov chain complex and
the above dual Khovanov chain complex are dual each other.)
\\

When we make $\mathcal F$ above,
it is important to analyze the ladybug configuration \cite[\S5.4.2]{LSk}.

\begin{fact}\label{konnya}{\bf\rm(\cite[\S6]{LSk}.)}
The stable homotopy type of  Khovanov-Lipshitz-Sarkar stable homotopy type
for  link diagrams  in $S^2$
does not depend on
whether we use the right pair or the left pair of each ladybug configuration.
%in the case of links  in $S^3$.
\end{fact}

Fact \ref{ari} is used in the proof of Fact \ref{konnya}.

\begin{fact}\label{ari} {\bf(This is written in \cite[Proof of Proposition 6.5]{LSk}.)}
Fix a classical link diagram $L$ in $S^2$, and let $L'$ be the result of reflecting
$L$ across the $y$-axis, say, and reversing all of the crossings.
Then $L$ and $L'$ represent the same link in $S^3$.
\end{fact}

\begin{pr}\label{kya}  {\bf
(This follows from results in \cite{LSk}. See the comments below.)}
The stable homotopy type of  the Khovanov-Lipshitz-Sarkar stable homotopy type
for link diagrams in $S^2$
does not depend on
the choice of the coherent framing, which is defined in \cite[Definition 3.18]{LSk}.
\end{pr}

%Of course Proposition \ref{kya} does not hold in the general case of construction of CW complexes.  See Example \ref{exam} below. However, Proposition \ref{kya}  is true in this case.
Proposition \ref{kya} is the same as
\cite[(4) in the first part of section six]{LSk}, which is proved in the proof of \cite[Proposition 6.1]{LSk}:
In three lines above \cite[Definition 3.4]{LSs},
it is written, ``all such framings lead to the same Khovanov homotopy type
\cite[Proposition 6.1]{LSk}''.
See also \cite[Lemma 4.13, which is cited in the proof of Proposition 6.1]{LSk}.
In short,
in the case of Khovanov homotopy type for links in $S^3$,
all modulis are contractible so we do not need to check framings on them.

Of course,
the same set of modulis and  different set of framings
construct CW complexes of different stable homotopy types, in general.
 See Example \ref{exam} below.

\begin{exa}\label{exam}
Let $S^2\vee S^4$ denote the one point union of $S^2$ and $S^4$.
Both
$\Sigma^k(S^2\vee S^4)$
and
$\Sigma^k(\C P^2)$,
 where $\Sigma^k$ denotes the $k$-times suspension and $k$ is large,
have a natural CW decomposition
 (the base point)$\cup e^{2+k}\cup e^{4+k}$.
Consider a set of moduli spaces associated with \\
$\Sigma^k(S^2\vee S^4)$
(respectively, $\Sigma^k(\C P^2)$).
In $\partial e^{2+k}$, there is no moduli space.
In $\partial e^{4+k}$, take an embedded circle.
It is a moduli space.
Take the normal bundle of the circle in  $\partial e^{4+k}$,
and take the trivial (respectively, nontrivial) framing.
It is a framing on the moduli space.
 \end{exa}

\begin{defn}\label{todo}
In \cite{LSk} it is proved that
if two link diagrams $L$ and $L'$ in $S^2$ represent the same link in $S^3$,
then $Y(L)$ and $Y(L')$ (see Definition \ref{poy})
are stable homotopy type equivalent.
Thus we obtain a unique stable homotopy type,
{\it Khovanov-Lipshitz-Sarkar stable homotopy type},
and the
{\it Khovanov-Lipshitz-Sarkar spectrum}
for  links in $S^3$.
They are link type invariants.
\end{defn}

It is an outstanding property of
the Khovanov chain complex
and Khovanov stable homotopy type for  link diagrams in $S^2$
 that,
if $\mathcal M_{\mathcal C_K(L)}(\mathbf{x},\mathbf{y})\neq\emptyset$,
each connected component of
$\mathcal M_{\mathcal C_K(L)}(\mathbf{x},\mathbf{y})$
is determined only by
gr$_h\mathbf{x}-$gr$_h\mathbf{y}$.
Chain complexes in other cases do not have this property in general.

\section{\bf Ladybug %configurations
and quasi-ladybug configurations for link diagrams in surfaces}\label{hiyoko}

In this section we explain
why the torus case is special,
 and
 why we concentrate in  the higher genus case in this paper.
We will discuss the torus case in a sequel of this paper \cite{zoku}.

\begin{defn}\label{LQ}
Let $D$ be a resolution configuration
which is made of one circle and two m-arcs.  \\

Stand at a point in the circle where you see an arc to your right. %left.
Go ahead along the circle. Go around one time.
Assume that you encounter the following pattern:
 In the order of travel you next touch the other arc.
Then you touch the  first arc.
Then you touch the other arc again.
Finally, you come back to the point at the beginning. \\

Since both arcs are m-arcs,
both satisfy the following property:
At both endpoints of each arc,
you see the arc in the same side -- either on
the right hand side or on the left hand side.\\

If you see the arcs
both in the right hand side
and in the left hand side
(respectively, only in the right %left
hand side)
while you go around one time,
we call $D$ a {\it ladybug configuration}
(respectively,
 {\it quasi-ladybug configuration}).

If $F$ is the 2-sphere,
our definition of ladybug configurations
is the same as
that
in \S\ref{tamago}.

\bb
Let $D$ be
a ladybug (respectively, quasi-ladybug) configuration.
Then $Z(D)$ have only one circle and $A(D)$ have only two arcs.
Make $s(D)$.
Give $D$ (respectively, $s(D)$) a labeling $x$  (respectively, $y$).
%$x_+$ (respectively, $x_-$).
%Call the resultant labeled resolution configuration $(D,x)$ (respectively, $(s(D),y)$).
We call the decorated resolution configuration $(D,y,x)$
a {\it decorated resolution configuration associated with
the ladybug $($respectively, quasi-ladybug$)$ configuration $D$}.
\end{defn}

 Note that  $(D,y,x)$ may be empty %the empty set
 as explained below.

Since each of $D$ and $s(D)$ has only one circle,
we can let $x_+$ or $x_-$ denote $x$ (respectively, $y$).

\begin{figure}
\includegraphics[width=170mm]{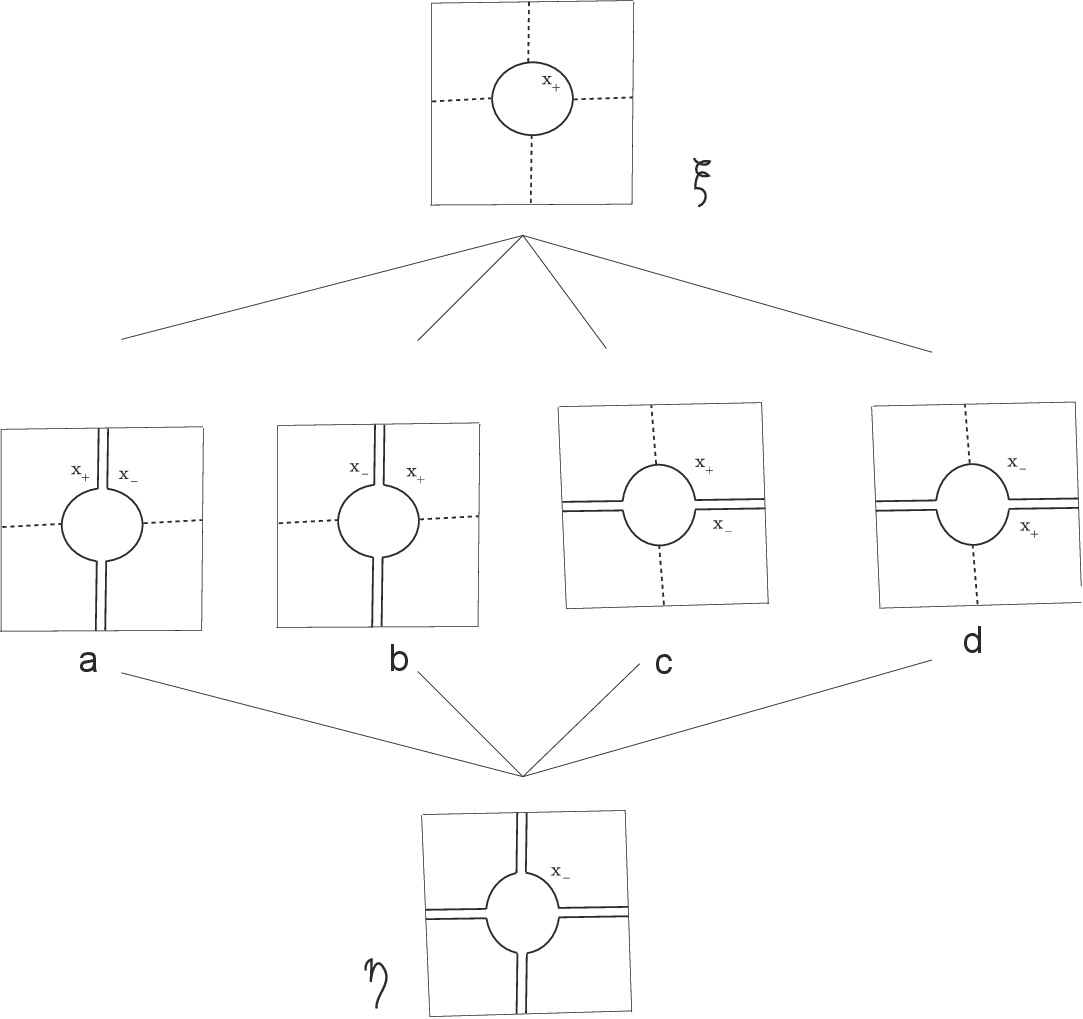}
\caption{{\bf
The poset for
a decorated resolution configuration $(D,x_-,x_+)$ associated with
a quasi-ladybug configuration on $T^2$:
We envelope $T^2$ along two circles as usual,
and draw six labeled resolution configurations.
Here, we have $[\xi;a]\cdot[a;\eta]=[\xi;b]\cdot[b;\eta]=-[\xi;c]\cdot[c;\eta]=-[\xi;d]\cdot[d;\eta].$
}\label{quasiT2}}
\end{figure}

Proposition \ref{daijida} below explains
what difference there is between the torus case and the higher genus surface case.
In the higher genus surface case,
all modulis associated with quasi-ladybug configurations are the empty set.
In the torus case, it does not hold.
Therefore the construction of CW complexes in the torus case is
more complicated than the higher genus surface  case.
See also
`comments below Theorem \ref{tsukiyo}' and
Notes \ref{konki} and \ref{kikon}.

\begin{pr}\label{daijida}
\h$(1)$
Let $F$ be the torus.
Let $D$ be  a quasi-ladybug configuration in $F$.
%Let $(D,y,x)$ be a decorated resolution configuration  associated with $D$.
Assume that the only one circle in $D$ is contractible.
Then there is a non-vacuous decorated resolution configuration
$(D, x_-,x_+)$ associated with $D$.

\bs
\h$(2)$
Let $F$ be a higher genus surface.
Let $D$ be  a quasi-ladybug configuration  in a surface $F$.
Assume that the only one circle in $D$ is contractible.
Let $(D,y,x)$ be a decorated resolution configuration  associated with $D$.
Then $(D,y,x)$ is empty %the empty set
for arbitrary $x$ and $y$.

\bs
\h$(3)$
Let $F$ be an arbitrary surface.
 Let $D$ be a ladybug $($respectively, quasi-ladybug$)$ configuration in $F$.
 Let $(D,y,x)$ be a decorated resolution configuration  associated with $D$.
Assume that the only one circle in $D$ is non-contractible.
Then $(D,y,x)$ is empty %the empty set
for arbitrary $x$ and $y$.

\bs
\h$(4)$
Let $F$ be an arbitrary surface.
There is a ladybug configuration $D$ in $F$
such that
a decorated resolution configuration $(D,y,x)$ associated with $D$
is non-empty. %not the empty set.
\end{pr}

%\h{\bf Note.} Compare Proposition\ref{daijida}.(1) with Figure \ref{quasiT2}.

\bb
\h{\bf Proof of Proposition \ref{daijida}.}
The proof of Proposition \ref{daijida}.(1).
In Figure \ref{quasiT2}, we draw an example.
%We use a net of surfaces in which link diagrams exist, and draw labeled resolution configurations.
%See e.g. Figure \ref{quasiT2}.
% of a non-vacuous decorated resolution configuration $(D,x_-,x_+)$ associated with a quasi-ladybug configuration $D$ in $T^2$.
\\

\h The proof of Proposition \ref{daijida}.(2).
Let $D$ be a quasi-ladybug configuration in $F$.
Let $(D,y,x)$ be the decorated resolution configuration  associated with $D$.

Let $C$ be only one circle in $D$.
Let $A$ and $A'$ be just two arcs in $A(D)$.

%Suppose that $C$ is non-contractible. The number of the set \\ $\{{\bf z}|{\bf z}$ is a labeled resolution configuration. ${\bf z}\prec(D,x), {\bf z}\neq(D,x)\}$ \\is just two by the rule in Definitions \ref{2.10} and \ref{korekos}. Therefore the decorated resolution configuration $(D,y,x)$ associated with $D$ is the empty set.

%Therefore $D$ is contractible. \\

Let $A_p$ and $A_q$ (respectively, $A'_p$ and $A'_q$) be
the endpoints of $A$ (respectively, $A'$).
We can suppose that,
when we go around $C$ one time in an orientation,
we meet $A_p$, $A'_p$, $A_q$, and $A'_q$ in this order,
without loss of generality.
We obtain four circles in $F$, as below:
\\

\h
A circle made of the following two:
The arc $A$.
An arc which is a part of $C$,
whose boundary is $A_p\amalg A_q$, and
which includes $A'_p$.
\\

\h
A circle made of the following two:
The arc $A$.
An arc which is a part of $C$,
whose boundary is $A_p\amalg A_q$, and
which includes $A'_q$.
\\

\h
A circle made of the following two:
The arc $A'$.
An arc which is a part of $C$,
whose boundary is $A'_p\amalg A'_q$, and
which includes $A_p$.
\\

\h
A circle made of the following two:
The arc $A'$.
An arc which is a part of $C$,
whose boundary is $A'_p\amalg A'_q$, and
which includes $A_q$.
\\

Since $D$ is a quasi-ladybug configuration,
the four circles above are non-contractible.

%Therefore we have the following:
The two circles in $s_A(D)$ (respectively, $s_{A'}(D)$)
divide $F$ into two connected compact surfaces,
$W_A$ and $W'_A$
(respectively, $W_{A'}$ and $W'_{A'}$),
with boundary.

Note that the only one circle in $D$ is contractible,
and that the genus of $F$ is greater than one.
Therefore one of
$W_A$ and $W'_A$
(respectively, $W_{A}$ and $W'_{A'}$)
is an annulus and the other has the genus greater than one.

Since the four circles above are non-contractible,
$A'$ (respectively, $A$) is not included in the annulus.

Therefore the only one circle in $s(D)$
divided $F$ into two components: \\
a compact surface with boundary $S^1$ and with genus one,
and \\
a compact surface with boundary $S^1$ and with genus greater than zero.

Therefore the only one circle in $s(D)$ is non-contractible.

Use the rule in Definitions \ref{2.10} and \ref{korekos} again.
Since the two circles in $s_A(D)$ (respectively, $s_{A'}(D)$) and
the only one circle in $s(D)$
are non-contractible,
 $(D,y,x)$ is empty. %the empty set.

%We arrived at a contradiction.

This completes the proof of Proposition \ref{daijida}.(2).
\\

\h
The proof of Proposition \ref{daijida}.(3).
 By the rules in Definitions \ref{2.10} and \ref{korekos}.
\\

\h
The proof of Proposition \ref{daijida}.(4).
The ladybug configuration in Figure \ref{uenp} is put in a 2-disc in $S^2$.
Put it in a 2-disc embedded in $F$.

Note that, of course, there is a ladybug configuration
which is not put in a 2-disc in $F$:
Let $D$ be a ladybug configuration in a surface $F$.
%%We draw Figure \ref{uenp} in \S\ref{tamago} associated with $D$.
Let $A$ be an arc in $A(D)$. Note that $s_A(D)$ may include a non-contractible circle.
However
there is a ladybug configuration
which is put in a 2-disc in $F$. \qed\\

Therefore we must divide our discussion
into three cases: $F=S^2$, $F=T^2$, and $F$ is a higher genus surface.
In \cite{LSk}, Lipshitz and Sarkar did the $S^2$ case.
In this paper, we obtain new results mainly about the higher genus surface case,
and point out the difficulties of the $T^2$ case.
In a sequence \cite{zoku} of this paper, we will write the detail of the $T^2$ case.
%\vskip7mm

\begin{defn}\label{tentomushi}
Let $D$ be a ladybug configuration.
Let $C$ be only one circle in $Z(D)$.
Recall the round trip of $C$, used above when we define  ladybug configurations.
Cut the circle at the four points where the %dotted
arcs meet the endpoints.
The circle is then divided into four pieces.
Recall that, at the beginning point, you see an arc %a dotted line
on the right %left
hand side.
We call  the first and third %second and fourth
pieces of the four, which you are in while your trip,
the {\it right pair},
and call the other two the {\it left pair.}
Note that the orientation of your trip
and the place where you stand at the beginning of your trip
do not change the right and the left pair.
Note also that,
if $F$ is the 2-sphere,
this definition is the same as the one in \cite[\S5.4.2]{LSk}
and in \S\ref{tamago}.
\end{defn}

It is important that
we cannot determine the right and left pair
in the case of  quasi-ladybug configurations
by this method.
(We pose the question: Can one find a method to define the right and the left pair
for quasi-ladybug configurations,
to be compatible with the construction of Khovanov homotopy type?)
\bb

By using the right and left pairs introduced above,
we determine the right and left pair of the labeled resolution configurations
in the middle row of
the poset for a given decorated resolution configuration
associated with a ladybug configuration
%(Figure \ref{uenp} is an example of such decorated resolution configuration.).
(See Figure \ref{uenp} for an example.).
%How to determine
The determination
is  explained in \cite[Figure 5.1 and its explanation in \S5.4.2]{LSk}.
%and \S\ref{tamago}.
\\

%Note. In Figure \ref{quasiT2}, we have $[\xi;a]\cdot[a;\eta]=[\xi;b]\cdot[b;\eta]=-[\xi;c]\cdot[c;\eta]=-[\xi;d]\cdot[d;\eta].$\\

In this paper, we take the right pair
when we construct a CW complex if there is a ladybug configuration.
(If we take the left pair, we can construct a CW complex in  a parallel method.)
\bb

However,
in the case of quasi-ladybug configurations,
we cannot distinguish the two cases.
%Therefore we make two CW complexes if there is a quasi-ladybug configuration. \\

It means that, in general,  in the case of link diagrams in $T^2$,
we may associate more than one CW complex to a single link diagram.

In the case of the higher genus surfaces,
we give only one CW complex to a single link diagram.
{\it Reason.} By Proposition  \ref{daijida},
the decorated resolution configuration associated with an arbitrary quasi-ladybug configuration is empty. %the empty set.
Therefore the moduli associated with it is the empty set.

%\begin{pr}\label{gek}  Let $D$ be a ladybug $($respectively, quasi-ladybug$)$ configuration in a surface $F$. Assume that the only one circle in $D$ is non-contractible. Then any decorated resolution configuration $(D,y,x)$ associated with $D$ is the empty set.  \end{pr}

%\h{\bf Note.} The genus of $F$ in the above proposition may be one.  Some of the above proposition is included in Proposition \ref{daijida}. \\

%\h{\bf Proof of Proposition \ref{gek}.} By the rules in Definitions \ref{2.10} and \ref{korekos}. \qed\\

\section{\bf A  moduli
in the case of link diagrams in the torus,
which never appears
in the case of those in $S^2$ nor in the case of those in higher genus surfaces
}\label{chigau}

\begin{figure}
\vskip-10mm
\includegraphics[width=168mm]{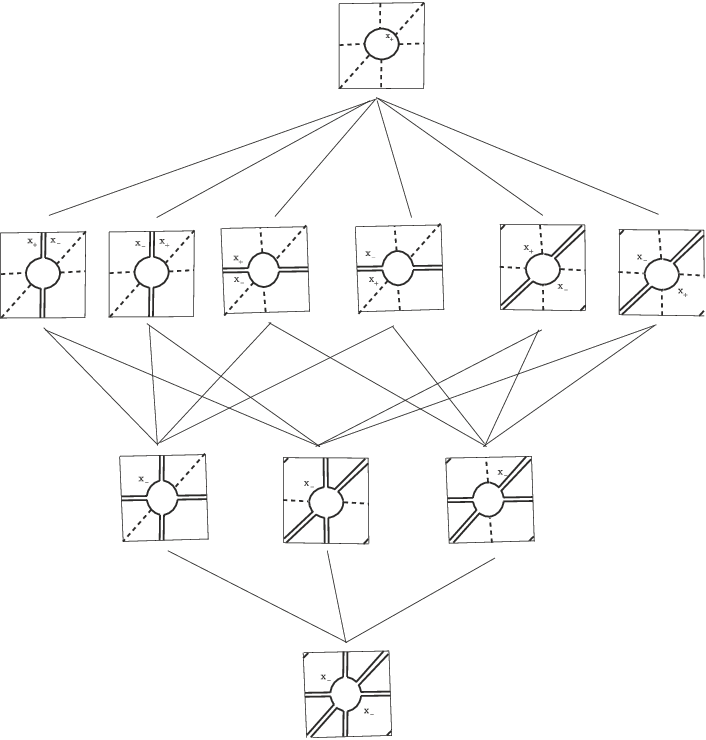}
\caption{{\bf
The poset for the decorated resolution configuration
associated with a quasi-ladybug configuration in $T^2$:
We use a net of $T^2$ with a labeled resolution configuration.
Each one in the second row includes just two circles.
Each one in the third row includes only one circle.
The lowest one includes two circles.
}\label{uncon}}
\end{figure}

\h
In Figure \ref{uncon}, we draw the poset for a decorated resolution configuration.
Each labeled resolution configuration in it is put in $T^2$.
Assume that we give more than one moduli for a single decorated resolution configuration in general,
and make many CW complexes for a single link diagram,
as explained in \S\ref{hiyoko}.
Then one of moduli spaces for the decorated resolution configuration
in Figure \ref{uncon} is a dodecagon.
It is not the 3-dimensional cube moduli.
Of course, it is also not a trivial covering of the 3-dimensional cube moduli.
This is a new phenomenon which we do not have in the $S^2$ case.
We also do not have it in the higher genus case.
That is a reason why the torus case is difficult.

\bs
\h{\bf Review.} See \cite[Proposition 5.2]{LSk}:
In the $S^2$ case, each moduli space is
the empty set or a trivial covering of the $n$-dimensional cube moduli.
\bs

%{\bf Should we write more explicitly?}

%\begin{figure}[H] \includegraphics[width=180mm]{sugoi2.jpg}    \caption{{\bf }\label{sugoi2}}   \end{figure}\np

\bb
\section{\bf The Khovanov-Lipshitz-Sarkar stable homotopy type
for links in thickened higher genus surfaces
}\label{KLShomotopy}

\h
In this section,
we define
the Khovanov-Lipshitz-Sarkar stable homotopy type
for links in thickened higher genus surfaces
%(Definition \ref{Khflow}).
(\S\ref{te-gi}).

\subsection{Moduli spaces for decorated resolution configurations}\label{modu}\hskip10mm

\h
\begin{defn}\label{cucu}
Take an $n$-dimensional cube in a coordinate space $\R^n$
whose vertices are points with coordinates $(a_1,...,a_n)$,
where each $a_i$ is 0 or 1.
Let $C_n$ be the set of all of these vertices.
Elements of $C_n$ can be regarded as vectors.
$C_n$ has a partial order $\prec$ as follows:
Let $u=(u_1,...,u_n), v=(v_1,...,v_n)\in C_n$.
Let $k\in\{1,...,n\}$.
If
$u_k=0,v_k=1$ and $u_i=v_i$ for $i\neq k$,
then $u\prec v$.
It follows from this definition that
if $u \prec v$ and $v \prec w$ then $u \prec w$.

We define the {\it norm} $|u|$ of $u$ to be   $u_1+...+u_n$.
(Note that $u_1+...+u_n=(u_1)^2+...+(u_n)^2$ in this situation,
since $0^2=0$ and $1^2=1$.)
\end{defn}

Recall the norm for resolution configurations in Defintion \ref{2.2}.
The {\it $n$-dimensional cube flow category} $\mathcal C_C(n)$
in \cite[\S4]{LSk} is associated with $C_n$.
We use the moduli
$\mathcal M_{{\mathcal C}_C(n)} (\bar 0,\bar 1)$
for $\mathcal C_C(n)$,
defined there.

\bb

%We define a CW complex for the Khovanov chain for links in thickened  higher genus surfaces as in \cite[\S5]{LSk}.

%Let $\mathcal L$ be a link in a thickened higher genus surface. Let $L$ be a link diagram which represents $\mathcal L$.

%The set of all labeled resolution configurations of $L$ is a basis of the Khovanov chain complex of $L$. We call this set, the {\it Khovanov basis}. Each element of the Khovanov  basis is called a {\it Khovanov basis element}. Note that a Khovanov basis element is a labeled resolution configuration.

Let $(D, x, y)$ be an index $n$ basic decorated resolution configuration
in a higher genus surface.
%(See  Definition \ref{2.4} for the term, ``basic''. See \cite[Definiion 2.11]{LSk} for decorated resolution configurations.)
Let ${\bf y}=(D,y)$ and ${\bf x}=(s(D),x)$.
%For ${\bf x}$ and ${\bf y}$,
We associate to each non-vacuous $(D,x,y)$
an $(n-1)$-dimensional $<n-1>$-manifold,
$\mathcal M(D,x,y)$ or $\mathcal M({\bf x,y})$,
together with an $(n-1)$-map

$$\mathcal F:\mathcal M(D,x,y)\to \mathcal M_{{\mathcal C}_C(n)} (\bar 0,\bar 1)$$

%See \cite[\S3.1]{LSk} for $<m>$-manifolds and $m$-maps, where $m$ is an integer.
%we make a moduli space $\mathcal M({\bf x,y})$. % as in \cite[\S5.1]{LSk}.
%See Definition \ref{decor} below.
\h below, as done in \cite[all of \S5. In particular, see \S5.1 and Proposition 5.2.]
%, especially, Definition 5.3]
{LSk}.
%See \cite[\S4]{LSk} for the definition of  $\mathcal M_{{\mathcal C}_C(n)} (\bar 0,\bar 1)$.
%If $(D,x,y)\\=\emptyset$,
If $(D,x,y)$ is empty,
we associate to $\mathcal M({\bf x, y})$ the empty set.
%If $(D,x,y)=\emptyset$, then $\mathcal M({\bf x, y})$ is the empty set.
%The cell corresponding {\bf x} does not touch {\bf y}.
We suppose that $(D,x,y)$ is not
empty below.
%the empty set below.

%See \cite[Definition 3.3]{LSk} for the definition of $\bar 0,\bar 1$.

%Note that  we regard each index $n$ basic decorated resolution configurationas a flow category. See \cite[Definition 3.12]{LSk} for the definition of flow category.

This means that we make a moduli space for
any pair of Khovanov basis elements in the Khovanov chain complex of $L$.
%  a moduli space for each pair of Khovanov basis elements %generators,
%and a CW complex for a Khovanov chain complex for links in thickened surfaces.

\bb

\begin{defn}\label{pi}
Take any element ${\bf r}\in P(D,x,y)$.
Let $v$ be the vector of ${\bf r}$.
Define
$\pi$ to be the map $P(D,x,y)\to C_n$ such that $\pi({\bf r})=v$.
\end{defn}

This map $\pi:P(D,x,y)\to C_n$ keeps the partial order
because $\text{gr}_h((D_L(u), x))\\=-n_-+|u|$.
\\

We define the map $\mathcal F$ to be associated with $\pi$,
as done in \cite[all of \S5. In particular, see \S5.1 and Proposition 5.2.]{LSk}.

\begin{pr}\label{toro}
The map $\pi$ in Definition \ref{pi} is onto. % if $n\geqq3$.
\end{pr}

\begin{thm}\label{hanami}
The map $\mathcal F$ is a trivial covering map.
\end{thm}

Theorem \ref{hanami} corresponds to \cite[Proposition 5.2 and \S5.1]{LSk}.
See \cite[\S3.4.1]{LSk} for the definition of trivial covering maps in this case.

While we prove Proposition \ref{toro} and Theorem \ref{hanami}  in the following subsections,
we prove Theorem \ref{tsukiyo} together.
In the proof of Theorem \ref{hanami},
Theorem \ref{tsukiyo} plays a crucial role.

\begin{thm}\label{tsukiyo}
Let $(D, x, y)$ be an index $n$ basic decorated resolution configuration
in a higher genus surface above.
Suppose that $(D, x, y)$ is non-vacuous, %not the empty set,
as above.
Take any element $A$ in $P(D,x,y)$.
$($Recall that $A$ is a labelled resolution configuration.$)$
Then all arcs in $A$ are mc arcs.
\end{thm}

In the case of links in $S^3$,
the map $\pi$ is an epimorphism for all natural numbers $n$,
 and
the above theorems are true (\cite{LSk}).

In the case of links in the thickened torus and
in the case of virtual links,
we have different features.
See %\S\ref{konnyaku} and
\S\ref{chigau} and \cite{KauffmanOgasasq}.
\\

We prove Proposition \ref{toro}, Theorems \ref{hanami} and \ref{tsukiyo} below.

\bb
\subsection{The case $n=1$}\label{n1} \hskip10mm

\h{\bf Proof of the case $n=1$ of Theorem \ref{hanami} and
that of the case $n=1$ of Proposition \ref{toro}.}
 We %can
 associate to
$\mathcal M({\bf x, y})$ one point
because the coefficient  of the right hand side of Equality (\ref{bibun})
is $+1,0,$ and $-1$.
Therefore Theorem \ref{hanami} and Proposition \ref{toro} hold in this case.
\\

\h{\bf Proof of the $n=1$ case of Theorem \ref{tsukiyo}.}
%$(D,y)$ is the minimal element of the decorated resolution configuration $(D,x,y)$.
Suppose that only one arc in $D$ is a scs arc.
%By Proposition \ref{ichigo}, $\delta(D,x)=0.$ Hence
By the rule in Definitions \ref{2.10} and \ref{korekos},
the decorated resolution configuration is empty. %the empty set.
We arrived at a contradiction.
\qed

\bb
\subsection{The case $n=2$}\label{n2} \hskip10mm

\h{\bf Proof of the case $n=2$ of Theorem \ref{hanami} and
that of the case $n=2$ of Proposition \ref{toro}.}
 Let $P$ be a partial ordered set.
Let $\alpha, \beta\in P$.
Define $Q(\alpha,\beta)$ associated with $P$ to be
$\{q|q\in P, q\neq\alpha, q\neq\beta, \alpha\prec q, q\prec\beta\}$.
We omit the words,  associated with $P$, when it is clear from the context.
Let $A$ be a finite set. Let $\# A$ be the number of all elements of $A$.

Take $Q({\bf x,y})$ associated with $P(D,x,y)$.
%Let $\{p|p\in (D,x,y), {\bf y}\prec p, p\prec {\bf x}, p\neq{\bf x}, p\neq{\bf y}\}$.
%
Then $\#Q({\bf x,y})$ is 2 or 4 because of Definitions \ref{2.10} and \ref{korekos}
since we assume that $(D,x,y)$ is non-vacuous. %not the empty set.
 % See also Theorem \ref{sikon} and its proof. %See also \cite[\S 5.4]{LSk}.
If $\#Q({\bf x,y})$ is 2,
we %can
associate to $\mathcal M({\bf x, y})$  one segment.

Suppose that  $\#Q({\bf x,y})$ is 4,
$(D,x,y)$ is a ladybug configuration and $D$ includes only one contractible circle
because of Proposition \ref{daijida}.
%Propositions \ref{daijida} and \ref{gek}.
We %can
associate to
$\mathcal M({\bf x, y})$
  a disjoint union of two segments.
  We use the right pair in \S\ref{hiyoko}.

Therefore Theorem \ref{hanami} and Proposition \ref{toro} hold in this case.
\\

\begin{note}\label{insight}
 If we choose the left pair, we can also construct the Khovanov stable homotopy type.
 In the case of links in $S^3$,
 both choices give the same  Khovanov stable homotopy type (Fact \ref{konnya}).
We pose the question: Do both choices give the same  Khovanov stable homotopy type
in the case of links in thickened (higher genus) surfaces?
\end{note}

\h{\bf Proof of the $n=2$ case of Theorem \ref{tsukiyo}.}
Suppose that there is a scs arc in
a labeled resolution configuration in $P(D,x,y)$.
 Then the decorated resolution configuration is empty. %the empty set.
{\it Reason.} Consider index 2 decorated resolution configurations
associated with
 Figures  \ref{Hp}-\ref{HqZ}.
 Other cases are easy.

We arrived at a contradiction.
\qed

\begin{note}\label{konki}
In the case of virtual links in \cite{DKK, KauffmanOgasasq, Man},
 Theorem \ref{tsukiyo} is not true
and the map $\pi$ is not an epimorphism,
even if $n=2$.
%See Figure \ref{eK4}.
See \cite{KauffmanOgasasq}.
See Note \ref{kikon}.
\end{note}

\bb
\subsection{The first part of the case $n\geqq3$
}\label{n34} \hskip10mm

\h {\bf Proof of Proposition \ref{toro}.}
The case $n=1,2$ holds by \S\ref{n1} and \S\ref{n2}.
We prove the case $n\geqq3$.
Take two arbitrary labeled resolution configurations,
${\bf w}=(D_\xi(L),w)$ and ${\bf z}=(D_\zeta(L),z)$,
in $P(D,x,y)$
such that $|\xi|+2=|\zeta|$,
and such that ${\bf w}$$\prec$${\bf z}$.
Since
%$(D,x,y)\neq\emptyset$,
$(D,x,y)$ is not empty,
there is such a pair.
Assume that $Q(\xi,\zeta)$ associated with $C_n$ is $\{\eta, \eta'\}$.

Take $Q({\bf w,z})$ associated with $P(D,x,y)$.
By \S\ref{n2} and \cite[\S 5.4]{LSk},
we have that
$\#Q({\bf w,z})$ is 2 or 4.
If it is 2, one element of $Q({\bf w,z})$
is over $D_\eta(L)$ and the other is over  $D_{\eta'}(L)$
by the rule in Definitions \ref{2.10} and \ref{korekos}.
If it is 4, two elements of $Q({\bf w,z})$
are over  $D_\eta(L)$ and the other two are over $D_{\eta'}(L)$
by the rule in Definitions \ref{2.10} and \ref{korekos}.
Repeat this procedure.
Hence $\pi$ is onto. \qed\\

\h{\bf Proof of Theorem \ref{tsukiyo}.}
We have proved the $n<3$ case in the previous subsections.
We prove the $n\geqq3$ case. % of Theorem \ref{tsukiyo}.
We prove by reductio ad absurdum.
Assume that ${\bf g}\in P(D,x,y)$ has a scs arc $A$.
Let $G$ be a (non-labeled) resolution configuration under ${\bf g}$.
Carry out a surgery along $A$, and obtain a (non-labeled) resolution configuration $G'$.
Let $v$ be the vector of $G'$.
Then $v\in C_n$.
By Proposition \ref{toro}, we have $\pi^{-1}(v)\neq\emptyset.$

Any labeled resolution configuration on $G'$
has a different quantum degree from that of $g$ %(respectively, ${\bf x, y}$).
(See Proposition \ref{ichigo}.). %Fact \ref{qde}.).
(Note that gr$_q{\bf g}=$gr$_q{\bf x}=$gr$_q{\bf y}$.)
By the definition of decorated resolution configurations,
all labeled resolution configurations in $P(D,x,y)$
have the same quantum degree.
Hence  $\pi^{-1}(v)=\emptyset.$
We arrived at a contradiction.
\qed\\

\h{\bf Review.}
In \cite[Proof of Proposition 5.2]{LSk},
the case of three arcs
 is more complicated than the case of greater than three arcs.
(In our way of description in this paper,
the case of three arcs is the case $n=3$.)
 %Here, note that $n=3$ is described by using the notation $n$ in this paper. In \cite{LSk}, there is used another notation.)
The reason is as follows:
Let $f$ be a local diffeomorphism map from $X$ to $S^m$.
If $m>1$, $f:X\to S^m$ is a trivial covering map.
If $m=1$, $X$ is not a trivial covering map in general.
An example is the connected double covering map $S^1\to S^1$.
(Here, $m+2$ is the above $n$.)

The case of three arcs
in \cite[Proof of Proposition 5.2]{LSk}
is proved  in \cite[\S5.5]{LSk} by
checking  all  resolution configurations with three arcs.

In this paper, the case $n=3$ is also more complicated than the case $n>3$.

\bb
We split the case $n\geqq3$ into the case $n=3$ and the case $n\geqq4$ below.

\bb
\subsection{The proof of Theorem \ref{hanami} in the case $n=3$}\label{n3} \hskip10mm

 There are just two cases.

\bb
\h{\bf Case 1.}
Suppose that all circles in all labeled resolution configurations in $P(D,x,y)$
are contractible.

By  Theorem \ref{tsukiyo},
all arcs in all labeled resolution configurations in $P(D,x,y)$
are mc arcs.
Hence
we can prove  Theorem \ref{hanami} in this case,
as proved in \cite[\S 5.1, especially Proposition 5.2, and 5.5]{LSk}.
\bb

\h{\bf Case 2.}
Assume that
there is a non-contractible circle
in a labeled resolution configuration in $P(D,x,y)$.

%By  Theorem \ref{tsukiyo}, all arcs in each labeled resolution configuration in $(D,x,y)$ are mc arcs.

Recall that $\partial\mathcal M_{\mathcal C_{(3)}}(\bar1,\bar0)$ is a circle.
By Proposition \ref{toro} %\ref{katsuo}
and
\cite[Definition 3.12, especially (M-3) in it, and Proposition 5.2]{LSk},
we have the following:
$\partial\mathcal M({\bf x,\bf y})$ is a disjoint union of circles.
Furthermore, $\partial\mathcal M({\bf x,\bf y})$ is a covering space of
$\partial\mathcal M_{\mathcal C_{(3)}}(\bar1,\bar0)$ by $\mathcal F$.

\begin{fact}\label{hamachi}
Any circle which is a connected component of $\partial\mathcal M({\bf x,\bf y})$
covers
$\partial\mathcal M_{\mathcal C_{(3)}}(\bar1,\bar0)$
by a degree one map.
\end{fact}

\h{\bf Proof of Fact \ref{hamachi}.}
By  Theorem \ref{tsukiyo},
all arcs in all labeled resolution configurations in $P(D,x,y)$
are mc arcs.
We check all basic resolution configurations with three mc arcs.
We can use almost the same methods as those in \cite[\S5.5]{LSk}.

We must take care of the difference between
 the rule of the partial order in Definition \ref{2.10}
 and
 that in \cite[Definition 2.10]{LSk}.

The result (\cite[Lemma 2.14]{LSk})
on leaves and co-leaves is used in \cite[\S5.5]{LSk}.
This result is generalized easily, and
also holds in our case.
Therefore we only have to concentrate on
resolution configurations without a leaf or a co-leaf, as  in \cite[\S5.5]{LSk}.

The duality theorem, \cite[Lemma 2.13]{LSk}, is used in \cite[\S5.5]{LSk}.
This result is generalized easily, and
also holds in our case.
It helps our purpose below.

%Therefore it is enough to check the four cases, \cite[Figures 5.3.d-5.3.g]{LSk}, directly.

%We must take care of the following facts.  The circles in \cite[Figures 5.3.d-5.3.g]{LSk}  may be non-contractible in our case.  There may appear a non-contractible circle after a resolution in our case.

%We use Propositions \ref{daijida} and \ref{gek} about ladybug configurations and quasi-ladybug ones.

If the union of all arcs and circles in a basic resolution configuration is disconnected,
there is a leaf or a co-leaf. Note that it is basic.
Hence we assume that it is connected.

\begin{figure}
\centering\targ{4}{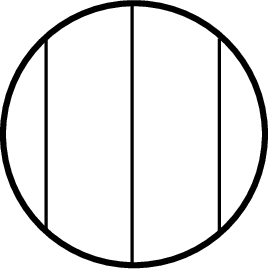}{1} \quad \targ{4}{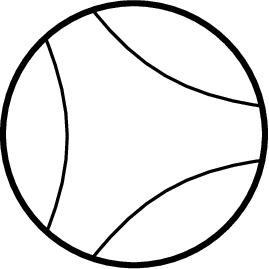}{2} \quad \targ{4}{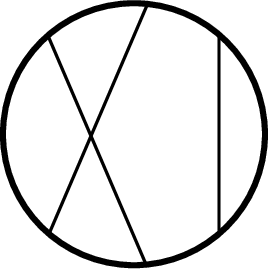}{3} \quad \targ{4}{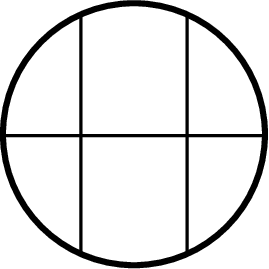}{4} \quad \targ{4}{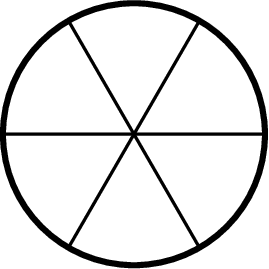}{5} \\
\targ{3.5}{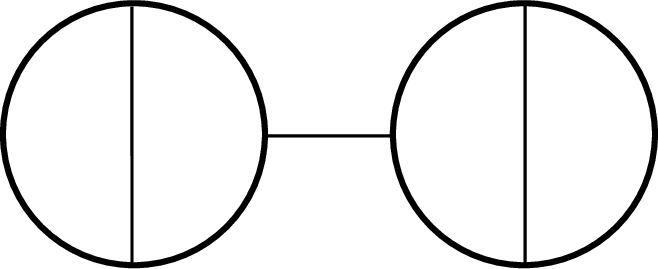}{6} \quad \targ{3.5}{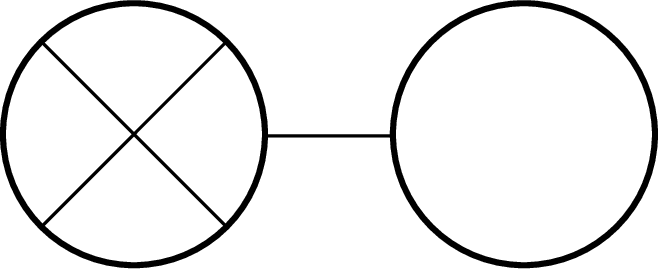}{7} \quad \targ{3.5}{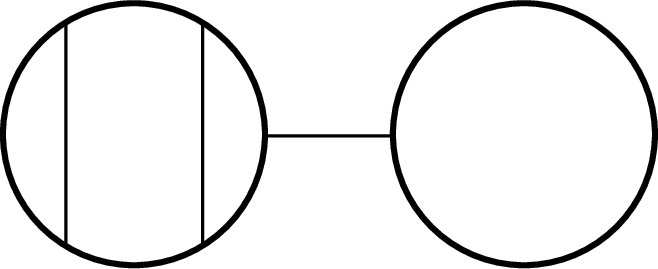}{8} \\
\targ{3.5}{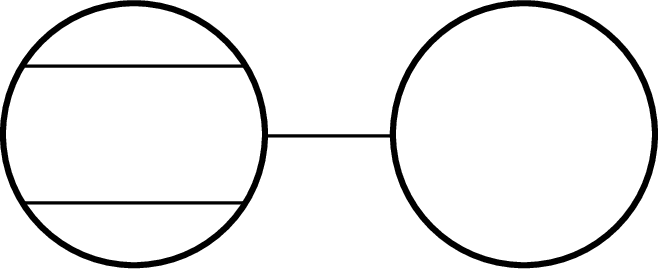}{9}\quad \targ{3.5}{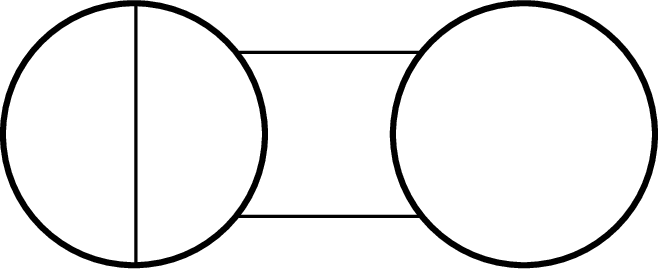}{10}\quad \targ{3.5}{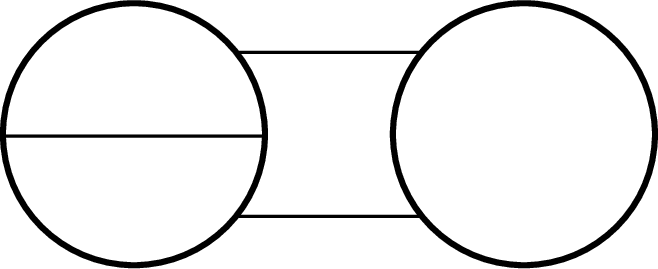}{11} \\
\targ{3.5}{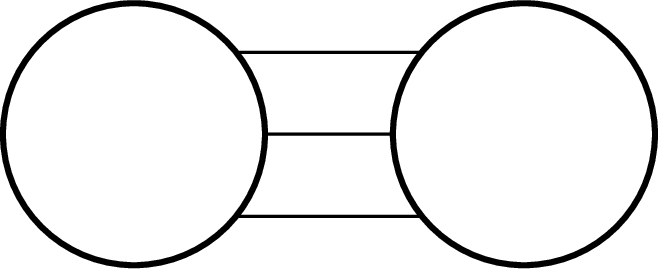}{12} \quad \targ{4}{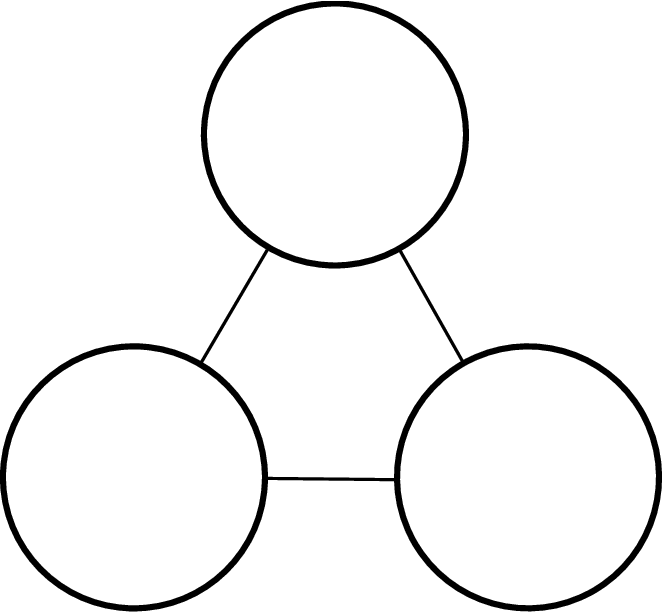}{13} \quad \targ{3.5}{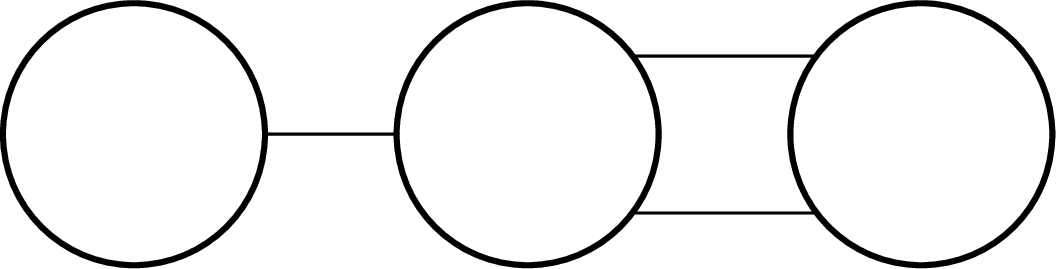}{14} \\
\targ{3.5}{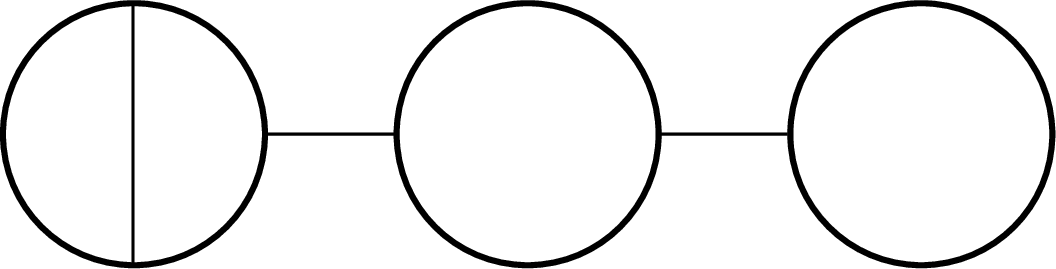}{15}\quad \targ{3.5}{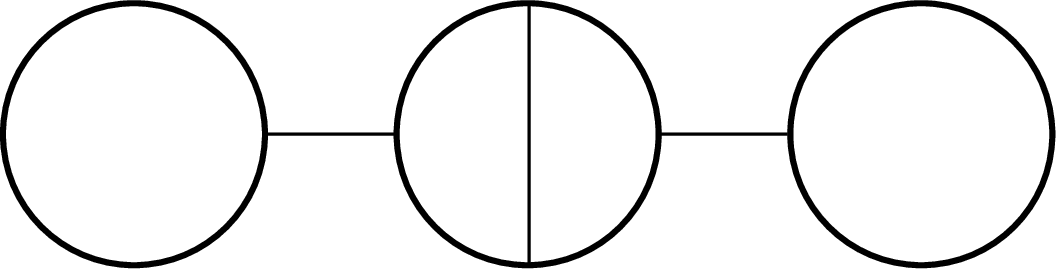}{16} \\ \targ{3}{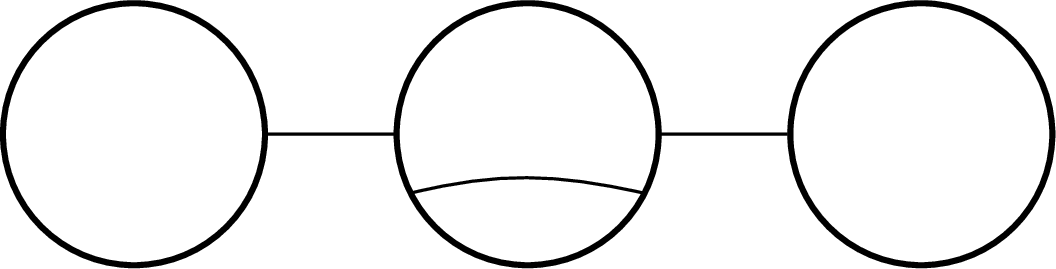}{17} \quad \targ{3}{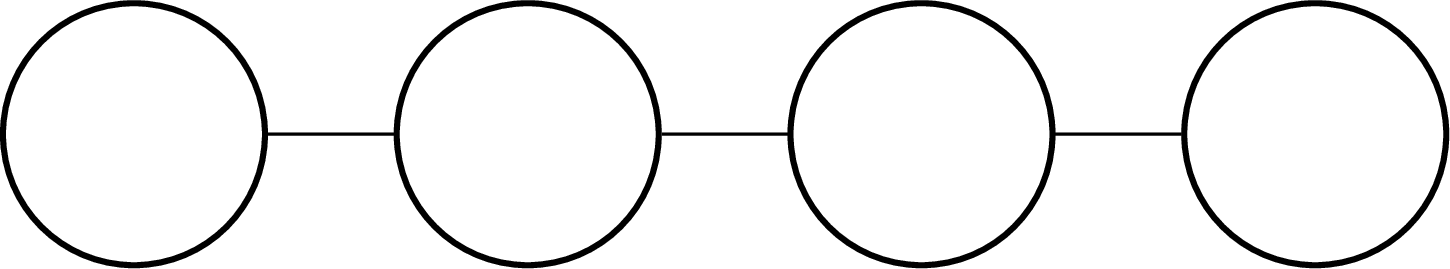}{18} \\ \targ{7}{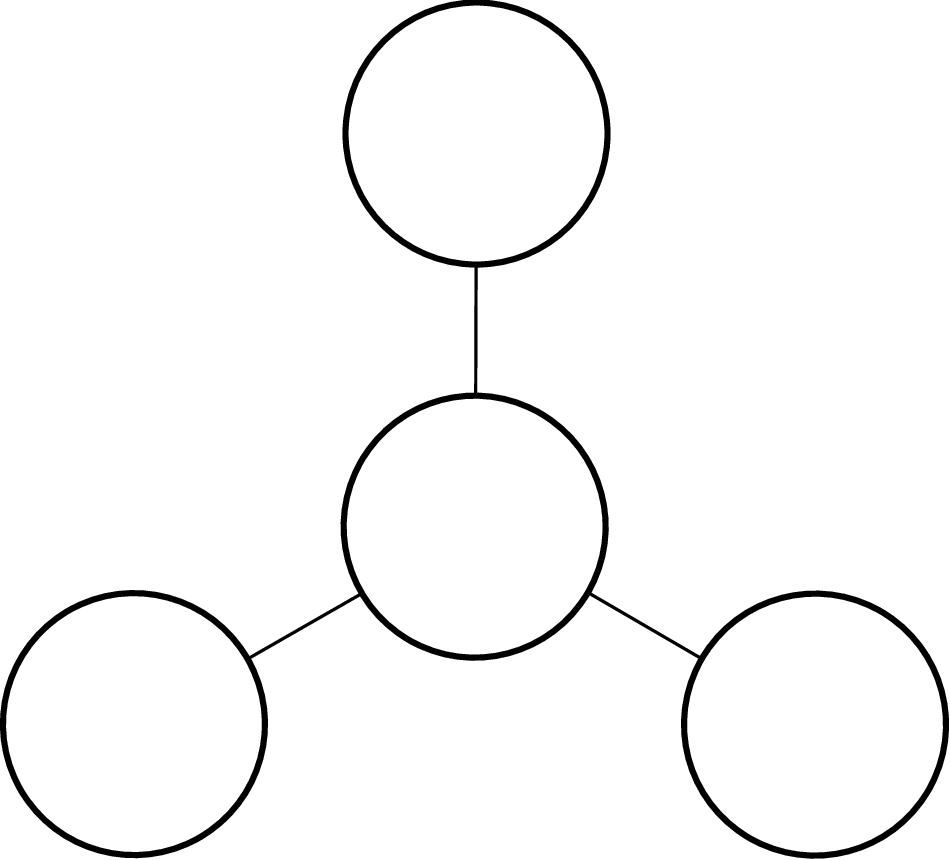}{19}
\caption{{\bf
Connected graphs of the resolution configurations of index 3:
%All cases of $Z(D)\cup A(D)$ are drawn like abstract graphs in the case where
% $D$ is a resolution configuration with three arcs such that $Z(D)\cup A(D)$ is connected.
The segments denote arcs. We do not use dotted segments here.
We draw only $Z(D)\cup A(D)$ abstractly, like abstract graphs.
}}\label{zenbu}
\end{figure}

Let $D$ be any resolution configuration with three arcs such that
$Z(D)\cup A(D)$ is connected.
We draw all such cases of $Z(D)\cup A(D)$ like abstract graphs in Figure \ref{zenbu}.

\begin{figure}
\includegraphics[width=80mm]{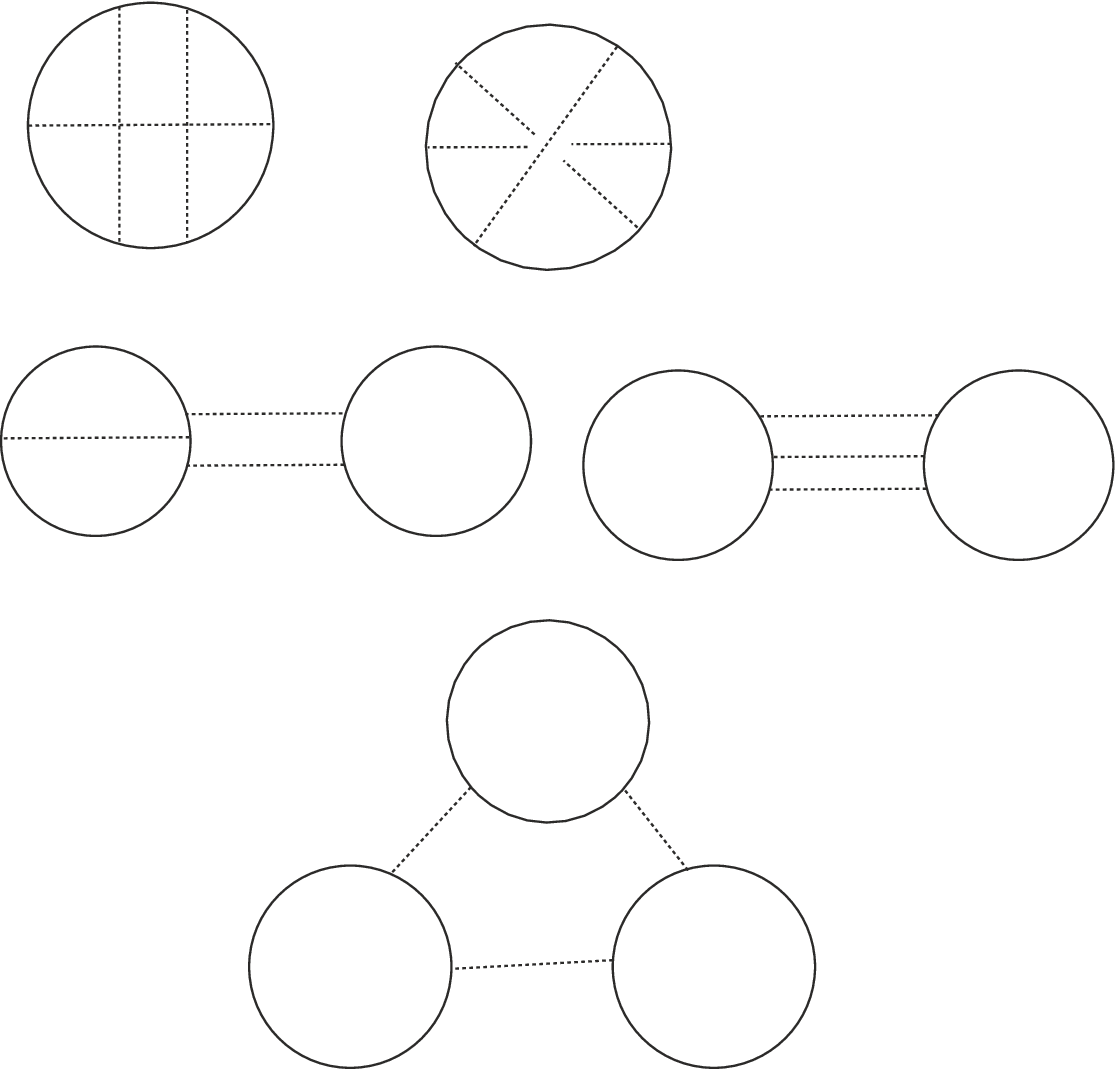}
\caption{{\bf
Graphs of the connected resolution configurations
%We draw abstractly $Z(D)\cup A(D)$ of an arbitrary connected basic resolution configuration $D$
with three arcs, without a leaf or a co-leaf:
Recall that we consider the case where arcs are mc arcs.
}\label{3arc}}   \end{figure}

We choose $Z(D)\cup A(D)$  without a leaf or a co-leaf from  Figure \ref{zenbu},
and draw them in Figure \ref{3arc}.
%In Figure \ref{3arc}, we draw abstractly $Z(D)\cup A(D)$ of an arbitrary connected basic resolution configuration $D$ with three arcs, without a leaf nor a co-leaf.

Recall the following facts  associated with Figures \ref{zenbu} and \ref{3arc}.
In a surface, each $Z(D)\cup A(D)$ has a neighborhood $N$
such that $N$ is a compact surface and
such that the inclusion map of $Z(D)\cup A(D)$ to $N$
 is a homotopy type equivalence map.
There are many homeomorphism types of $N$ in general.
%By Theorem \ref{tsukiyo},
We can assume that no scs arc appears.

\bb
\bb

Let $D$ be a resolution configuration made from
one
%figure
diagram
of Figure \ref{3arc}.
Assume that,
if we choose two arcs and one circle from $D$,
then they make a ladybug configuration (or a quasi-ladybug configuration).
Then we say that $D$
{\it includes a ladybug configuration}
(or {\it a  quasi-ladybug configuration$)$}.
\bb

Let $D$ induce a quasi-ladybug configuration.
By Proposition \ref{daijida},
any decorated resolution configuration starting from $D$ is empty. %the empty set.
\bb

%Let $D$ be a resolution configuration in Figure \ref{Ig1}.
The following two conditions are equivalent.

\bs
\h(1)
$D$ does not include a ladybug configuration.
Let $E$ be any resolution configuration obtained from $D$ by a single surgery.
$E$  is not a ladybug configuration.

\bs
\h(2)
Neither $D$ nor  the dual  resolution configuration  $D^\ast$
includes a ladybug configuration.
\bs

If we have the above condition (1) (respectively, (2)),
the moduli of any decorated resolution configuration starting from $D$
is the empty set or the single 3-dimensional cube moduli.
\bb

Note that, even if $D$ does not  include a ladybug configuration,
$D^\star$ may include a ladybug configuration.
An example is the case where $D$ is \cite[Figure 5.3.g]{LSk}.
\bb

Let $D$ be a resolution configuration made from
one figure of Figure \ref{3arc}.
Then $D^\ast$ is also made from one figure of Figure \ref{3arc}.
{\it Reason.}
Let $G$ be a resolution configuration in Figure \ref{zenbu}.
If $G$ has a leaf (respectively, co-leaf),
then $G^\ast$ has a co-leaf (respectively, leaf).

Apply Propositions \ref{daijida}. % and \ref{gek}.
If the moduli of
a decorated resolution configuration
starting from $D$
is not the empty set nor a single 3-dimensional cube moduli,
$D$ or $D^\ast$ satisfies the condition:
%following condition $(\ast)$:
It includes a ladybug configuration,
and does not include a quasi-ladybug configuration.
The circle in the ladybug configuration is contractible.

Therefore we only have to check
%investigate decorated resolution configurations starting from $D$
%made from  one figure of Figure \ref{3arc},
%which satisfy the condition $(\ast)$ and which does not include a leaf or a co-leaf.
%Recall that we consider the case where arcs are mc arcs.
decorated resolution configurations
starting from $D$  made from the left upper figure in Figure \ref{3arc}.
The   moduli  of each is a disjoint union of the 3-dimensional cube moduli
or the empty set.

This completes the proof of Fact \ref{hamachi}.
\qed\\
 %the following cases.
%make a decorated resolution configuration starting from $D$ with a labeling if $D$ satisfies the following condition: $D$ includes a ladybug configuration and the circle which generates the ladybug configuration is contractible.
%See Figure \ref{Ig2}.

%\bb\begin{figure}
%%\includegraphics[width=150mm]{.jpg}
%\caption{{\bf {\bf \Large Igor's figures.} Decorated resolution configurations which we must check  }\label{Ig2}}   \end{figure}  \bb

Henceforth we have Theorem \ref{hanami} in this case.

\subsection{The proof of Theorem \ref{hanami} in the case $n\geqq4$}\label{n4}\hskip10mm

\h %By proposition \ref{toro},
%we can prove Theorem \ref{hanami} in this case,
% as in \cite[Proposition 5.2]{LSk}.
The case $n\leqq3$ is true by \S\ref{n1}-\ref{n3}. %{n4}.
Therefore the case $n\geqq4$
 is proved by the same method as one in  \cite[Proposition 5.2]{LSk}.
 %true for all natural numbers $n$. \qed
 \bb

\subsection{Definition of the Khovanov-Lipshitz-Sarkar stable homotopy type for links in thickened surfaces}\label{te-gi}\hskip10mm

\h
Definition \ref{Khflow} is made by
adding the fact about the quantum grading to \cite[Definition 5.3]{LSk}.
%We prove that Definition \ref{Khflow} is well defined (Theorem \ref{well}).

\begin{defn}\label{Khflow}
Let $L$ be a link diagram in a higher genus surface.
%The set of all labeled resolution configurations of $L$ is a basis of the Khovanov chain complex of $L$. We call this set, the {\it Khovanov basis}. Each element of the Khovanov  basis is called a {\it Khovanov basis element}. Note that a Khovanov basis element is a labeled resolution configuration.
The {\it Khovanov flow category} ${\mathcal C}_K(L)$ has one object for
each Khovanov basis element.
%of the basis elements of the Khovanov homology.
%  cf. \cite[Definition 2.15]{LSk}.
That is, an object of
  ${\mathcal C}_K(L)$ is a labeled resolution configuration of the form
  ${\bf x}=(D_L(u),x)$ with $u\in\{0,1\}^n$.
  The grading
  on the objects is the homological grading gr$_h$; % from;
%  \cite[Definition 2.15]{LSk};
  the quantum grading  gr$_q$ and the homotopical grading gr$_\mathfrak h$
  are additional gradings on the objects.
  We need the orientation of $L$ in order to define these gradings, but the rest of
  the construction of ${\mathcal C}_K(L)$ is independent of the
  orientation.

  Consider objects
  ${\bf x}=(D_L(u),x)$
   and
  ${\bf y}=(D_L(v),y)$
  of ${\mathcal C}_K(L)$.
  The space
  $\mathcal M_{{\mathcal C}_K(L)}(\bf{x,y})$
  is defined to be empty unless
  ${\bf y}\prec {\bf x}$ with respect to the partial order from
  Definition \ref{2.10}.
  So, assume that ${\bf y}\prec{\bf x}$.
  Let $x|$ denote the restriction of $x$ to
  $s(D_L(v)-D_L(u))\\=D_L(u)-D_L(v)$
  and let $y|$ denote the restriction of $y$
  to $D_L(v)-D_L(u)$.
  Therefore,
  $(D_L(v)-D_L(u), x|, y|)$
  is a basic decorated resolution configuration.
  Recall that we defined $\mathcal M(D_L(v)-D_L(u), x|, y|)$ in \S\ref{modu}.
  %Definition \ref{decor}.
  Define
  $$\mathcal M_{{\mathcal C}_K(L)}({\bf x,y})=  \mathcal M(D_L(v)-D_L(u), x|, y|)$$
   as smooth manifolds with corners.
   The composition maps for the resolution configuration moduli spaces
  (see \cite[(RM-1) in \S5.1]{LSk}) induce composition maps
  $$\mathcal M_{{\mathcal C}_K(L)}(\bf{z,y})\x
  \mathcal M_{{\mathcal C}_K(L)}(\bf{x,z})\to
  \mathcal M_{{\mathcal C}_K(L)}(\bf{x,y})$$

  Given a flow category $\mathcal C$ and an integer $n$,
  let $\mathcal C[n]$ be
  the flow category obtained from $\mathcal C$ by increasing the grading of
  each object by $n$.

  The Khovanov flow category ${\mathcal C}_K(L)$ is equipped with a functor
  $\mathcal F$ to $\mathcal C_C(n)[-n_-]$,
  %\footnote{Henceforth, we will usually suppress the grading shifts from the notation.}
  which is a cover in the sense of
\cite[Definition 3.28]{LSk}.
%  \end{defn} \end{document}

  On the objects,
  $\mathcal F:{\text{Ob}}_{{\mathcal C}_K(L)} \rightarrow {\text{Ob}}_{\mathcal C_C(n)}$
  is defined as
$$\mathcal F: \mathcal M_{{\mathcal C}_K(L)}((D_L(u),x),(D_L(v),y))\to
\mathcal M_{\mathcal C_C(n)}(u,v)$$
  is defined to be composition
$$\mathcal M(D_L(v)-D_L(u), x|, y|)
\stackrel{\mathcal F}{\to}
\mathcal M_{\mathcal C_C(|u|-|v|)}(\overline{1},\overline{0})
\stackrel{{\mathcal I}_{u,v}}{\to}
\mathcal M_{\mathcal C_C(n)}(u,v).
$$
We can say that $\mathcal F$ %in Definition \ref{Khflow}
is associated with $\pi$ in Definition \ref{pi}.

%Since Theorem \ref{well} is true, Definition \ref{Khflow} is well-defined.  That is, w
%We have defined
%By using the moduli spaces defined above, as explained in \S\ref{Khpi},
As explained in \S\ref{Khpi},
use the moduli spaces which are defined above.
Thus we construct
the {\it Khovanov-Lipshitz-Sarkar stable homotopy type} for $L$
and
the  {\it Khovanov-Lipshitz-Sarkar spectrum}
$\mathcal X_{Kh}(L)
=\displaystyle\bigvee_{q,\mathfrak h}\mathcal X^{q,\mathfrak h}_{Kh}(L)$ for $L$.
%\bb
%The Khovanov stable homotopy type for links in thickened surfaces
%Note that it has two kinds of degree, $q,\mathfrak h$.
\end{defn}

\h{\bf Note.}
In general, we need framings over moduli spaces when we construct a CW complex.
%See \cite[Definitions 3.18, 3.20, and 3.23]{LSk}.
However, in the case of Khovanov homotopy type,
we do not need to take care of framings.
See Proposition \ref{kya} and the comment below it
%\cite[(4) in the first part of section six]{LSk}, which is proved in the proof of \cite[Proposition 6.1]{LSk}. See  the three lines written above \cite[Definition 3.4]{LSs}: there is written, ``all such framings lead to the same Khovanov homotopy type \cite[Proposition 6.1]{LSk}''.  See \cite[Lemma 4.13, which is cited in the proof of Proposition 6.1]{LSk}.
%
%In short, in the case of Khovanov homotopy, all modulis are contractible so we do not need to check framings on them.
where it is pointed out that the construction is independent of the choice of framing.

%\vskip5mm  \begin{thm}\label{well}  Definition \ref{Khflow} is well-defined. \end{thm}

%\h{\bf Proof of Theorem \ref{well}.}  Because Theorem \ref{hanami} is true.\qed

\begin{defn}\label{tokei}
Let $L$ be a link diagram of a link $\mathcal L$ in a thickened higher genus surface.
We define
the {\it Khovanov-Lipshitz-Sarkar homotopy type}
(respectively,
the  {\it Khovanov-Lipshitz-Sarkar spectrum}
$\mathcal X_{Kh}(\mathcal L)
=\displaystyle\bigvee_{q,\mathfrak h}\mathcal X^{q,\mathfrak h}_{Kh}(\mathcal L)$)
for $\mathcal L$,
to be
the Khovanov-Lipshitz-Sarkar homotopy type
(respectively, the  Khovanov-Lipshitz-Sarkar spectrum
$\mathcal X_{Kh}(L)\\
=\displaystyle\bigvee_{q,\mathfrak h}\mathcal X^{q,\mathfrak h}_{Kh}(L)$)
for $L$,
defined in Definition \ref{Khflow}.
\end{defn}

We show examples in \S\ref{examp}.

\begin{thm}\label{korena}
Definition \ref{tokei} above is well-defined.
\end{thm}

\h{\bf Proof of Theorem \ref{korena}.}
Use almost the same method as that in \cite[\S 5]{LSk}.
%we can prove Theorem \ref{korena}.
Note, in particular,
that
\cite[the comment
between
the end of the proof Proposition 6.3,
and
 Proposition 6.4]{LSk}
also holds
in our case:
 use
the rule of partial order in Definition \ref{2.10}
instead of that in \cite[Definition 2.10]{LSk}.
\qed\\

\bb
\h{\bf Proof of Main theorem \ref{main}.}
In \cite{Seed}, it is proved that there is a pair of  knots,
$\mathcal J$ and $\mathcal J'$,
in $S^3$ (respectively, $D^3$)
with the following properties:
The second Steenrod square of $\mathcal J$ and that of $\mathcal {J'}$ are different.
The Khovanov homotopy type of $\mathcal J$ and that of $\mathcal {J'}$ are different.
The Khovanov homology of $\mathcal J$ and that of $\mathcal J'$ are the same.

Take the knots in $D^3$ in thickened higher genus surfaces.
Therefore  Main theorem \ref{main} holds.

Furthermore we have the following.
%in the case of links in thickened higher genus surfaces.
%
%Take a knot diagram $K=S^1\x{\ast}\x\{0\}$ in $S^1\x S^1\x[-1,1]$.
Let $F$ be a surface. Let $C$ be a  circle in $F$ which represents
a
%basis
nontrivial
element of
$H_1(F;\Z)$.
Regard $C$ as a knot in $F\x[-1,1]$.
Take $\mathcal J$ and $\mathcal J'$ in a 3-ball $B$ embedded in $F\x[-1,1]$.
Assume that $C\cap B=\emptyset$.
Make a disjoint 2-component link
which is made from $C$ and
$\mathcal {J}$ (respectively, $\mathcal {J'}$).
By the above result in \cite{Seed}, %and \cite[Lemma 2.14]{LSk},
these two links have different Steenrod squares and the same Khovanov homology.
\qed

\begin{note}\label{kikon}
Our case is different from Lipshitz and Sarkar's case \cite{LSk}
in that
the circles in $D$ may be non-contractible circles, and arcs in $D$ may be scs arcs.
However, as we saw above,
by the property on the sign convention
in Definition \ref{korekos},
%explained in
%Note
%\S\ref{konnyaku},
%\S\ref{KMNhomology},
the Khovanov chain complexes for links in thickened surfaces
have similar theorems in \cite[section 5]{LSk}.
More precisely we have Theorem \ref{tsukiyo}.
%See also Notes \ref{konki} and \ref{kikon}.

%We have proved Theorem \ref{tsukiyo}.However, i
In the case of virtual links in \cite{DKK, KauffmanOgasasq, Man},
%this theorem
Theorem \ref{tsukiyo}
is not true.  as explained in Note \ref{konki}.
%See Figure \ref{eK4}.
See also \cite{KauffmanOgasasq}:
The reason of this difference is the sign convention of the differential.
%as we explain in Note \S\ref{konnyaku}.
%This is a difference between the case of the Khovanov homology for links in thickened surfaces in this paper and that for virtual links in \cite{DKK, KauffmanOgasasq, Man} as explained inNote \S\ref{konnyaku}.
See also
\S\ref{itokon}.
\end{note}

\section{\bf Examples of Khovanov homotopy type}\label{examp}

We show examples of Khovanov homotopy type.

\subsection{links whose link diagrams have one crossing}\label{ex1}
\begin{figure}%[h]
\centering\includegraphics[width=0.4\textwidth]{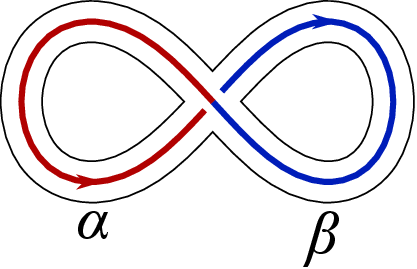}
\caption{A neighborhood of a link diagram $U$ with one crossing.}\label{fig:unknot_diagram}
\end{figure}

\hskip10mm
Consider an oriented link diagram $U$ with one crossing, see Fig.~\ref{fig:unknot_diagram}.

\begin{figure}%[h]
\centering\includegraphics[width=0.3\textwidth]{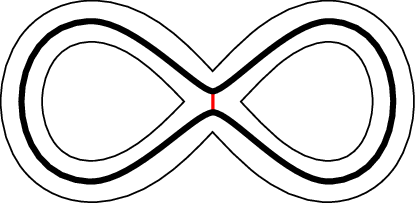}\qquad
\includegraphics[width=0.3\textwidth]{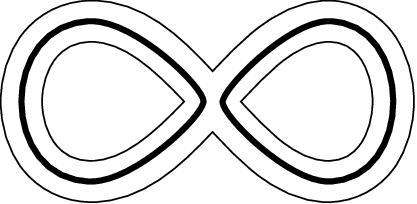}
\caption{Resolution configurations $D_U(0)$ (left) and $D_U(1)$ (right) of the link diagram $U$.}\label{fig:unknot_resolutions}
\end{figure}

The diagram has two resolution configurations, see Fig.~\ref{fig:unknot_resolutions}.

The resolution configurations correspond to six labeled resolution configurations, see Table~\ref{tab:unkont_labeled_resoutions}. The homotopical grading is expressed using the homotopical classes of the left and the right loops of the diagram, see Fig.~\ref{fig:unknot_diagram}.

\begin{table}
\caption{Labeled resolution configurations of the diagram $U$.}\label{tab:unkont_labeled_resoutions}
\centering
\begin{tabular}{ccccc}
  \hline
  Name & Generator & $\mathrm{gr}_h$ & $\mathrm{gr}_q$ & $\mathrm{gr}_{\mathfrak H}$\\
  \hline
  $\mathbf a$ & $(D_U(1),x_+x_+)$ & 0 & 1 & $[\alpha]+[\beta]$ \\
  $\mathbf b$ & $(D_U(1),x_+x_-)$ & 0 & -1 & $[\alpha]-[\beta]$ \\
  $\mathbf c$ & $(D_U(1),x_-x_+)$ & 0 & -1 & $-[\alpha]+[\beta]$ \\
  $\mathbf d$ & $(D_U(1),x_-x_-)$ & 0 & -3 & $-[\alpha]-[\beta]$ \\
  $\mathbf x$ & $(D_U(0),x_+)$ & -1 & -1 & $[\alpha\beta]$ \\
  $\mathbf y$ & $(D_U(0),x_-)$ & -1 & -3 & $-[\alpha\beta]$ \\
  \hline
\end{tabular}
\end{table}

The partial order on the set of labeled resolution configurations depends on the layout of the knot $U$ in the surface. If the loops $\alpha$ and $\beta$ are contractible then the situation does not differ from the classical case, cf. \cite[section 9.1]{LSk}. We have ${\mathcal X}_{Kh}(U)=S^0\vee S^0$ here.

There are three homotopically nontrivial cases.

\begin{figure}%[h]
\centering\includegraphics[width=0.3\textwidth]{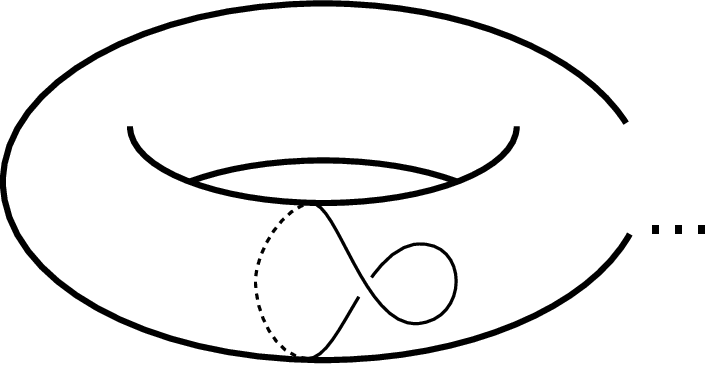}\qquad\quad
\includegraphics[width=0.3\textwidth]{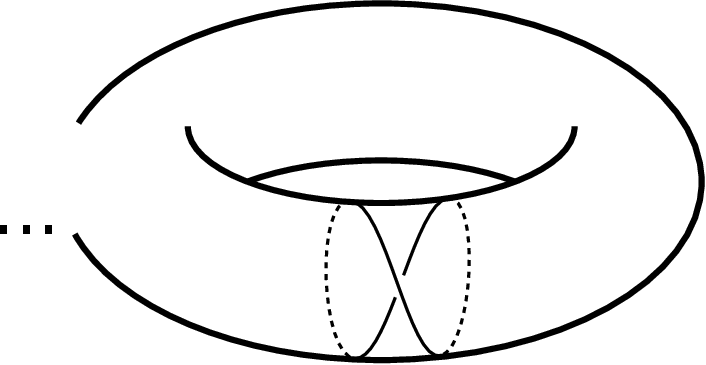}\\
\medskip\medskip
\includegraphics[width=0.6\textwidth]{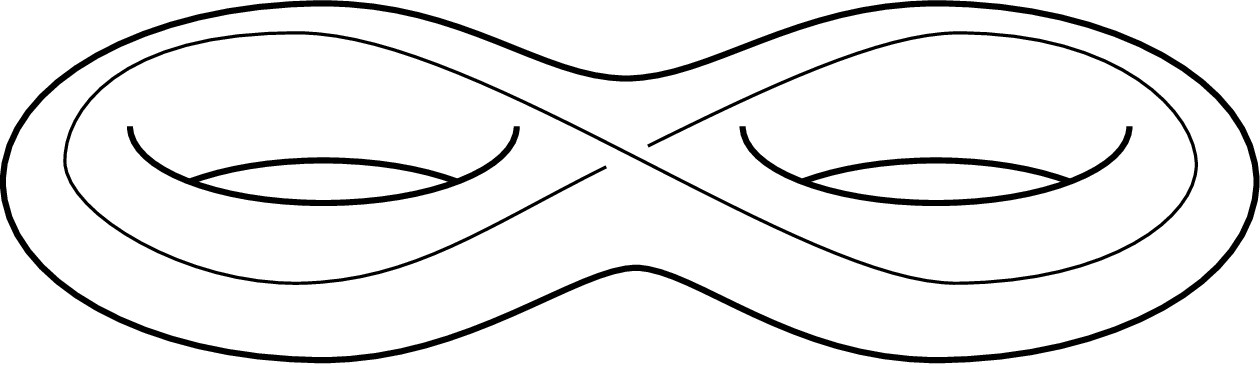}
\caption{Possible layouts of the knot $U$ in the surface.}\label{fig:unknot_diagram_cases}
\end{figure}

Case 1 (see Fig.~\ref{fig:unknot_diagram_cases}, upper left). Let one of the loops of $U$ (say $\alpha$) be non-contractible and the other be contractible. Then $[\beta]=0$ and $[\alpha\beta]=[\alpha]\ne 0$. Hence, we have $x\prec b$, $y\prec d$, and the other generators are incomparable because they have different quantum and homotopical gradings. We can treat ${\mathcal X}_{Kh}$ as a desuspension of the cell complex consisting of the basepoint $*$, two $0$-cells $x$ and $y$, and four $1$-cells $a,b,c,d$ where both ends of $a$ and $c$ are the basepoint $*$, $b$ is attached to $*$ and $x$, $d$ is attached to $*$ and $y$. Thus, ${\mathcal X}_{Kh}(U)= \Sigma^{-1}(S^1_a\vee S^1_c\vee D^1_b\vee D^1_d)=S^0\vee S^0$. In the splitting
$${\mathcal X}_{Kh}(U)=\bigvee_{q,{\mathfrak h}}{\mathcal X}^{q,{\mathfrak h}}_{Kh}(U)$$
we have ${\mathcal X}^{1,[\alpha]}_{Kh}(U)={\mathcal X}^{-1,-[\alpha]}_{Kh}(U)=S^0$, the other ${\mathcal X}^{q,{\mathfrak h}}_{Kh}(U)$ are trivial.

Case 2 (see Fig.~\ref{fig:unknot_diagram_cases}, upper left). Let the composition $\alpha\beta$ be contractible and $\alpha$ be non-contractible. Then $[\beta]=[\alpha^{-1}]=[\alpha]\ne 0$ and $[\alpha\beta]=0$. Hence, $x\prec b$, $x\prec c$, and the other generators are incomparable. Then we have
$${\mathcal X}^{1,2[\alpha]}_{Kh}(U)={\mathcal X}^{-1,0}_{Kh}(U)={\mathcal X}^{-3,-2[\alpha]}_{Kh}(U)=S^0,
{\mathcal X}^{-3,0}_{Kh}(U)=\Sigma^{-1}(S^0),$$
the other ${\mathcal X}^{q,{\mathfrak h}}_{Kh}(U)$ are trivial. Thus, ${\mathcal X}_{Kh}(U)=S^0\vee S^0\vee S^0\vee \Sigma^{-1}(S^0)$.

Case 3 (see Fig.~\ref{fig:unknot_diagram_cases}, lower). Let the $\alpha$, $\beta$ and $\alpha\beta$ be non-contractible. Then all the generators have different gradings and are incomparable. Hence, each generator yields a nontrivial space in the bouquet decomposition of ${\mathcal X}_{Kh}(U)=(\Sigma^{-1}S^0)^{\vee 2}\vee(S^0)^{\vee 4}$:
$$
{\mathcal X}^{1,[\alpha]+[\beta]}_{Kh}(U)={\mathcal X}^{-1,[\alpha]-[\beta]}_{Kh}(U)={\mathcal X}^{-1,-[\alpha]+[\beta]}_{Kh}(U)={\mathcal X}^{-3,-[\alpha]-[\beta]}_{Kh}(U)=S^0,$$

$$
{\mathcal X}^{-1,[\alpha\beta]}_{Kh}(U)={\mathcal X}^{-3,-[\alpha\beta]}_{Kh}(U)=\Sigma^{-1}(S^0).
$$

\subsection{links whose link diagrams have three crossings}\label{ex3}
\begin{figure}
\centering\includegraphics[width=0.4\textwidth]{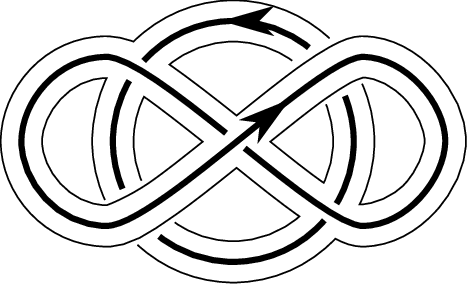}
\caption{A neighborhood of a link in a surface of genus 2. The complement to the neighborhood in the surface is a disk.}\label{fig:link3_11b}
\end{figure}

\hskip10mm
Consider the oriented link in a surface of genus 2 drawn in Fig.~\ref{fig:link3_11b}. The diagram has two positive and one negative crossings.

\begin{figure}
\centering\includegraphics[width=0.4\textwidth]{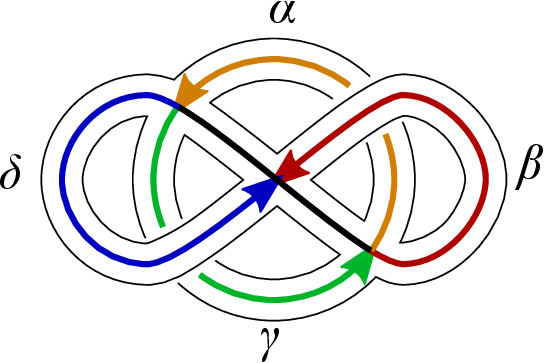}
\caption{The arcs of the link.}\label{fig:link3_arcs}
\end{figure}

The link as an embedded graph consists of a black segment which can be contracted to a point, and four arcs, that generate the fundamental group of the surface, see Fig.~\ref{fig:link3_arcs}. We denote these arcs as $\alpha,\beta,\gamma,\delta$.

\begin{figure}
\centering\includegraphics[width=0.9\textwidth]{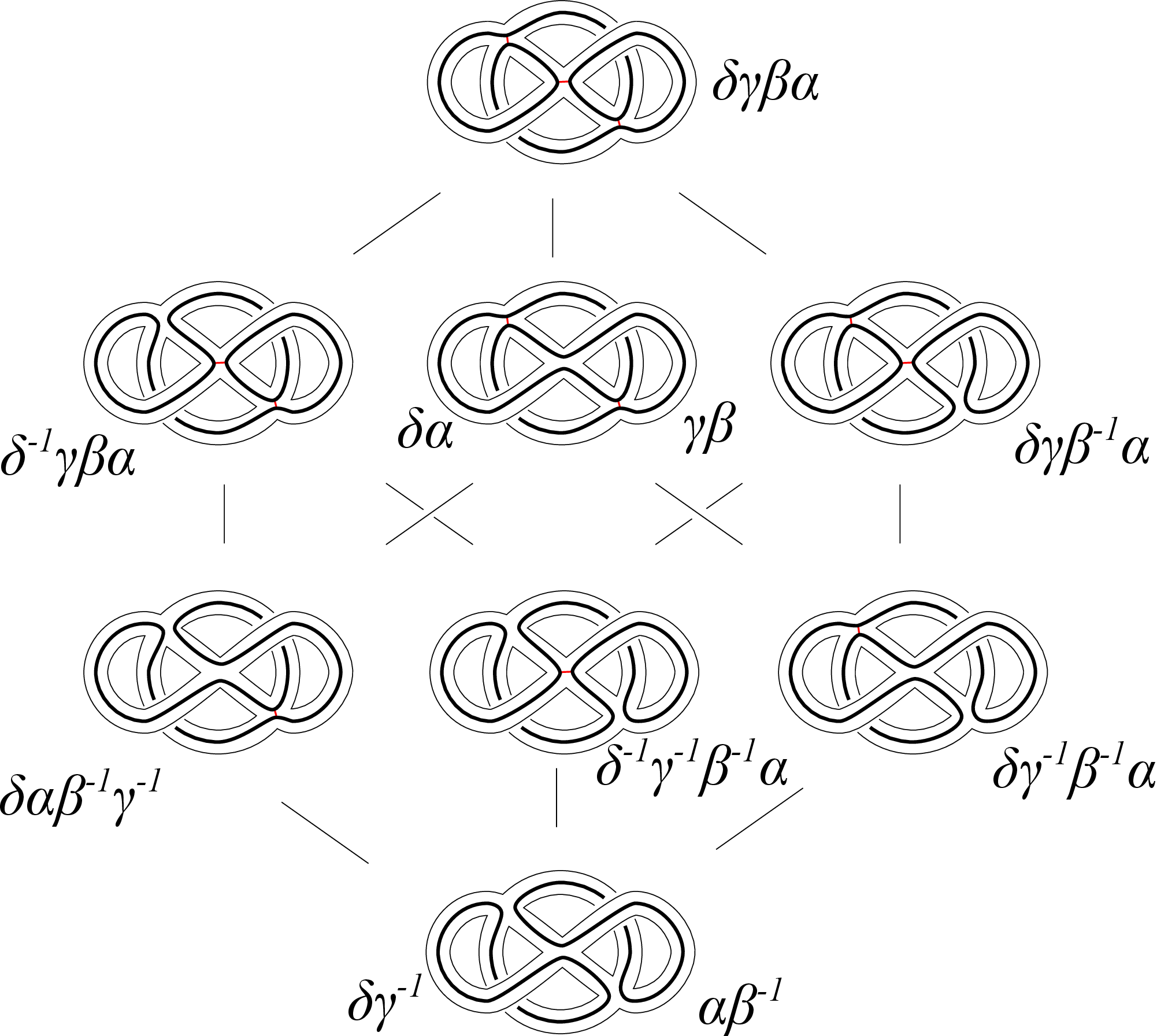}
\caption{The resolution cube of the link diagram. The circles are marked with their homotopy type.}\label{fig:link3_resolution_cube}
\end{figure}

Consider the resolution cube of the diagram (Fig.~\ref{fig:link3_resolution_cube}). The homotopy classes of the circles in the resolution configurations are all different. Hence, all the labeled resolution configurations are incomparable by the partial order, and there are no nontrivial decorated resolution configurations. Then the differential in the homotopical Khovanov complex is zero, and the homotopical Khovanov homology coincides with the chain complex. The Khovanov--Lipshitz--Sarkar homotopy type is a bouquet of spheres corresponding to the labeled resolution configurations
$$
{\mathcal X}_{Kh}=(\Sigma^{-1}S^0)^{\vee 2}\vee (S^0)^{\vee 8}\vee(S^1)^{\vee 6}\vee (S^2)^{\vee4}.
$$
The dimension of the spheres is determined by the homological grading of the resolution configuration.

The example above is opposite to the classical case in some sense. For classical links, the homotopical grading does not matter because all circles are contractible. And in this case the homotopical grading brakes all connections between resolution configurations. Note that among the twelve surgeries in the resolution cube six are single circle surgeries.

\bb
\section{\bf Open questions}
%\begin{note}
\label{itokon}

\h
%As we recall in \S\ref{konnyaku},
Links in thickened surfaces are regarded as virtual links
(\cite{Kauffman1, Kauffman, Kauffmani}).
Recalling Note \ref{notevirtual},
it is natural to ask a question: Compare the strength of the following invariants.

\begin{enumerate}
\item[(1)]
 The Khovanov homology for virtual links,

\item[(2)]
 The second Steenrod square for virtual links
in conjunction with Khovanov homology  for virtual links

\item[(3)]
 The Khovanov %-Manturov-Nikonov
 homology for links in thickened surfaces

\item[(4)]
The second Steenrod square for links in thickened surfaces
in conjunction with Khovanov homology for links in thickened surfaces

\item[(5)]
 The Khovanov-Lipshitz-Sarkar homotopy type for links in thickened surfaces
\end{enumerate}

%In \S\ref{konnyaku} we asked what relations there are between (1) and (3).

 In \cite{KauffmanOgasasq}, it is proved that $(2)$ is stronger than $(1)$.

In this paper we prove the following:
$(4)$ is stronger than $(3)$ if the genus is greater than one.
$(5)$
is stronger than $(3)$ if the genus is greater than one.

In the case of the thickened torus, we shall deal with it in a separate paper \cite{zoku}.
We have not defined
Khovanov-Lipshitz-Sarkar stable homotopy type for virtual links.
 %nor that for links in the thickened torus.  have we defined a homotopy type for links in the thickened torus.
How does one introduce them?

%
%
%Can we combine the Steenrod square (and Khovanov homotopy theory) in the case of links in thickened surfaces and that of virtual links and  define a new Steenrod square (and Khovanov homotopy theory)?
%
%
Can we combine the Steenrod square (and Khovanov homotopy theory) in a single theory that would apply to both links in thickened surfaces and to virtual links?

%\end{note}

\np
\noindent
Louis H. Kauffman

\noindent
Department of Mathematics, Statistics and Computer Science

\noindent
University of Illinois at Chicago

\noindent
851 South Morgan Street

\noindent
Chicago, Illinois 60607-7045

\noindent
USA

\noindent
and

\noindent
Department of Mechanics and Mathematics

\noindent
Novosibirsk State University

\noindent
Novosibirsk

\noindent
Russia

\noindent
kauffman@uic.edu
\\

%\h Vassily Olegovich Manturov,  \\

\h Igor Mikhailovich Nikonov

\h Department of Mechanics and Mathematics

\h Lomonosov Moscow State University

\h  Leninskiye Gory, GSP-1

\h Moscow, 119991

\h Russia

\h nikonov@mech.math.msu.su
\\

\noindent
Eiji Ogasa

\noindent
Meijigakuin University, Computer Science

\noindent
Yokohama, Kanagawa, 244-8539

\noindent
Japan

\noindent
pqr100pqr100@yahoo.co.jp

\noindent
ogasa@mail1.meijigkakuin.ac.jp

%}\Large uke
\end{document}